\documentclass[a4paper, 11pt, oneside]{article}

\pdfoutput=1
\usepackage{amsmath}
\usepackage[runin]{abstract}
\usepackage{hyperref}
\usepackage{fancyvrb}
\usepackage{amssymb}
\usepackage{mathtools}
\usepackage{amsthm}
\usepackage{amsfonts}
\usepackage{dsfont}
\usepackage{fancybox}
\usepackage{xcolor}
\usepackage{tikz}

\newcommand\nnfootnote[1]{%
	\begin{NoHyper}
		\renewcommand\thefootnote{}\footnote{#1}%
		\addtocounter{footnote}{-1}%
	\end{NoHyper}
}

\newtheorem*{algorithm*}{Algorithm}

\theoremstyle{plain}
\newtheorem{Theorem}{Theorem}[section]
\newtheorem{Lemma}[Theorem]{Lemma}

\newtheorem{Definition}[Theorem]{Definition}

\theoremstyle{definition}
\newtheorem{Remark}[Theorem]{Remark}

\numberwithin{equation}{section}

\newcommand{\supp}{{\rm supp \, }}
\newcommand{\sgn}{{\rm sgn \, }}

\author{Marc Hovemann$^{\ast}$}
\title{Adaptive Bivariate Quarklet Tree Approximation via Anisotropic Tensor Quarklets}
\date{\today}

\begin{document}
\maketitle

\nnfootnote{
	\hspace{-21.5pt}
	This work partly has been supported by Deutsche Forschungsgemeinschaft (DFG), grant HO 7444/1-1 with project number 528343051.

\noindent	
${}^\ast$Friedrich-Schiller-Universität Jena, Fakultät für Mathematik und Informatik, Ernst-Abbe-Platz 2, 07743 Jena, Germany, Email: \texttt{marc.hovemann@uni-jena.de}, Phone number: +49-3641-9-46193.
}

\noindent
\textbf{Abstract.} This paper deals with near-best approximation of a given bivariate function using elements of quarkonial tensor frames. For that purpose we apply anisotropic tensor products of the univariate B-spline quarklets introduced around 2017 by Dahlke, Keding and Raasch. We introduce the concept of bivariate quarklet trees and develop an adaptive algorithm which allows for generalized  hp-approximation of a given bivariate function by selected frame elements. It is proved that this algorithm is near-best, which means that as long as some standard conditions concerning local errors are fulfilled  it provides an approximation with an error close to that one of the best possible quarklet tree approximation. For this algorithm the complexity is investigated. Moreover, we use our techniques to approximate a bivariate test function with inverse-exponential rates of convergence. It can be expected that the results presented in this paper serve as important building block for the design of adaptive wavelet-hp-methods for solving PDEs in the bivariate setting with very good convergence properties. 
 
\vspace{0,3 cm}


\vspace{0,3 cm}

\section{Introduction}

Many problems in natural sciences, economics and public finance can be described by partial differential equations. Often a closed form of the unknown solution is not known, and hence numerical schemes in order to find a good approximation for it are required. Thereto a very popular approach is the finite element method (FEM). The well-known $h$-FEM is based on a space refinement of the domain of interest. Alternatively, when it comes to the so-called $p$-method, the polynomial degrees of the ansatz functions are increased. It is also possible to combine both methods in order to obtain $hp-$FEM techniques. When dealing with large-scale problems often it is advantageous to deploy adaptive strategies to increase the overall efficiency. The goal is to obtain a satisfactory approximation after a tolerable number of calculation steps. In particular for adaptive $h$-FEM there exists a huge amount of literature. Let us refer to \cite{bib:Cia02}, \cite{bib:Hac10}, \cite{NSV09}, \cite{bib:Sch98} and \cite{Verf} at least. In recent years the convergence analysis of adaptive $p$- and $hp$-methods attracted a lot of attention. It turned out that these schemes converge very fast and in many cases even show exponential convergence. However, concerning theoretical analysis and rigorous convergence proofs only a few results have been derived recently. Some state of the art results concerning the convergence of adaptive $hp$-strategies are \cite{bib:BPS13}, \cite{bib:BD11}, \cite{bib:DV20}, \cite{bib:DH07} and \cite{bib:CNSV14}, \cite{bib:CNSV17}, whereby the latter  also contain optimality results. 

Another approach is the use of wavelets. Wavelets have very strong analytical properties that can be utilized to attain adaptive methods that converge with the optimal order of the best $N$-term wavelet approximation. In connection with that let us refer to \cite{bib:CDD01} and \cite{bib:Ste09}. In the main adaptive wavelet schemes are space refinement methods and hence can be classified as $h$-methods. Then the natural question arises how $hp$-versions of adaptive wavelet schemes can be designed. At this juncture the approach of using \emph{quarklets} comes into play. Quarklets are polynomially enriched wavelets that have been introduced in the last decade in the pioneering paper  \cite{DaKRaa}. Univariate quarklets are constructed out of biorthogonal compactly supported Cohen-Daubechies-Feauveau spline wavelets, whereby the primal generator is a cardinal B-spline. The theory of these biorthogonal wavelets can be found in Section 6.A in \cite{CoDau}. In principle univariate quarklets are linear combinations of translated cardinal B-splines multiplied with some monomials. A precise definition can be found in Definition \ref{def_quarklet} below. The theoretical properties of univariate quarklets have been studied in detail in \cite{DaFKRaa,DaHoRaVo,DaKRaa,DaRaaS,HoDa,HoKoRaVo,SiDiss} and \cite{VoDiss}. Moreover, it turned out that quarklets can be used to design schemes that resemble $hp$-versions of adaptive wavelet methods. In \cite{DaHoRaVo} univariate quarklets have been used to approximate functions $ f \in L_{2}((0,1))  $ in a very efficient way. For that purpose an adaptive algorithm called {\bf NEARBEST\textunderscore TREE} is provided which allows for both space refinement and polynomial enrichment. A very important role for the theory developed in \cite{DaHoRaVo} plays the newly introduced concept of univariate quarklet trees. Using a proof technique developed by Binev in  \cite{Bin18} it is shown that the algorithm {\bf NEARBEST\textunderscore TREE} is near-best which means that it delivers an approximation with an error close  to the best tree approximation error for a given cardinality. Moreover, in \cite{DaHoRaVo} several numerical experiments are presented which show, that adaptive univariate quarklet tree approximation can be applied to approximate certain functions with inverse-exponential convergence rates. 

Recently in \cite{DaFKRaa} and \cite{Ho_bivqua} also bivariate and even multivariate quarklets have been introduced. They have been constructed out of univariate quarklets using anisotropic tensor products. A precise definition also is recalled in Section \ref{subsec_biv_tensorwavelets} below. The main goal of this paper is to design an adaptive algorithm {\bf BIVARIATE\textunderscore NEARBEST\textunderscore TREE} which allows to approximate bivariate functions via bivariate tensor quarklets in a very efficient way. In the long run this algorithm could serve as important building block for the design of adaptive quarklet-hp-methods for solving PDEs in the bivariate setting with very good convergence properties. Indeed, for univariate linear elliptic variational problems an optimal quarklet Galerkin scheme using univariate quarklet tree approximation already has been developed  successfully, see Chapter 5 in \cite{VoDiss}. 

The foundation of the theory presented in this paper is the concept of bivariate quarklet trees which is introduced in Section \ref{subsec_biv_qua_tre} below. Although univariate quarklet trees already have been defined in \cite{DaHoRaVo}, the specification of bivariate quarklet trees is a delicate job, since in the multivariate setting several new phenomena show up. Some of them, which will be discussed in detail throughout this paper, are the following: 

\begin{itemize}
\item[(i)] First of all in the bivariate setting we have to deal with reference rectangles instead of the univariate reference intervals. These reference rectangles can be highly anisotropic, since it is possible to work with different refinement levels in the two Cartesian directions. At the first glance this circumstance seems to make things more difficult, but in the long run it should be possible to design algorithms with better convergence properties when we also use anisotropic tensor quarklets for approximation. A precise definition of reference rectangles can be found in Section \ref{subsec_ref_rectangle1}.

\item[(ii)] A second challenging task is the determination of an unique parent-child relation between reference rectangles (or equivalently between wavelet indices) of neighboring refinement levels. For a given reference rectangle there are different refinement options due to the two Cartesian directions. Thereby, in a single refinement step refinement is carried out in only one direction, and as a consequence the reference rectangle is divided into two smaller child rectangles. However, very often we enable the incidence of a third child, which does not stand for a current refinement, but describes the possibility to carry out a refinement step in the other direction later on. All in all we only allow refinement steps, where two or three children are showing up. When we think on approximation in the context of large-scale problems, this strategy might pave the way to approximation schemes with dimension-independent rates of convergence. Much more information concerning parent-child relations in the context of adaptive bivariate quarklet tree approximation is given in Section \ref{subsec_wav_tree}.

\item[(iii)] Another delicate issue is the definition of bivariate wavelet and quarklet trees itself. The tree structure should be arranged in such a way that it can be used to design an adaptive algorithm to approximate bivariate functions which is near-best and has very good convergence properties. Below we utilize a proof technique stemming from \cite{DaHoRaVo} which only works if each wavelet node has an unique parent. Since in the bivariate setting this condition is not fulfilled automatically, in Section \ref{subsec_wav_tree} we introduce an additional refinement rule called unique parent condition. It restricts the number of possible refinement options such that each wavelet node has a unique parent. Thereby it is not too restrictive since still all possible refinement samples can be reached. When it comes to the definition of bivariate quarklet trees, another obstacle shows up. When we grow a quarklet tree also the inner nodes remain as active contributors to the approximation. Consequently, we have to assign a polynomial degree not only to the leaves of the tree, but also to all inner nodes. This problem already has been observed in \cite{DaHoRaVo} for the univariate setting. However, since bivariate quarklet trees obviously have a more complicated structure than their univariate counterparts, this issue becomes even more difficult in the bivariate setting. One possible solution is presented in Section \ref{subsec_biv_qua_tre}, where sets $  \Upsilon( \cdot )   $ are introduced, which provide a partition of the whole bivariate wavelet tree.

\end{itemize}

An important key result of this paper is the adaptive algorithm {\bf BIVARIATE\textunderscore NEARBEST\textunderscore TREE} which is presented in Section \ref{subsec_alg1}. For a given bivariate function it provides a bivariate wavelet tree. Using a trimming routine it can be transformed into a bivariate quarklet tree which delivers an approximation in terms of bivariate tensor quarklets for the input function. The algorithm {\bf BIVARIATE\textunderscore NEARBEST\textunderscore TREE} is based on its univariate forerunner {\bf NEARBEST\textunderscore TREE} presented in \cite{DaHoRaVo}. However, it requires some additional computation in order to decide in which Cartesian direction a refinement should be carried out in each step. The algorithm {\bf BIVARIATE\textunderscore NEARBEST\textunderscore TREE} can be seen as $hp$-method since it allows for both space refinement and polynomial enrichment, whereby the decision what to do next is found in an adaptive way. Depending on the input function the algorithm is able to produce highly anisotropic refinement meshes, see also Figures 3.3 and 3.4 in \cite{BuGrie}. This is advantageous if it comes to the approximation of functions with anisotropic singularities, see Section \ref{Sec_aniso_example}, for instance. Furthermore, in the long run when it comes to the approximation of multivariate functions more benefits of anisotropic refinement are expected. Indeed, we can hope to generalize the results of this paper to the multivariate setting in order to obtain adaptive quarklet approximation methods with dimension-independent convergence rates, since the tensor product quarklets can be interpreted as wavelet versions of sparse grids including polynomial enrichment.

In Lemma \ref{lem_complexity} we investigate the complexity of the algorithm {\bf BIVARIATE\textunderscore NEARBEST\textunderscore TREE}. Furthermore, in our main result Theorem \ref{theorem:1} we show that the approximations provided by the algorithm {\bf BIVARIATE\textunderscore NEARBEST\textunderscore TREE} are near-best. This means that we obtain approximations with an approximation error close to the error of the best possible bivariate quarklet tree approximation. To see this we use a proof technique already applied in \cite{DaHoRaVo} to deal with the univariate setting and modify it in order to treat the bivariate case. 

Both our algorithm {\bf BIVARIATE\textunderscore NEARBEST\textunderscore TREE} and Theorem \ref{theorem:1} are formulated in a very general way. Consequently, they can be utilized to approximate a very broad class of functions. However, in Section \ref{sec_practice} we explain how our approach can be used to approximate functions $  f \in L_{2}((0,1)^2)   $. For that purpose some local error functionals have to be defined according to the $  L_{2}((0,1)^2)$-setting. Finally, Theorem \ref{lem_L2_bestappr} shows that our algorithm {\bf BIVARIATE\textunderscore NEARBEST\textunderscore TREE} also is near-best in the case of $  L_{2}((0,1)^2)$-approximation. 

This paper is organized in the following way. In Section \ref{sec_def_quarks11} at first the concept of univariate quarklets on the real line is recalled. In addition we explain how they must be modified in order to obtain boundary adapted quarklets on intervals such as $ (0,1)  $. Moreover, we define bivariate tensor quarklets by using univariate quarklets and tensor product methods.
In Section \ref{sec_biv_tre_main111} we introduce bivariate quarklet trees. To prepare this at first we recall the concept of reference rectangles and explain what (enhanced) bivariate wavelet indices are. Once we have determined which refinement strategies are allowed, we can define bivariate wavelet trees in Definition \ref{def_wav_tre}. After that also (enhanced) bivariate quarklet indices are introduced. They show up in Definition \ref{def:quarklet_tree} when it comes to the specification of bivariate quarklet trees. 
The core part of this paper is Section \ref{sec:adap_ref}. Here the central algorithm {\bf BIVARIATE\textunderscore NEARBEST\textunderscore TREE} is provided. As a preparation for this at first some local and global error functionals are introduced. Moreover, a trimming routine called {\bf BIVARIATE\textunderscore TRIM} is developed. It transforms a bivariate wavelet tree produced by {\bf BIVARIATE\textunderscore NEARBEST\textunderscore TREE} into a bivariate quarklet tree with very good approximation properties. Finally we prove our main Theorem \ref{theorem:1} which verifies that the approximations found by the algorithm {\bf BIVARIATE\textunderscore NEARBEST\textunderscore TREE} are near-best indeed. Section \ref{sec_practice} is devoted to the special case of $ L_{2}((0,1)^2)$-approximation using bivariate tensor quarklets. In order to run the algorithm {\bf BIVARIATE\textunderscore NEARBEST\textunderscore TREE} in this setting we have to define some local errors in a suitable way which is explained in Definition \ref{def_locerr_L2_new}. Moreover, by proving Theorem \ref{lem_L2_bestappr} we see that {\bf BIVARIATE\textunderscore NEARBEST\textunderscore TREE} also is near-best in the case of $  L_{2}((0,1)^2)$-approximation. Finally, in Section \ref{Sec_aniso_example} we present a bivariate test function which can be approximated by using bivariate quarklet trees, whereby inverse-exponential rates of convergence are achieved.

\section{Quarks and Quarklets}\label{sec_def_quarks11}

In this paper we use bivariate tensor quarklets to approximate functions $f \in L_{2}((0,1)^2)$. To construct such quarklets, we have to carry out several substeps. At first we deal with univariate quarklets, defined either on $\mathbb{R}$ or on bounded intervals such as $(0,1)$.

\subsection{B-Splines, Quarks and Quarklets on the Real Line}\label{sec_inner_quark}

In the following section we recall the definition of univariate quarklets for the shift-invariant setting on $\mathbb{R}$. For that purpose we follow \cite{DaKRaa}. In a first step we repeat the definition of cardinal B-splines. The first order cardinal B-spline $  N_{1}  $ is just the characteristic function of the interval $ [0,1)  $, namely $ N_{1} := \chi_{[0,1)} $. Higher order cardinal B-splines of order $ m \in \mathbb{N}  $ with $ m \geq 2  $ are defined by induction using the convolution $ \ast $. So we have 
\begin{align*}
N_{m} := N_{m - 1} \ast N_{1} = \int_{0}^{1} N_{m-1}( \cdot - t ) dt.
\end{align*}
The cardinal B-splines possess some very nice properties, see for example Chapter 5.2 in \cite{DeLo} and \cite{Ch1}. In what follows for fixed $ m \in \mathbb{N} $ we will work with the symmetrized cardinal B-spline $ \varphi(x)  :=   N_{m}  ( x + \lfloor \frac{m}{2} \rfloor   )    $. We observe $ \supp \varphi = [ - \lfloor \frac{m}{2} \rfloor    ,  \lceil \frac{m}{2} \rceil    ]   $. The symmetrized cardinal B-spline shows up in the following definition where we explain the so-called quarks.

\begin{Definition}\label{Bquark}
Let $ m \in \mathbb{N} $ and $ p \in \mathbb{N}_{0}  $. Then the p-th cardinal B-spline quark $ \varphi_{p}  $ is defined by
\begin{equation}
\varphi_{p}(x)  := \Big ( \frac{x}{\lceil \frac{m}{2} \rceil } \Big )^{p}  N_{m} \Big ( x + \lfloor \frac{m}{2} \rfloor  \Big ) .
\end{equation}
\end{Definition}
\noindent
The quarks are very important in order to define the quarklets. Their properties have been studied in \cite{DaKRaa}. It is shown in \cite{CoDau} by Cohen, Daubechies and Feauveau that for a given $ \tilde{m} \in \mathbb{N}  $ with $ \tilde{m} \geq m    $ and $  m + \tilde{m} \in 2 \mathbb{N}   $ there exists a compactly supported biorthogonal spline wavelet $ \psi  $ (sometimes also called CDF-wavelet) with
\begin{equation}\label{def_CDF_wav}
\psi = \sum_{k \in \mathbb{Z}} b_{k}  \varphi ( 2 \cdot - k   )
\end{equation}
with expansion coefficients $ b_{k} \in \mathbb{R}   $. Only finitely many of them are not zero. Moreover $ \psi  $ has $ \tilde{m}   $ vanishing moments and the  system
\begin{align*}
\Big \{ \varphi (  \cdot - k )  \ : \ k \in \mathbb{Z}  \Big  \} \cup \Big \{ 2^{\frac{j}{2}} \psi (2^{j} \cdot - k) \ : \ j \in \mathbb{N}_{0} \ , \ k \in \mathbb{Z}  \Big \}
\end{align*}
is a Riesz basis for $ L_{2}(\mathbb{R})   $. To construct such a $ \psi   $ we have to work with a compactly supported dual generator $ \tilde{\varphi}   $ associated to the primal generator $ \varphi $ that fulfills
\begin{equation}\label{biorto1}
\left \langle  \varphi , \tilde{\varphi} (\cdot - k)   \right\rangle_{L_{2}(\mathbb{R})} = \delta_{0,k} , \qquad k \in \mathbb{Z}.
\end{equation}
Connected with that there is another compactly supported biorthogonal wavelet $ \tilde{\psi} \in L_{2}(\mathbb{R})  $  
\begin{equation}\label{def_biort_wav1}
\tilde{\psi} = \sum_{k \in \mathbb{Z}}  \tilde{b}_{k} \tilde{\varphi} ( 2 \cdot - k   ).
\end{equation}
Here only finitely many of the $ \tilde{b}_{k} \in \mathbb{R}   $  are not zero. Moreover, $ \tilde{\psi}  $ has $ m \in \mathbb{N}  $ vanishing moments and the system
\begin{align*}
\Big \{ \tilde{\varphi} (  \cdot - k )  \ : \ k \in \mathbb{Z}  \Big  \} \cup \Big \{ 2^{\frac{j}{2}} \tilde{\psi} (2^{j} \cdot - k) \ : \ j \in \mathbb{N}_{0} \ , \ k \in \mathbb{Z}  \Big \}
\end{align*}
is a Riesz basis for $ L_{2}(\mathbb{R})   $. For $ j \in \mathbb{N}_{0}   $ and $ k \in \mathbb{Z}   $ let us write
\begin{equation}\label{def_biort_wav2}
\psi_{j,k} = 2^{\frac{j}{2}} \psi ( 2^{j} \cdot - k ) \qquad \qquad \mbox{and} \qquad \qquad \tilde{\psi}_{j,k} = 2^{\frac{j}{2}} \tilde{\psi} ( 2^{j} \cdot - k ) .
\end{equation}
Moreover, for $ k \in \mathbb{Z}   $ we put $ \psi_{-1,k} = \varphi (  \cdot - k )   $ and $ \tilde{\psi}_{-1,k} = \tilde{\varphi} (  \cdot - k )   $. Then we observe
\begin{equation}\label{biorto2}
\langle \psi_{j,k} , \tilde{\psi}_{j',k'}     \rangle_{L_{2}(\mathbb{R})} = \delta_{j , j'} \delta_{k , k'} , \qquad j, j' \in \mathbb{N}_{0}, \quad k, k' \in \mathbb{Z} .
\end{equation}
For each $ f \in L_{2}(\mathbb{R})  $ we find
\begin{align}
f & = \sum_{k \in \mathbb{Z}} \langle f , \tilde{\psi}_{-1,k}     \rangle_{L_{2}(\mathbb{R})} \psi_{-1,k}  + \sum_{j \in \mathbb{N}_{0} ,k \in \mathbb{Z}} \langle f , \tilde{\psi}_{j,k}     \rangle_{L_{2}(\mathbb{R})} \psi_{j,k} \nonumber \\
& = \sum_{k \in \mathbb{Z}} \langle f , \psi_{-1,k}     \rangle_{L_{2}(\mathbb{R})} \tilde{\psi}_{-1,k}  +  \sum_{j \in \mathbb{N}_{0}, k \in \mathbb{Z}} \langle f , \psi_{j,k}     \rangle_{L_{2}(\mathbb{R})} \tilde{\psi}_{j,k}    \label{def_biort_wav3}
\end{align}
with convergence in $ L_{2}(\mathbb{R})  $. For details and proofs concerning the above construction we refer to \cite{CoDau}, see especially Section 6.A. Now we can use the CDF-wavelets $ \psi $ to define the quarklets.  
\begin{Definition}\label{def_quarklet}
Let $ p \in \mathbb{N}_{0}  $. Then the p-th quarklet $ \psi_{p} $ is defined by 
\begin{equation}
\psi_{p} := \sum_{k \in \mathbb{Z}} b_{k} \varphi_{p}(2 \cdot - k).
\end{equation}
Here the $ b_{k}  $ are the same as in \eqref{def_CDF_wav}. Furthermore, for $ j \in \mathbb{N}_{0}   $ and $ k \in \mathbb{Z}  $ we write 
\begin{equation}
\psi_{p,j,k} := 2^{\frac{j}{2}} \psi_{p}(2^{j} \cdot - k) \qquad \qquad \mbox{and} \qquad \qquad \psi_{p,-1,k} := \varphi_{p}( \cdot - k) .
\end{equation}
\end{Definition}

\begin{Remark}\label{rem_hist_of_qua}
The univariate quarklets given in Definition \ref{def_quarklet} have been introduced around 2017 in \cite{DaKRaa}. Later their properties have been studied in detail in \cite{DaFKRaa}, \cite{DaRaaS} and \cite{HoKoRaVo}. A systematic and comprehensive treatise can be found in \cite{SiDiss}. When defining the quarklets the main focus is on numerical applications. They are specially tailored for adaptive approximation of functions and the numerical treatment of PDEs with very good convergence properties. For that purpose let us refer to \cite{DaHoRaVo} and \cite{DaFKRaa}. For the definition of the quarklets B-spline wavelets are used for the following reasons. B-splines possess optimal smoothness properties compared to their support size. Moreover, explicit formulas exist which make point evaluations quite simple. This issue is important for the construction of suitable quadrature formulas, that are necessary for any numerical scheme for the treatment of PDEs. Of course, in principle quarkonial decompositions also can be provided using other wavelets such as orthonormal Daubechies wavelets. However, Daubechies wavelets are not symmetric which sometimes is disadvantageous. Moreover, for these wavelets no explicit formulas exist which makes point evaluations much more difficult. Using B-spline pre-wavelets would also be a possible choice when constructing quarkonial decompositions. However, the biorthogonal approach we used in Definition \ref{def_quarklet} has the advantage that the lenghts of all filters involved in the associated decomposition and reconstruction schemes are finite, which is usually not the case in the pre-wavelet setting. Of course, due to the polynomial enrichment, the quarklet dictionary is highly redundant. For the construction of adaptive wavelet $hp$-methods, this fact cannot be avoided. At the first glance, this might look as a disadvantage, but it seems to be clear that this is not the case. So the long-term goal is the development of adaptive numerical schemes based on quarklets. The art of adaptivity is to find a sparse expansion of an unknown object, namely the solution of a PDE. Now if we work with a very rich dictionary, then the chance to find such a sparse expansion is much higher compared to the basis case where the expansion is unique. From this point of view, redundancy is very helpful. Indeed, in \cite{DaHoRaVo} for the univariate setting some numerical experiments showed, that our quarklets can be used for adaptive $hp$-tree approximation of functions $f \in L_{2}((0,1))$ with inverse-exponential convergence rates. Moreover, a rigorous proof that certain model singularities showing up in the solution theory of elliptic PDEs can be approximated via quarklets with inverse-exponential rates can be found in \cite{DaRaaS}. The present paper can be seen as a continuation of \cite{DaHoRaVo} to the bivariate setting. We will see that bivariate tensor quarklets can be used to approximate functions $f \in L_{2}((0,1)^2)$ with very good convergence properties.  
\end{Remark}

\subsection{Boundary Adapted Quarks and Quarklets on the Interval}\label{Subsec_bound_qua}

When we deal with univariate functions defined on bounded intervals such as $I := (0,1) \subset \mathbb{R}$ we require special boundary adapted quarks and quarklets. Their construction is explained in \cite{SiDiss} and \cite{DaFKRaa} and will be summarized in the following section. The foundation of this construction is given by a wavelet basis designed by Primbs, see \cite{Pri1}. In a first step we recall the definition of the so-called Schoenberg B-splines. Again let $ m, \tilde{m} \in \mathbb{N}_{0}   $ with $ \tilde{m} \geq m \geq 2   $ and $ m + \tilde{m} \in 2 \mathbb{N}   $. Let $ j_{0} \in \mathbb{N}  $ be a fixed number that depends on $m$ and $\tilde{m}$ and is sufficiently large, see Chapter 4.4 in \cite{Pri1} for further explanations. For $ j \in \mathbb{N} $ with $ j \geq j_{0}  $ let $ \Delta_{j} := \{ -m+1, \ldots , 2^{j} - 1   \} $. We define the knots 
\[
t^{j}_{k}:= \left\{ \begin{array}{lll}
0    & \quad & \mbox{for} \qquad k = -m+1, \ldots , 0 ;
\\  
2^{-j} k & \quad & \mbox{for}\qquad k = 1 , \ldots , 2^{j} - 1 ;
\\
1 & \quad & \mbox{for}\qquad k = 2^{j} , \ldots , 2^{j} + m - 1 .
\\
\end{array}
\right.
\]
Now the Schoenberg B-splines $ B^{m}_{j,k}  $ are defined by
\begin{equation}\label{eq_def_SchoenbergBs}
B^{m}_{j,k}(x) := ( t^{j}_{k+m} - t^{j}_{k} ) ( \cdot - x )^{m-1}_{+} [ t^{j}_{k}, \ldots , t^{j}_{k+m} ] ,  \qquad k \in \Delta_{j} , x \in I .
\end{equation}
Here the symbol $ ( \cdot - x )^{m-1}_{+} [ t^{j}_{k}, \ldots , t^{j}_{k+m} ]    $ stands for the $m-$th divided difference of the function $  ( \cdot - x )^{m-1}_{+}   $. The generating functions of the Primbs basis are 
\begin{equation}\label{eq_qua_gen_jjj}
\varphi_{j,k} := 2^{\frac{j}{2}} B^{m}_{j,k}, \qquad k \in \Delta_{j}.
\end{equation}
The Schoenberg B-splines are generalizations of the cardinal B-splines $ N_{m} $ and have some useful properties, see for example \cite{Pri1}. Recall that the Primbs basis is a biorthogonal wavelet basis. Therefore a dual multiresolution analysis with dual generators $\tilde{\varphi}_{j,k}$ is necessary for the construction. If the generators are represented as column vectors $\Phi_{j} := \{ \varphi_{j,k} : k \in  \Delta_{j}  \} $ and $\tilde{\Phi}_{j} := \{ \tilde{\varphi}_{j,k'} : k' \in  \Delta_{j}  \} $, they fulfill the duality relation
\begin{align*}
\langle \Phi_{j} , \tilde{\Phi}_{j} \rangle := ( \langle \varphi_{j,k} , \tilde{\varphi}_{j,k'} \rangle_{L_{2}(I)})_{k , k' \in \Delta_{j}  } = Id_{|\Delta_{j}  |}.
\end{align*}
For the construction of the Primbs wavelets the following index set is defined:
\[
\nabla_{j}:= \left\{ \begin{array}{lll}
\{ 0, 1, \ldots , 2^{j} - 1 \}   & \quad & \mbox{for} \qquad j \geq j_{0} ;
\\  
\Delta_{j_{0}} & \quad & \mbox{for}\qquad j = j_{0} - 1 .
\\
\end{array}
\right.
\]
To construct the Primbs wavelets suitable matrices $M_{j,1}^{b}$, $\tilde{M}_{j,1}^{b}$ are defined, that contain the two-scale coefficients of the wavelet column vectors $\Psi = \{ \psi_{j,k}^{b} : k \in  \nabla_{j}  \}$. We have
\begin{equation}\label{eq_primbs_wav1}
\Psi_{j} := ( M_{j,1}^{b} )^{T} \Phi_{j+1}, \qquad j \geq j_{0} , 
\end{equation}
with $( M_{j,1}^{b} )^{T} := (b_{k,l}^{j,b})_{k \in \nabla_{j}, l \in \Delta_{j+1}} \in \mathbb{R}^{|\nabla_{j}| \times |\Delta_{j+1}|} $. Similar relations hold for $\tilde{\Psi}_{j}$. Then $\langle \Psi_{j} , \tilde{\Phi}_{j} \rangle = 0$, $\langle \Phi_{j} , \tilde{\Psi}_{j} \rangle = 0$ and $\langle \Psi_{j} , \tilde{\Psi}_{j} \rangle = Id_{|\nabla_{j}|}$. Now let us turn to the definitions of boundary adapted quarks and quarklets. We start with the Schoenberg B-spline quarks. 

\begin{Definition}\label{def_bound_quark}
Let $ m \in \mathbb{N}   $, $ j \in \mathbb{N}  $ with $ j \geq j_{0}  $ and $  p \in \mathbb{N}_{0}    $. Then the p-th Schoenberg B-spline quark $  \varphi_{p,j,k}  $ is defined by
\[
\varphi_{p,j,k} := \left\{ \begin{array}{lll}
\Big ( \frac{2^{j} \cdot}{k + m}  \Big )^{p} \varphi_{j,k}  & \quad & \mbox{for} \qquad k = -m+1, \ldots , -1 ;
\\  
\Big ( \frac{2^{j} \cdot - k - \lfloor \frac{m}{2} \rfloor }{\lceil \frac{m}{2} \rceil}  \Big )^{p} \varphi_{j,k}  & \quad & \mbox{for}\qquad k = 0 , \ldots , 2^{j} - m ;
\\
\varphi_{p,j,2^{j}-m-k}(1 - \cdot ) & \quad & \mbox{for}\qquad k = 2^{j} - m +1 , \ldots , 2^{j} - 1 .
\\
\end{array}
\right.
\]
\end{Definition}
\noindent
Notice that the inner Schoenberg B-spline quarks are translated copies of the cardinal B-spline quarks, see Definition \ref{Bquark}. Now let us turn to the construction of the quarklets. For the inner quarklets we can proceed as in Subsection \ref{sec_inner_quark}. Let $  p \in \mathbb{N}_{0}   $ and $ j \in \mathbb{N}   $ with $ j \geq j_{0}  $. Let $ k \in  \nabla_{j}  $ with $ m - 1 \leq k \leq 2^{j} - m    $. We use the inner wavelets constructed by Primbs, see \cite{Pri1} and \cite{Pri2010}. They can be written as
\begin{equation}\label{wav_prim_bjk}
\psi_{j,k}^{b} := \sum_{l \in \Delta_{j+1}} b^{j,b}_{k,l} \varphi_{j+1,l} .
\end{equation}
To obtain the inner quarklets for $ m - 1 \leq k \leq 2^{j} - m    $ and $ l \in \Delta_{j+1}   $ we use the numbers $ b^{j ,b}_{k,l} $ given in \eqref{wav_prim_bjk}. We put $ b^{p, j ,b}_{k,l} := b^{j ,b}_{k,l}    $ and define the inner quarklets by
\begin{equation}\label{def_inner_quark}
\psi_{p,j,k}^{b} := \sum_{l \in \Delta_{j+1}}  b^{p, j ,b}_{k,l} \varphi_{p,j+1,l} .
\end{equation} 
If the inner Primbs wavelets have $ \tilde{m}  $ vanishing moments, also the inner quarklets defined in \eqref{def_inner_quark} have $ \tilde{m}  $ vanishing moments. This result can be found in \cite{DaKRaa}, see Lemma 2. In a next step we construct the boundary quarklets. Later on it will be very important that also they have vanishing moments. Therefore in general we can not use the boundary wavelets constructed by Primbs in \cite{Pri1}. A simple counterexample to illustrate the lack of vanishing moments can be found in \cite{SiDiss}, see page 67. Hence our strategy to construct the boundary quarklets reads as follows. We fix that they have $\tilde{m}$ vanishing moments and obtain a system of linear equations from that determination. Let us deal with the left boundary quarklets first. Then right boundary quarklets can be obtained via reflection. We work with $ k = 0, 1, \ldots, m-2    $.  We assume that each left boundary quarklet consists of  $ \tilde{m} + 1 $ quarks, which are either left boundary or inner quarks as given in Definition \ref{def_bound_quark}. Furthermore the $k$-th quarklet representation should begin at the leftmost but $k$-th quark. This leads to a $ \tilde{m} \times (\tilde{m} + 1)      $ linear system of equations. It can be written as
\begin{equation}\label{eq_bou_quarklet_b}
\sum_{l = -m+1 + k}^{-m+1 + k + \tilde{m}} b^{p,j,b}_{k,l} \int_{\mathbb{R}} x^{q} \varphi_{p,j+1,l}(x) dx = 0 , \qquad q = 0, 1, \ldots , \tilde{m} - 1 .
\end{equation} 
The resulting coefficient matrix is of size $ \tilde{m} \times (\tilde{m} + 1)      $ and has a nontrivial kernel, see Chapter 4.3 in \cite{SiDiss}. So we can find nontrivial solutions for \eqref{eq_bou_quarklet_b}. Consequently we are able to construct boundary quarklets with vanishing moments. There is the following definition, see also Definition 4.20 in \cite{SiDiss}. 
\begin{Definition}\label{def_bound_quarklet}
Let $  k = 0, 1, \ldots, m-2   $ and $ \tilde{m} \in \mathbb{N}   $ with $ \tilde{m} \geq m   $. Let $ j \in \mathbb{N}   $ with $ j \geq j_{0}  $ and $ p \in \mathbb{N}_{0}   $. If the vector $ \textbf{b}^{p,j,b}_{k} = (  b^{p,j,b}_{k,-m+1+ k}  , \ldots , b^{p,j,b}_{k,-m+1+ k + \tilde{m}}     ) \in \mathbb{R}^{\tilde{m} + 1}    $ with $ \textbf{b}^{p,j,b}_{k} \not = 0   $ is a solution for \eqref{eq_bou_quarklet_b}, then we define the $k$-th left boundary quarklet by 
\begin{equation}\label{eq_kleftboundaryquar}
\psi^{b}_{p,j,k} := \sum_{l = -m+1+ k}^{-m+1+ k + \tilde{m}} b^{p,j,b}_{k,l}  \varphi_{p,j+1,l} .
\end{equation}
\end{Definition}
\noindent
Here the parameter $k$ refers to the fact that we have $ m - 1  $ left boundary quarklets and $ m - 1  $ right boundary quarklets.
\begin{Remark}
The boundary quarklets given in Definition \ref{def_bound_quarklet} have been constructed in \cite{DaFKRaa}, see Section 2.4. A detailed study of their properties can be found in \cite{SiDiss}, see Chapter 4.3.
\end{Remark}
\noindent
Later it will be convenient to use a uniform notation that refers to both quarks and quarklets at the same time. For that purpose for $p \in \mathbb{N}_{0}$ and $k \in \nabla_{j_{0}-1}$ we write 
\begin{equation}\label{eq_def_qua_as_qualet}
\psi^{b}_{p,j_{0}-1,k} := \varphi_{p,j_{0},k} .
\end{equation}
The quarklets constructed in this section can be used to assemble quarklet systems that are frames for $ L_{2}((0,1))  $. To see this let us introduce some additional notation. We define the index set for the whole quarklet system by
\begin{equation}\label{index_quarkl_full1}
\nabla :=  \{  (p,j,k) : p,j \in \mathbb{N}_{0}, j \geq j_{0}-1, k \in \nabla_{j}     \} .
\end{equation}
It contains the Primbs basis index set
\begin{equation}\label{index_Primbs1_full1}
\nabla^{P} :=  \{  (0,j,k) : j \in \mathbb{N}_{0}, j \geq j_{0}-1, k \in \nabla_{j}     \} .
\end{equation}
The whole quarklet system itself based on the index set $ \nabla $ is given by
\begin{equation}\label{def_quarklsy_ful1}
\Psi^{b} := \{  \psi^{b}_{p,j,k} : (p,j,k) \in   \nabla    \} .
\end{equation}
Recall that for $j = j_{0}-1$ the system $\Psi^{b}$ contains the Schoenberg B-spline quarks, see Definition \ref{def_bound_quark} and \eqref{eq_def_qua_as_qualet}. For $j \geq j_{0}$ it consists of inner and boundary quarklets. Thereby for $k \in \nabla_{j} $ with $ m - 1 \leq k \leq 2^{j} - m    $ it refers to the inner quarklets given in \eqref{def_inner_quark}. Otherwise if $k \in \{ 0, 1, \ldots , m - 2   \}$ or $k \in \{ 2^{j} - m + 1 , \ldots , 2^{j} - 1   \}$ the system $ \Psi^{b} $ consists of left boundary quarklets or right boundary quarklets, respectively. The quarklets collected in the system $ \Psi^{b}  $ form a frame for $ L_{2}((0,1))   $, see Theorem 2.7 in \cite{DaFKRaa} and Theorem 4.23 in \cite{SiDiss}. 

\begin{Theorem}\label{mainresult1_Hsr_d=1}
Let $  \nabla  $ be the index set defined in \eqref{index_quarkl_full1} and $  \delta > 1 $. Then the weighted quarklet system 
\begin{equation}\label{set_L2_quark_weight1_d1}
\Psi^{b}_{L_{2}((0,1))} := \Big \{ ( p + 1 )^{- \frac{\delta}{2}}  \psi^{b}_{p,j,k} : (p,j,k) \in   \nabla \Big  \}
\end{equation}
is a frame for $ L_{2}((0,1))  $.
\end{Theorem}

\subsection{Bivariate Quarklets via Tensor Products}\label{subsec_biv_tensorwavelets}

In what follows we construct bivariate quarklets out of the univariate quarklets obtained in Section \ref{Subsec_bound_qua} via tensor product methods. To this end we follow \cite{DaFKRaa}, see Section 3.3. Hereinafter $i \in \{ 1, 2  \}$ always refers to the $i-$th Cartesian direction. Recall that each univariate quarklet $\psi_{p,j,k}^{b}$ can be identified via a triple $\lambda = (p,j,k)$ with $p \in \mathbb{N}_{0}$, $j \geq j_{0}-1$ and $ k \in \nabla_{j}$. Sometimes $\lambda$ also is called a quarklet index. To obtain bivariate quarklets we require quarklet indices for each Cartesian direction $i$. They are denoted by $\lambda_{i} = (p_{i},j_{i},k_{i})$. We put $\boldsymbol{\lambda} := ( \lambda_{1}, \lambda_{2} )$. The index set for the whole quarklet system concerning direction $i$ again is given by $\nabla$, see  \eqref{index_quarkl_full1}. Now for given quarklet indices $\lambda_{i}$ we define bivariate quarklets using tensor products of univariate quarklets. We put
\begin{equation}\label{def_tensor_quarklet}
\boldsymbol{\psi_{\lambda}} :=  \psi_{\lambda_{1}}^{b} \otimes \psi_{\lambda_{2}}^{b} .
\end{equation}
To address these bivariate tensor quarklets we define the index set
\begin{equation}\label{def_tensorqua_index}
\boldsymbol{\nabla} :=  \nabla \times \nabla .
\end{equation}
The collection of all bivariate quarklets that can be obtained via tensor products as described above is given by
\begin{equation}\label{def_collect_multqua}
\boldsymbol{\Psi} :=  \Psi^{b} \otimes \Psi^{b} = \Big \{ \boldsymbol{\psi_{\lambda}} : \boldsymbol{\lambda}  \in  \boldsymbol{\nabla} \Big \} .
\end{equation}
Here $\Psi^{b}$ refers to the whole system of univariate quarklets as defined in \eqref{def_quarklsy_ful1}. The bivariate quarklets collected in the set $  \boldsymbol{\Psi} $ can be used to construct tensor frames for $  L_{2}((0,1)^2)  $. Starting point for this is the observation
\begin{equation}
L_{2}((0,1)^2) = L_{2}((0,1)) \otimes_{2} L_{2}((0,1)) ,
\end{equation} 
see Theorem 1.39 in \cite{LightChe}, and Lemmas 1.34 - 1.36 in \cite{LightChe} for further explanations concerning $  \otimes_{2}   $. Based on this identity we can use the bivariate tensor quarklets to obtain frames for $  L_{2}((0,1)^2)  $. The following result can be found in \cite{DaFKRaa}, see Theorem 3.10.

\begin{Theorem}\label{thm_biv_tensor_frame_L2}
Let $ m \geq 2   $ and $ \tilde{m} \in \mathbb{N} $  with $  \tilde{m} \geq m  $ and $ m + \tilde{m} \in 2 \mathbb{N}   $. Let $  \Psi^{b}_{L_{2}((0,1))}       $ be the weighted  quarklet system given in \eqref{set_L2_quark_weight1_d1}. Let $ \delta > 1   $.  Then the family 
\begin{equation}
\boldsymbol{\Psi}_{L_{2}((0,1)^2)} :=  \Psi^{b}_{L_{2}((0,1))} \otimes \Psi^{b}_{L_{2}((0,1))} = \Big \{  w_{\boldsymbol{\lambda}}^{-1} \boldsymbol{\psi_{\lambda}} : \boldsymbol{\lambda}  \in  \boldsymbol{\nabla} := \nabla \times \nabla \Big \} 
\end{equation}
with weights 
\begin{equation}\label{tensor_frame_weights1}
w_{\boldsymbol{\lambda}} := ( p_{1} + 1)^{\frac{\delta}{2}}  ( p_{2} + 1)^{\frac{\delta}{2}} 
\end{equation}
is a quarkonial tensor frame for $ L_{2}((0,1)^2)  $.
\end{Theorem}
It is possible to represent every $ f \in L_2((0,1)^2)  $ in terms of the bivariate tensor quarklets introduced above. Indeed, using Theorem \ref{thm_biv_tensor_frame_L2} and the properties of the frame operator we find that for each $ f \in L_2((0,1)^2)  $ there exists at least one sequence $ \{  c_{\boldsymbol{\lambda}}  \}_{\boldsymbol{\lambda} \in  \boldsymbol{\nabla} } \in \ell_{2}( \boldsymbol{\nabla} ) $ such that 
\begin{equation}\label{eq_canon_dual_frame_2_vorl1}
f = \sum_{ \boldsymbol{\lambda} \in  \boldsymbol{\nabla}  }   c_{\boldsymbol{\lambda}} w_{\boldsymbol{\lambda}}^{-1}  \boldsymbol{\psi_{\lambda}} . 
\end{equation}

%Therefore we introduce the synthesis operator 
%\begin{align*}
%T : \ell_{2}( \boldsymbol{\nabla} ) \rightarrow L_{2}((0,1)^2) , \qquad \{  c_{\boldsymbol{\lambda}}  \}_{\boldsymbol{\lambda} \in  \boldsymbol{\nabla} } \mapsto \sum_{ \boldsymbol{\lambda} \in  \boldsymbol{\nabla}  } c_{\boldsymbol{\lambda}} w_{\boldsymbol{\lambda}}^{-1} \boldsymbol{\psi_{\lambda}}
%\end{align*}
%and its adjoint, the analysis operator, which is given by 
%\begin{align*}
%T^{\star} : L_{2}((0,1)^2) \rightarrow \ell_{2}( \boldsymbol{\nabla} ) , \qquad f \mapsto \{ \langle f , w_{\boldsymbol{\lambda}}^{-1}  \boldsymbol{\psi_{\lambda}}  \rangle_{L_{2}((0,1)^2)}   \}_{\boldsymbol{\lambda} \in  \boldsymbol{\nabla} } .
%\end{align*}
%A combination of both leads to the so-called frame operator
%\begin{align*}
%S : L_{2}((0,1)^2) \rightarrow  L_{2}((0,1)^2), \qquad f \mapsto Sf := T T^{\star} f = \sum_{ \boldsymbol{\lambda} \in  \boldsymbol{\nabla}  }\langle f , w_{\boldsymbol{\lambda}}^{-1}  \boldsymbol{\psi_{\lambda}}  \rangle_{L_{2}((0,1)^2)} w_{\boldsymbol{\lambda}}^{-1}  \boldsymbol{\psi_{\lambda}} . 
%\end{align*}
%The frame operator is bounded and invertible and we have
%\begin{equation}\label{eq_canon_dual_frame-1}
%f = \sum_{ \boldsymbol{\lambda} \in  \boldsymbol{\nabla}  }\langle f , S^{-1} ( w_{\boldsymbol{\lambda}}^{-1} \boldsymbol{\psi_{\lambda}})  \rangle_{L_{2}((0,1)^2)} w_{\boldsymbol{\lambda}}^{-1}  \boldsymbol{\psi_{\lambda}} 
%\end{equation}
%for all $ f \in L_2((0,1)^2)  $. However, this representation is not unique.
%\end{comment}

\section{The Concept of Bivariate Quarklet Trees}\label{sec_biv_tre_main111}

It is the main goal of this paper to approximate functions $  f \in L_2((0,1)^2)  $ via bivariate tensor quarklets using tree approximation techniques. For that purpose we have to introduce the concept of bivariate quarklet trees. Univariate quarklet trees have been introduced in \cite{DaHoRaVo}, see Section 2. However, it turns out, that in the case of two dimensions the situation is much more complicated, since then several new phenomena show up. Consequently, also our definition of bivariate quarklet trees is much more intricate than the univariate counterpart given in \cite{DaHoRaVo}. As an important intermediate step, we explain the concept of bivariate wavelet trees. For that purpose in a first step we recall the idea of reference rectangles.

\subsection{Reference Rectangles}\label{subsec_ref_rectangle1}

Reference rectangles are two-dimensional generalizations of reference intervals, which have been recalled in \cite{DaHoRaVo}, see Section 2.2. For $ j_{1} , j_{2} \in \mathbb{N}_{0}  $ and $ k_{1} \in \{ 0, 1, \ldots , 2^{j_{1}} - 1 \}   $,  $ k_{2} \in \{ 0, 1, \ldots , 2^{j_{2}} - 1 \}   $ reference rectangles are defined by:
\begin{equation}
R^{(j_{1},k_{1})}_{(j_{2},k_{2})} := [ 2^{-j_{1}} k_{1} , 2^{-j_{1}} ( k_{1} + 1 ) ) \times  [ 2^{-j_{2}} k_{2} , 2^{-j_{2}} ( k_{2} + 1 ) ) .
\end{equation}
The reference rectangles have some interesting  properties, which will be important for us later on. Some of them are collected in the following list:
\begin{itemize}

\item[(i)] Let $ j_{1} , j_{2} \in \mathbb{N}_{0}  $ be fixed. Then we observe
\begin{align*}
\bigcup_{k_{1} = 0}^{2^{j_{1}} - 1} \bigcup_{k_{2} = 0}^{2^{j_{2}} - 1} R^{(j_{1},k_{1})}_{(j_{2},k_{2})} = [0,1) \times [0,1) = [0,1)^{2} .
\end{align*}
Moreover, the above partition of $ [0,1)^{2}  $ is disjoint.

\item[(ii)] Let $ j_{1} , j_{2} \in \mathbb{N}_{0}  $, $ k_{1} \in \{ 0, 1, \ldots , 2^{j_{1}} - 1 \}   $ and  $ k_{2} \in \{ 0, 1, \ldots , 2^{j_{2}} - 1 \}   $ be fixed. Then there are two different possibilities to disassemble $ R^{(j_{1},k_{1})}_{(j_{2},k_{2})} $ into two finer reference rectangles of the next higher level. More precisely we have
\begin{equation}\label{eq_dec_twopos}
R^{(j_{1},k_{1})}_{(j_{2},k_{2})} = R^{(j_{1}+1, 2 k_{1})}_{(j_{2},k_{2})} \cup R^{(j_{1}+1, 2 k_{1} + 1)}_{(j_{2},k_{2})} \qquad \mbox{and} \qquad R^{(j_{1},k_{1})}_{(j_{2},k_{2})} = R^{(j_{1},k_{1})}_{(j_{2}+1, 2 k_{2})} \cup R^{(j_{1},k_{1})}_{(j_{2}+1, 2 k_{2} + 1)} .
\end{equation} 
These partitions are disjoint. We call 
\begin{align*}
R^{(j_{1}+1, 2 k_{1})}_{(j_{2},k_{2})}  \; \mbox{and} \;   R^{(j_{1}+1, 2 k_{1} + 1)}_{(j_{2},k_{2})} \;  \mbox{the children of} \; R^{(j_{1},k_{1})}_{(j_{2},k_{2})}   
\end{align*}
in direction $i = 1$. Similar we call
\begin{align*}
R^{(j_{1},k_{1})}_{(j_{2}+1, 2 k_{2})} \; \mbox{and} \;  R^{(j_{1},k_{1})}_{(j_{2}+1, 2 k_{2} + 1)} \;  \mbox{the children of}  \; R^{(j_{1},k_{1})}_{(j_{2},k_{2})} 
\end{align*}
in direction $ i = 2$.

\end{itemize}
Let $ j_{1} , j_{2} \in \mathbb{N}_{0}  $ and $ k_{1} \in \{ 0, 1, \ldots , 2^{j_{1}} - 1 \}   $,  $ k_{2} \in \{ 0, 1, \ldots , 2^{j_{2}} - 1 \}   $ be as above. Then each reference rectangle $ R^{(j_{1},k_{1})}_{(j_{2},k_{2})}   $ refers to a bivariate wavelet index $  \boldsymbol{\lambda} := ( (j_{1},k_{1}), (j_{2},k_{2}) ) $ and vice versa. Recall, that for $j_{1} \geq j_{0} - 1$, $j_{2} \geq j_{0} - 1$ we recover the bivariate wavelet indices in the index set $  \boldsymbol{\nabla}  $, when we put $ p_{1} = p_{2} = 0   $. Consequently, there is a connection between reference rectangles and bivariate tensor quarklets of the lowest polynomial degree. Much more details concerning this topic can be found in Section \ref{sec_ref_rec_wavl} below.

\subsection{Bivariate Wavelet Trees}\label{subsec_wav_tree}

In order to introduce bivariate quarklet trees, it is an important intermediate step to deal with bivariate wavelet trees. Below we generalize the theory explained in \cite{DaHoRaVo}, see Section 2.2, where the univariate case has been investigated. To define bivariate wavelet trees, at first we require an ancestor-descendant relation concerning the reference rectangles. For that purpose let $ j_{1} , j_{2} \in \mathbb{N}_{0}  $, $ k_{1} \in \{ 0, 1, \ldots , 2^{j_{1}} - 1 \}   $ and  $ k_{2} \in \{ 0, 1, \ldots , 2^{j_{2}} - 1 \}   $. Moreover, let $ \tilde{j}_{1} \geq j_{1} , \tilde{j}_{2} \geq j_{2}   $, $ \tilde{k}_{1} \in \{ 0, 1, \ldots , 2^{\tilde{j}_{1}} - 1 \}   $ and  $ \tilde{k}_{2} \in \{ 0, 1, \ldots , 2^{\tilde{j}_{2}} - 1 \}   $ be such that
\begin{equation}\label{eq_example_situation1}
R^{(\tilde{j}_{1},\tilde{k}_{1})}_{(\tilde{j}_{2},\tilde{k}_{2})} \subseteq R^{(j_{1},k_{1})}_{(j_{2},k_{2})} .
\end{equation}
Then there exists a sequence of decompositions of the form \eqref{eq_dec_twopos} to obtain a partition of 
\begin{align*}
R^{(j_{1},k_{1})}_{(j_{2},k_{2})} \; \mbox{which also contains} \; R^{(\tilde{j}_{1},\tilde{k}_{1})}_{(\tilde{j}_{2},\tilde{k}_{2})}.
\end{align*}
Given an arbitrary reference rectangle, there are always two refinement options, one for each Cartesian direction. Consequently, when carrying out a sequence of several refinement steps, it is possible to obtain the same partition following different refinement strategies, by interchanging the order of refinements in direction $i = 1$ or $ i = 2$. However, in order to obtain an efficient quarklet tree algorithm which is near-best, it becomes necessary to restrict the available refinement options. For that purpose we introduce an additional parameter $\alpha \in \{ 0, 1, 2 \}$. It will be associated to a wavelet index $\boldsymbol{\lambda} $ and describes the refinement options of the corresponding reference rectangle. So $\alpha = 0$ means that both refinement options given in \eqref{eq_dec_twopos} are permitted. In the case $\alpha = 1$ only space refinement in direction $ i = 1  $ is allowed. Finally, $\alpha = 2$ means that only refinement in direction $i = 2$ is permitted. In consequence we define an enhanced wavelet index $\tilde{\boldsymbol{\lambda}} := ((j_{1},k_{1}), (j_{2},k_{2}), \alpha )$ including refinement options. For a given enhanced wavelet index $ \tilde{\boldsymbol{\lambda}}  $ with refinement options $\alpha$ we use the notation $  \alpha(\tilde{\boldsymbol{\lambda}}) := \alpha   $. Now we are prepared to introduce an ancestor-descendant relation concerning the reference rectangles. If for bivariate wavelet indices $  \tilde{\boldsymbol{\lambda}} = ( (j_{1},k_{1}), (j_{2},k_{2}), \alpha ) $ and $  \tilde{\boldsymbol{\mu}} = ( (\tilde{j}_{1},\tilde{k}_{1}), (\tilde{j}_{2},\tilde{k}_{2}), \tilde{\alpha} ) $ we have 
\begin{align*}
\mbox{either} \qquad R^{(\tilde{j}_{1},\tilde{k}_{1})}_{(\tilde{j}_{2},\tilde{k}_{2})} \subset R^{(j_{1},k_{1})}_{(j_{2},k_{2})} \qquad \mbox{or} \qquad R^{(\tilde{j}_{1},\tilde{k}_{1})}_{(\tilde{j}_{2},\tilde{k}_{2})} = R^{(j_{1},k_{1})}_{(j_{2},k_{2})} \; \mbox{with} \; \sgn{\tilde{\alpha}} > \sgn{\alpha}
\end{align*}
we will use the notation $\tilde{\boldsymbol{\mu}} \succ \tilde{\boldsymbol{\lambda}} $ and say that $  \tilde{\boldsymbol{\mu}} $ is a descendant of $ \tilde{\boldsymbol{\lambda}}  $. Conversely, we will call $ \tilde{\boldsymbol{\lambda}}  $ an ancestor of $  \tilde{\boldsymbol{\mu}}  $. By $ \tilde{\boldsymbol{\mu}} \succeq \tilde{\boldsymbol{\lambda}}$ we mean that $\tilde{\boldsymbol{\mu}}$ is either a descendant of $  \tilde{\boldsymbol{\lambda}} $ or equal to $\tilde{\boldsymbol{\lambda}}$. One key tool for adaptive bivariate quarklet tree approximation is the possibility to carry out local space refinement. Due to the two Cartesian directions there are different options, which can be found in the listing below, whereby we use a distinction of cases concerning the parameter $\alpha$.

\begin{itemize}
\item[(LSR.a)] \underline{The case $  \alpha = 0   $.}

\begin{itemize}
\item[(LSR.a.1)] \underline{Space refinement in direction $i = 1$.}

Let $ \tilde{\boldsymbol{\lambda}} = ( (j_{1},k_{1}), (j_{2},k_{2}), 0 )   $ be given. Then $ \tilde{\boldsymbol{\lambda}} $ can be refined by adding 3 children, namely
\begin{align*}
\{ ( (j_{1} + 1, 2 k_{1}), (j_{2},k_{2}), 0 ) ,  ( (j_{1} + 1 , 2 k_{1} + 1), (j_{2},k_{2}), 0 )  ,  ( (j_{1},k_{1}), (j_{2},k_{2}), 2 ) \} .
\end{align*} 

\item[(LSR.a.2)] \underline{Space refinement in direction $i = 2$.}

Let $ \tilde{\boldsymbol{\lambda}} = ( (j_{1},k_{1}), (j_{2},k_{2}), 0 )   $ be given. Then $ \tilde{\boldsymbol{\lambda}} $ can be refined by adding 3 children, namely
\begin{align*}
\{ ( (j_{1} , k_{1}), (j_{2} + 1 , 2 k_{2}), 0 )  ,  ( (j_{1}  , k_{1} ), (j_{2} + 1 , 2 k_{2} + 1), 0 ) ,  ( (j_{1},k_{1}), (j_{2},k_{2}), 1 ) \}  .
\end{align*} 

\end{itemize}
\item[(LSR.b)] \underline{The case $  \alpha = 1 $. Space refinement in direction $i = 1$.}

Let $ \tilde{\boldsymbol{\lambda}} = ( (j_{1},k_{1}), (j_{2},k_{2}), 1 )   $ be given. Then $ \tilde{\boldsymbol{\lambda}}  $ can be refined by adding 2 children, namely
\begin{align*}
\{  ( (j_{1} + 1, 2 k_{1}), (j_{2},k_{2}), 0 )  ,  ( (j_{1} + 1 , 2 k_{1} + 1), (j_{2},k_{2}), 0 ) \} .
\end{align*}

\item[(LSR.c)] \underline{The case $  \alpha = 2 $. Space refinement in direction $i = 2$.} 

Let $ \tilde{\boldsymbol{\lambda}} = ( (j_{1},k_{1}), (j_{2},k_{2}), 2 )   $ be given. Then $  \tilde{\boldsymbol{\lambda}} $ can be refined by adding 2 children, namely
\begin{align*}
\{   ( (j_{1} , k_{1}), (j_{2} + 1 , 2 k_{2}), 0 )  , ( (j_{1}  , k_{1} ), (j_{2} + 1 , 2 k_{2} + 1), 0 ) \} .
\end{align*}

\end{itemize}

Here the abbreviation LSR stands for local space refinement. In the following sections the notations (LSR.a.1), (LSR.a.2), (LSR.b) and (LSR.c) will show up many times, and always refer to the refinement options described above. Looking at the cases (LSR.a.1) and (LSR.a.2) there are three children, at which the third one does not stand for a refinement, but gives us the possibility to carry out a refinement in the other direction later. In comparison to the univariate case this is a substantial difference, since there we always have two children, see Section 2.2 in \cite{DaHoRaVo}. In the cases (LSR.b) and (LSR.c) due to $\alpha \not = 0$ the number of children reduces to two. In the Figures \ref{Fig_1} and \ref{Fig_2} the different refinement options (LSR.a.1), (LSR.a.2), (LSR.b) and (LSR.c) are illustrated.

\vspace{0,5 cm}

\begin{figure}[h]
\centering
\begin{tikzpicture}[thick]

\draw (0,2) -- (0,0) ;
\draw (2,2) -- (2,0) ;
\draw (0,2) -- (2,2) ;
\draw (0,0) -- (2,0) ;

\node at (2.5,1) {$ \tilde{\boldsymbol{\lambda}}_{1} $} ;

\draw[->] (0, -0.2 ) -- (-1, - 0.8) ;
\draw[->] (0, -0.2 ) -- (0, - 0.8) ;
\draw[->] (2, -0.2 ) -- (2.5, - 0.8) ;

\draw (-1.5,-1) -- (-1.5,-3) ;
\draw (0.5,-1) -- (0.5,-3) ;
\draw (-1.5,-1) -- (0.5,-1) ;
\draw (-1.5,-3) -- (0.5,-3) ;
\draw ( -0.5 , -1 ) -- ( -0.5  , -3  )  ;

\node at ( -1  , -3.5 ){$ \tilde{\boldsymbol{\lambda}}_{2} $} ;
\node at ( 0  , -3.5 ){$ \tilde{\boldsymbol{\lambda}}_{3} $} ;
\node at ( 2.5  , -3.5 ){$ \tilde{\boldsymbol{\lambda}}_{4} $} ;

\draw (1.5,-1) -- (1.5,-3) ;
\draw (3.5,-1) -- (3.5,-3) ;
\draw (1.5,-1) -- (3.5,-1) ;
\draw (1.5,-3) -- (3.5,-3) ;

\draw (7,2) -- (7,0) ;
\draw (9,2) -- (9,0) ;
\draw (7,2) -- (9,2) ;
\draw (7,0) -- (9,0) ;

\node at (9.5,1) {$ \tilde{\boldsymbol{\mu}}_{1} $} ;

\draw[->] (7, -0.2 ) -- (6, - 2.2) ;
\draw[->] (7, -0.2 ) -- (7, - 0.8) ;
\draw[->] (9, -0.2 ) -- (9.5, - 0.8) ;

\draw (5.5,-1) -- (5.5,-3) ;
\draw (7.5,-1) -- (7.5,-3) ;
\draw (5.5,-1) -- (7.5,-1) ;
\draw (5.5,-3) -- (7.5,-3) ;
\draw ( 5.5 , -2 ) -- ( 7.5 , -2  )  ;

\node at ( 5  , -2.5 ){$ \tilde{\boldsymbol{\mu}}_{2} $} ;
\node at ( 5  , -1.5 ){$ \tilde{\boldsymbol{\mu}}_{3} $} ;
\node at ( 9.5  , -3.5 ){$ \tilde{\boldsymbol{\mu}}_{4} $} ;

\draw (8.5,-1) -- (8.5,-3) ;
\draw (10.5,-1) -- (10.5,-3) ;
\draw (8.5,-1) -- (10.5,-1) ;
\draw (8.5,-3) -- (10.5,-3) ;

\end{tikzpicture}
\caption{Refinement Strategie (LSR.a.1) with $ \tilde{\boldsymbol{\lambda}}_{1} =  ( (j_{1},k_{1}), (j_{2},k_{2}), 0 ) $, $ \tilde{\boldsymbol{\lambda}}_{2} =  ( (j_{1} + 1, 2 k_{1}), (j_{2},k_{2}), 0 )  $, $ \tilde{\boldsymbol{\lambda}}_{3} =   ((j_{1} + 1 , 2 k_{1} + 1), (j_{2},k_{2}), 0 )    $, $ \tilde{\boldsymbol{\lambda}}_{4}     =  ( (j_{1},k_{1}), (j_{2},k_{2}), 2 )  $ and Refinement Strategie (LSR.a.2) with $ \tilde{\boldsymbol{\mu}}_{1} = ( (j_{1},k_{1}), (j_{2},k_{2}), 0 )  $, $ \tilde{\boldsymbol{\mu}}_{2} = ( (j_{1} , k_{1}), (j_{2} + 1 , 2 k_{2}), 0 )   $, $ \tilde{\boldsymbol{\mu}}_{3} =    ( (j_{1}  , k_{1} ), (j_{2} + 1 , 2 k_{2} + 1), 0 )    $, $ \tilde{\boldsymbol{\mu}}_{4}     = ( (j_{1},k_{1}), (j_{2},k_{2}), 1 )   $.  }
\label{Fig_1}
\end{figure}
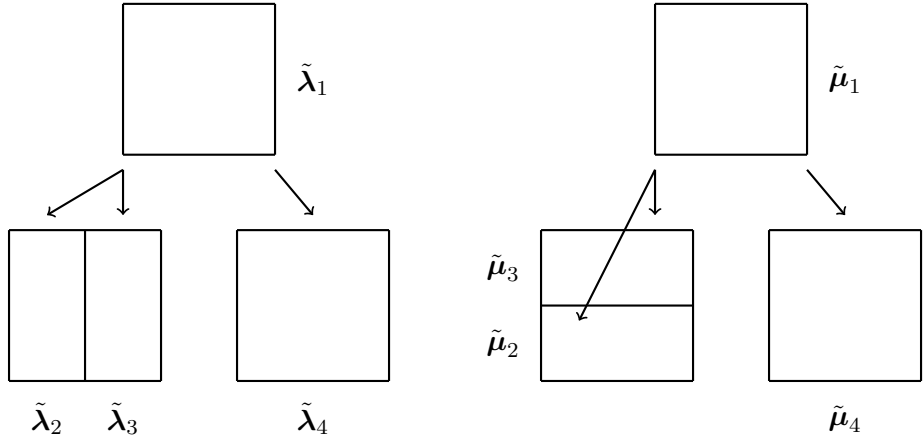

\begin{figure}[h]
\centering
\begin{tikzpicture}[thick]

\draw (0,2) -- (0,0) ;
\draw (2,2) -- (2,0) ;
\draw (0,2) -- (2,2) ;
\draw (0,0) -- (2,0) ;

\node at (2.5,1) {$ \tilde{\boldsymbol{\lambda}}_{1} $} ;

\draw[->] (1, -0.2 ) -- (0.5, - 0.8) ;
\draw[->] (1, -0.2 ) -- (1.5, - 0.8) ;

\draw (0,-1) -- (0,-3) ;
\draw (2,-1) -- (2,-3) ;
\draw (0,-1) -- (2,-1) ;
\draw (0,-3) -- (2,-3) ;
\draw ( 1 , -1 ) -- ( 1  , -3  )  ;

\node at ( 0.5  , -3.5 ){$ \tilde{\boldsymbol{\lambda}}_{2} $} ;
\node at ( 1.5 , -3.5 ){$ \tilde{\boldsymbol{\lambda}}_{3} $} ;

\draw (7,2) -- (7,0) ;
\draw (9,2) -- (9,0) ;
\draw (7,2) -- (9,2) ;
\draw (7,0) -- (9,0) ;

\node at (9.5,1) {$ \tilde{\boldsymbol{\mu}}_{1} $} ;

\draw[->] (8, -0.2 ) -- (7.5, - 2.2) ;
\draw[->] (8, -0.2 ) -- (8.5, - 0.8) ;

\draw (7,-1) -- (7,-3) ;
\draw (9,-1) -- (9,-3) ;
\draw (7,-1) -- (9,-1) ;
\draw (7,-3) -- (9,-3) ;
\draw (7 , -2 ) -- ( 9 , -2  )  ;

\node at ( 9.5  , -2.5 ){$ \tilde{\boldsymbol{\mu}}_{2} $} ;
\node at ( 9.5  , -1.5 ){$ \tilde{\boldsymbol{\mu}}_{3} $} ;

\end{tikzpicture}
\caption{Refinement Strategie (LSR.b) with $ \tilde{\boldsymbol{\lambda}}_{1} =  ( (j_{1},k_{1}), (j_{2},k_{2}), 1 )  $, $ \tilde{\boldsymbol{\lambda}}_{2} =  ( (j_{1} + 1, 2 k_{1}), (j_{2},k_{2}), 0 )  $, $ \tilde{\boldsymbol{\lambda}}_{3} =  ( (j_{1} + 1 , 2 k_{1} + 1), (j_{2},k_{2}), 0 )   $ and Refinement Strategie (LSR.c) with $ \tilde{\boldsymbol{\mu}}_{1} = ( (j_{1},k_{1}), (j_{2},k_{2}), 2 )  $, $ \tilde{\boldsymbol{\mu}}_{2} = ( (j_{1} , k_{1}), (j_{2} + 1 , 2 k_{2}), 0 )  $, $ \tilde{\boldsymbol{\mu}}_{3} =  ( (j_{1}  , k_{1} ), (j_{2} + 1 , 2 k_{2} + 1), 0 )   $.  }
\label{Fig_2}
\end{figure}
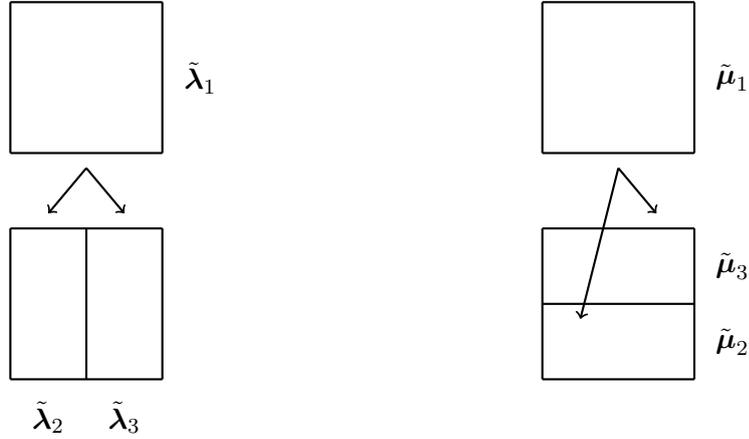

When looking at the situation described in \eqref{eq_example_situation1} in many cases there exist various refinement strategies only using (LSR.a.1), (LSR.a.2), (LSR.b) and (LSR.c), such that
\begin{align*}
R^{(\tilde{j}_{1},\tilde{k}_{1})}_{(\tilde{j}_{2},\tilde{k}_{2})} \; \mbox{is obtained out of} \; R^{(j_{1},k_{1})}_{(j_{2},k_{2})}  .
\end{align*} 
Consequently, it is possible to find refinement strategies such that certain reference rectangles have more than one parent. However, this causes grave difficulties for the theory we will develop below. Therefore we have to introduce an additional refinement rule which guarantees that each wavelet index has exactly one parent. There are different possibilities how this can be done, whereby we select the following. For an enhanced wavelet index $\tilde{\boldsymbol{\lambda}} := ((j_{1},k_{1}), (j_{2},k_{2}), \alpha )$ we introduce the notation $ | \tilde{\boldsymbol{\lambda}} | :=  | \boldsymbol{\lambda} | :=  j_{1} + j_{2}  $. We can formulate the following refinement rule. 

\begin{itemize}
\item[(UPC)] Let an enhanced wavelet index $\tilde{\boldsymbol{\lambda}} := ((j_{1},k_{1}), (j_{2},k_{2}), \alpha )$ with $ j_{1} < | \tilde{\boldsymbol{\lambda}} |$ be given. Then for this wavelet index only a refinement in direction $i = 2$ is allowed. 
\end{itemize} 

Here (UPC) stands for unique parent condition.

\begin{Remark}
The condition (UPC) implies that each wavelet index has exactly one parent. This is very important for the theory developed below. Condition (UPC) is well suited for the approximation of the bivariate function given in Section \ref{Sec_aniso_example} below. Therefore we stick with (UPC) in what follows. However, for some other test functions it seems to be reasonable to replace (UPC) by other unique parent conditions in order to obtain balanced quarklet trees.    
\end{Remark}

The (UPC) implies the following very important observation. 

\begin{Lemma}
Let an enhanced wavelet index $\tilde{\boldsymbol{\lambda}} := ((j_{1},k_{1}), (j_{2},k_{2}), 0  )$ with $ j_{1} , j_{2} \in \mathbb{N}_{0}  $, $ k_{1} \in \{ 0, 1, \ldots , 2^{j_{1}} - 1 \}   $ and  $ k_{2} \in \{ 0, 1, \ldots , 2^{j_{2}} - 1 \}   $ be given. Then there exists exactly one sequence of $|\tilde{\boldsymbol{\lambda}} |  $ refinement steps of the form (LSR.a.1), (LSR.a.2), (LSR.b) and (LSR.c) taking into account (UPC) and starting at $((0,0),(0,0),0)$, such that $ \tilde{\boldsymbol{\lambda}} $ is produced.    
\end{Lemma}

\begin{proof}
For the proof let $\tilde{\boldsymbol{\lambda}} := ((j_{1},k_{1}), (j_{2},k_{2}), 0  )$ be given. At first we construct a sequence of exactly $| \tilde{\boldsymbol{\lambda}} | $ refinement steps starting at $((0,0),(0,0),0)$ such that $ \tilde{\boldsymbol{\lambda}}  $ is obtained. For that purpose we start with $j_{1} $ refinement steps of kind (LSR.a.1) in direction $i = 1$ to obtain $  ((j_{1},k_{1}), (0,0), 0 )   $. Then we carry out $ j_{2}  $ refinement steps of the form (LSR.a.2) in direction $i = 2$ to get $   ((j_{1},k_{1}), (j_{2},k_{2}), 0 )   $. This strategy uses exactly $ | \tilde{\boldsymbol{\lambda}} |  $ steps and does not violate condition (UPC). Now assume that there exists another strategy with exactly $ | \tilde{\boldsymbol{\lambda}} |  $ steps such that $   ((j_{1},k_{1}), (j_{2},k_{2}), 0 )   $ is obtained. However, then the order of refinements in directions $i=1$ and $i=2$ must be different. This contradicts (UPC). Consequently the strategy described above is unique. 
\end{proof}

In what follows we only work with refinement strategies where (UPC) is fulfilled in each step. We observe that (UPC) affects the availability and shape of the refinement options (LSR.a.1), (LSR.a.2), (LSR.b) and (LSR.c).  So there exist $ \tilde{\boldsymbol{\lambda}}   $ for that (LSR.a.1) and (LSR.b) are not available at all. Moreover, for some $ \tilde{\boldsymbol{\lambda}}   $ due to (UPC) the number of children showing up in (LSR.a.2) reduces. Nevertheless, also if (UPC) holds, it is not possible to generally reformulate the refinement options (LSR.a.1), (LSR.a.2), (LSR.b) and (LSR.c), since there also exist $ \tilde{\boldsymbol{\lambda}}   $ where they remain unchanged.  In what follows, if for given $ \tilde{\boldsymbol{\lambda}}   $ a refinement option (or the occurrence of a child) contradicts (UPC), we call it \textit{not available} and ignore it for our considerations. Let us remark, that if (UPC) or a comparable rule which ensures the uniqueness of a parent does not hold, major parts of the theory presented below would not work any more. Now we have collected all tools to establish tree structured index sets. For that purpose we put 
\begin{align*}
\boldsymbol{\Lambda}_{0} := \Big  \{ ( (j_{1},k_{1}), (j_{2},k_{2}) ) : j_{1} , j_{2} \in \mathbb{N}_{0},  k_{1} \in \{ 0, 1, \ldots , 2^{j_{1}} - 1 \} , k_{2} \in \{ 0, 1, \ldots , 2^{j_{2}} - 1 \}  \Big   \} 
\end{align*}
and
\begin{align*}
\tilde{\boldsymbol{\Lambda}}_{0}  := & \Big  \{ ( (j_{1},k_{1}), (j_{2},k_{2}), \alpha ) : j_{1} , j_{2} \in \mathbb{N}_{0},  k_{1} \in \{ 0, 1, \ldots , 2^{j_{1}} - 1 \} , k_{2} \in \{ 0, 1, \ldots , 2^{j_{2}} - 1 \}, \\
& \alpha \in \{ 0, 1, 2 \}  \Big   \} . 
\end{align*}
Then bivariate wavelet trees can be defined in the following way.

\begin{Definition}\label{def_wav_tre}
Let $\mathcal{T} \subset \tilde{\boldsymbol{\Lambda}}_0 $ be an index set. The set $  \mathcal{T} $  is called a \emph{tree (of bivariate wavelet indices)} if the following conditions are fulfilled:
\begin{itemize}
\item[(i)] There exists an index $ \mathcal{R} = \tilde{\boldsymbol{\lambda}} \in \mathcal{T}   $ such that for all $ \tilde{\boldsymbol{\eta}} \in \mathcal{T}   $ we have $\tilde{\boldsymbol{\eta}} \succeq \tilde{\boldsymbol{\lambda}} $. This index is called \emph{root} of $ \mathcal{T} $.

\item[(ii)] The set $\mathcal{T} \subset \tilde{\boldsymbol{\Lambda}}_0 $ can be generated via a sequence of space refinements starting at $\mathcal{R} $, whereby in each step only space refinements as described in (LSR.a.1), (LSR.a.2), (LSR.b) or (LSR.c) are carried out and in each refinement step condition (UPC) is fulfilled.
\end{itemize}
 
\end{Definition}
Concerning bivariate wavelet trees the following terms will be important for us later.

\begin{Definition}\label{def_tre_terms}
Let $\mathcal{T} \subset \tilde{\boldsymbol{\Lambda}}_0   $ be a bivariate wavelet tree.
	\begin{itemize}
		\item[(i)]
		An index $\tilde{\boldsymbol{\lambda}} \in \mathcal{T}$ is called \emph{node} of $\mathcal{T}$. 
		\item[(ii)]
		We set
		\begin{align*}
		\mathcal{V}(\mathcal{T}) \coloneqq & \{ \tilde{\boldsymbol{\lambda}} \in \mathcal{T}  : \exists \mbox{ refinement strategy (LSR.a.1), (LSR.a.2), (LSR.b) or (LSR.c) } \\
& \mbox{for $ \tilde{\boldsymbol{\lambda}} $ such that its application to $ \tilde{\boldsymbol{\lambda}} $ generates at least 1 new child $\tilde{\boldsymbol{\eta}} \not \in \mathcal{T} $  }  \}. 
		\end{align*}
		The elements of $\mathcal{V}(\mathcal{T})$ are called \emph{leaves} of $  \mathcal{T} $. The set $\mathcal{T} \backslash \mathcal{V}(\mathcal{T})$ refers to the \emph{inner nodes}.

	\end{itemize}
\end{Definition}
Let us remark that the set $ \mathcal{V}(\mathcal{T})  $ collects all nodes of the bivariate wavelet tree $  \mathcal{T} $, which have not been refined during the creation process of $  \mathcal{T} $. Consequently, for these nodes a later refinement is still possible. On the other hand, $\mathcal{T} \backslash \mathcal{V}(\mathcal{T})$ gathers all nodes of $ \mathcal{T}  $ which have been refined during the creation of $  \mathcal{T} $. For these nodes no further refinement is allowed. For later use for $ \tilde{\boldsymbol{\lambda}} \in \tilde{\boldsymbol{\Lambda}}_{0}  $ we also define an infinite set $\mathcal{J}_{\tilde{\boldsymbol{\lambda}}}$ of bivariate wavelet indices which is given by 
\begin{align*}
\mathcal{J}_{\tilde{\boldsymbol{\lambda}}}  := & \{ \tilde{\boldsymbol{\mu}} \succeq \tilde{\boldsymbol{\lambda}} :  \exists \mbox{ sequence of refinement strategies (LSR.a.1), (LSR.a.2), (LSR.b) and} \\
& \mbox{(LSR.c) starting at $  \tilde{\boldsymbol{\lambda}}   $, where (UPC) holds in each step, such that $ \tilde{\boldsymbol{\mu}} $ is obtained} \} . 
\end{align*}
In other words the set $\mathcal{J}_{\tilde{\boldsymbol{\lambda}}}$ collects all nodes which belong to some infinite complete bivariate wavelet tree rooted at $  \tilde{\boldsymbol{\lambda}}  $. Notice that for the description of $\mathcal{J}_{\tilde{\boldsymbol{\lambda}}}$ we have to investigate infinitely many bivariate wavelet trees, since choosing either (LSR.a.1) or (LSR.a.2) in a refinement step leads to different trees. Later on we also will need the following notation. For a given enhanced bivariate wavelet index $ \tilde{\boldsymbol{\lambda}} \in \tilde{\boldsymbol{\Lambda}}_{0}   $ and a bivariate wavelet tree $ \mathcal{T}  $ with $ \tilde{\boldsymbol{\lambda}} \in  \mathcal{T} $ we define the set $  \mathcal{C}(\tilde{\boldsymbol{\lambda}}, \mathcal{T}) $ which collects all direct children of $   \tilde{\boldsymbol{\lambda}} $ according to $\mathcal{T}$. For each $  \tilde{\boldsymbol{\lambda}} $ the set $  \mathcal{C}(\tilde{\boldsymbol{\lambda}}, \mathcal{T}) $ is either empty (if $ \tilde{\boldsymbol{\lambda}} $ is a leaf of $ \mathcal{T}  $) or can be described by one of the refinement strategies (LSR.a.1), (LSR.a.2), (LSR.b) or (LSR.c). If it is clear from the context, which bivariate wavelet tree $\mathcal{T}$ we use, sometimes we only write $  \mathcal{C}(\tilde{\boldsymbol{\lambda}}) $ instead of $  \mathcal{C}(\tilde{\boldsymbol{\lambda}}, \mathcal{T}) $. Moreover, when no wavelet tree is given, sometimes we use the notation $   \mathcal{C}(\tilde{\boldsymbol{\lambda}},  \clubsuit )   $ with $ \clubsuit \in \{ $ (LSR.a.1), (LSR.a.2), (LSR.b), (LSR.c) $  \}  $. Then $   \mathcal{C}(\tilde{\boldsymbol{\lambda}},  \clubsuit )   $ refers to all children of $ \tilde{\boldsymbol{\lambda}} $ according to one of the refinement options (LSR.a.1), (LSR.a.2), (LSR.b) or (LSR.c).

\subsection{Reference Rectangles and Associated Bivariate Wavelets }\label{sec_ref_rec_wavl}

In this section we will see that each reference rectangle and therefore also each bivariate wavelet index $  \boldsymbol{\lambda} = ( (j_{1},k_{1}), (j_{2},k_{2}) )  $ can be associated with a bivariate tensor wavelet and vice versa. Since for enhanced wavelet indices $\tilde{\boldsymbol{\lambda}} = ((j_{1},k_{1}), (j_{2},k_{2}), \alpha )$ the additional parameter $ \alpha  $ only describes the refinement options, the subsequent considerations are independent of $ \alpha  $ and can be done for pure wavelet indices $\boldsymbol{\lambda}$. Looking at the construction of our bivariate tensor wavelets, see Section \ref{subsec_biv_tensorwavelets}, we find that they are only defined for $ j_{1} \geq j_{0} -1 $ and $ j_{2} \geq j_{0} -1 $. On the other hand reference rectangles exist for all $ j_{1}, j_{2} \in \mathbb{N}_{0}    $. Moreover, for $ j_{1} = j_{0} - 1   $ or $ j_{2} = j_{0} - 1   $ we have to deal with tensors of the generator functions given in \eqref{eq_qua_gen_jjj}.  Hence we have to pay special attention to the cases $ j_{1} \leq j_{0} - 1  $ and $ j_{2} \leq j_{0} - 1    $. In order to take into account all these issues, in what follows a distinction of cases becomes necessary. 

\vspace{0,2 cm}

\underline{Case 1: $j_{1} \geq j_{0} + 1  $ and $j_{2} \geq j_{0} + 1 $}

\vspace{0,2 cm}

Let $j_{1}, j_{2} \geq j_{0} + 1 $, $ k_{1} \in \{ 0, 1, \ldots , 2^{j_{1}} - 1 \}   $ and  $ k_{2} \in \{ 0, 1, \ldots , 2^{j_{2}} - 1 \}   $. Then each reference rectangle $ R^{(j_{1},k_{1})}_{(j_{2},k_{2})}   $ is associated with a  bivariate tensor wavelet of the form 
\begin{align*}
\boldsymbol{\psi_{\lambda}} =  \boldsymbol{\psi}_{( (0,j_{1},k_{1}), (0,j_{2},k_{2}) ) } =  \psi_{(0,j_{1},k_{1})}^{b} \otimes \psi_{(0,j_{2},k_{2})}^{b} 
\end{align*}
and vice versa, see Section \ref{subsec_biv_tensorwavelets} for more details. The motivation for that reads as follows. For each bivariate wavelet index $  \boldsymbol{\lambda} = ( (j_{1},k_{1}), (j_{2},k_{2}) )  $ with $j_{1}, j_{2} \geq j_{0} + 1 $, $ k_{1} \in \{ 0, 1, \ldots , 2^{j_{1}} - 1 \}   $ and  $ k_{2} \in \{ 0, 1, \ldots , 2^{j_{2}} - 1 \}   $ there exist constants $ C_{1} , C_{2} > 0   $ independent of $ \boldsymbol{\lambda}   $ such that $ \supp      \psi_{(0,j_{1},k_{1})}^{b} \subset C_{1} [ 2^{-j_{1}} k_{1} , 2^{-j_{1}} ( k_{1} + 1 ) )  $ and $ \supp      \psi_{(0,j_{2},k_{2})}^{b} \subset C_{2} [ 2^{-j_{2}} k_{2} , 2^{-j_{2}} ( k_{2} + 1 ) )$. Consequently, there exists a constant $  C_{3} > 0 $ independent of $ \boldsymbol{\lambda}  $ such that
\begin{align*}
\supp \Big ( \psi_{(0,j_{1},k_{1})}^{b} \otimes \psi_{(0,j_{2},k_{2})}^{b} \Big ) \subset C_{3} R^{(j_{1},k_{1})}_{(j_{2},k_{2})}  .
\end{align*}
Let us remark, that in principle for the case $ j_{1} = j_{0} $ or $ j_{2} = j_{0}   $ a similar approach could be used. However, when we match the bivariate functions resulting out of the generators given in \eqref{eq_qua_gen_jjj} to reference rectangles, this will also have consequences for the level $j_{0}$. Therefore the cases $ j_{1} = j_{0} $ or $ j_{2} = j_{0}   $ will be treated separately below. 

\vspace{0,2 cm}

\underline{Case 2: $ j_{1} < j_{0}   $ and/or $ j_{2} < j_{0}   $}

\vspace{0,2 cm}

In Section \ref{Subsec_bound_qua} we only have constructed univariate boundary wavelets and corresponding generator functions for $ j \geq j_{0} - 1   $. Since our bivariate wavelets are defined via tensor products of these univariate functions, see \eqref{def_tensor_quarklet}, the wavelet system contained in $  \boldsymbol{\Psi}  $ only includes functions with  $ j_{1} \geq j_{0} - 1   $ and $ j_{2} \geq j_{0} - 1   $. Hence, if $  j_{1} < j_{0} - 1     $ or $ j_{2} < j_{0} - 1  $, we define 
\begin{equation}\label{eq_zerowav_put}
\boldsymbol{\psi}_{( (0,j_{1},k_{1}), (0,j_{2},k_{2}) ) } := 0 .
\end{equation}
A similar strategy also is used in the univariate setting, see Chapter 4.5 in \cite{VoDiss}. Thanks to \eqref{eq_zerowav_put} later on we will be able to work with bivariate wavelet and quarklet trees that have only one root. The special cases $ j_{1} = j_{0} - 1   $ and $ j_{2} = j_{0} - 1     $ are connected with  the bivariate functions resulting out of the generators given in \eqref{eq_qua_gen_jjj}. Below we will see, that they can be assigned to reference rectangles with $ j_{1} = j_{0}   $ or $ j_{2} = j_{0}     $. Consequently, also the reference rectangles with $ j_{1} = j_{0} - 1   $ or $ j_{2} = j_{0} - 1     $ can be associated with the zero function.

\vspace{0,2 cm}

\underline{Case 3: The special case $ j_{1} = j_{0}   $ and $ j_{2} \geq j_{0} + 1  $ (or $ j_{2} = j_{0}    $ and $ j_{1} \geq j_{0} + 1  $)}

\vspace{0,2 cm}

It remains to incorporate the bivariate functions resulting out of the generators constructed in Section \ref{Subsec_bound_qua} into the concept of reference rectangles. Here we have the difficulty that in the univariate setting there are $ |\Delta_{j_{0}}  |  =  2^{j_{0}} - 1 + m  $ functions we have to deal with. In most of the cases this number is not a power of two. To overcome this problem we have to invent a rule how each element of the index set $ \nabla_{j_{0} - 1} =  \Delta_{j_{0}}  $ can be assigned to an element of  $ \nabla_{j_{0}} = \{ 0, 1, \ldots , 2^{j_{0}} - 1 \}   $. In principle there are different possibilities how to reach this goal. One that is especially valuable and balanced is the following, see Chapter 4.5.1 and especially equation (4.37) in \cite{VoDiss}. For each $ k \in \nabla_{j_{0}-1}  $ we define the number 
\begin{equation}\label{eq_def_l_k_opt2}
	\ell_k := \left [ \frac{(2^{j_0}-1)(k+m-1)}{2^{j_0} +m  - 2} \right ].
\end{equation}
Here $ [ \cdot ]  $ denotes the nearest integer function. Now we assign $(j_0-1,k)$ to $(j_0,\ell_k)$. The idea behind \eqref{eq_def_l_k_opt2} is to map all $k \in \nabla_{j_0-1}$ to real numbers contained in the interval $[0,2^{j_0}-1]$ such that they are distributed uniformly and such that the supports of the corresponding wavelets and generator functions match roughly. The function $ [ \cdot ]  $ is applied to end up with an integer. For $\hat{k} \in \nabla_{j_{0}}$ we put
\begin{equation}\label{index_square_dissdor}
\square_{j_{0}, \hat{k}} := \{ k \in \nabla_{j_0-1} : \ell_k = \hat{k}     \} .
\end{equation}
This set refers to the translation parameters of the generator functions in one direction that are assigned to a pair $ ( j_{0} , \hat{k} )    $ referring to the lowest wavelet level in this direction. Now let $ j_{1} = j_{0}  $, $ j_{2} \geq j_{0} + 1  $, $ k_{1} \in \{ 0, 1, \ldots , 2^{j_{0}} - 1 \}   $ and  $ k_{2} \in \{ 0, 1, \ldots , 2^{j_{2}} - 1 \}   $. Then on the one hand each reference rectangle $ R^{(j_{1},k_{1})}_{(j_{2},k_{2})}   $ will be associated with the bivariate wavelet 
\begin{align*}
\boldsymbol{\psi}_{( (0,j_{0},k_{1}), (0,j_{2},k_{2}) ) } =  \psi_{(0,j_{0},k_{1})}^{b} \otimes \psi_{(0,j_{2},k_{2})}^{b} .
\end{align*}
On the other hand the same reference rectangle $ R^{(j_{1},k_{1})}_{(j_{2},k_{2})}   $ will be associated with all functions 
\begin{align*}
  \varphi_{(0,j_{0},\tilde{k})} \otimes \psi_{(0,j_{2},k_{2})}^{b} .
\end{align*}
Here $ \tilde{k} $ runs through each element of the set  $  \nabla_{j_0-1}   $ such that $ \ell_{\tilde{k}} =  k_{1}  $. In other words the reference rectangles with $ j_{1} = j_{0}   $ are connected with more than one function in order to incorporate the generator functions. The case $ j_{2} = j_{0}    $ can be treated with similar methods.

\subsection{Bivariate Quarklet Trees}\label{subsec_biv_qua_tre}

It is one main goal of this paper to approximate bivariate functions by using bivariate quarklet trees. For that purpose in what follows we generalize the concept of bivariate wavelet indices to the more advanced concept of bivariate quarklet indices. Let $   j_{1} , j_{2} \in \mathbb{N}_{0},  k_{1} \in \{ 0, 1, \ldots , 2^{j_{1}} - 1 \} $ and $ k_{2} \in \{ 0, 1, \ldots , 2^{j_{2}} - 1 \}    $. Then to each bivariate wavelet index $  ( (j_{1},k_{1}), (j_{2},k_{2}) )   $ we match additional parameters $ p_{1}, p_{2} \in \mathbb{N}_{0}     $ in order to obtain bivariate quarklet indices $ \boldsymbol{\lambda} := ( ( p_{1}, j_{1},k_{1}), ( p_{2},  j_{2},k_{2}) ) $. Each bivariate quarklet index refers to a bivariate tensor quarklet. To see this recall the previous Section \ref{sec_ref_rec_wavl} where we have found that each wavelet index $  ( (j_{1},k_{1}), (j_{2},k_{2}) ) =   ( (0, j_{1},k_{1}), (0, j_{2},k_{2}) )     $ refers to a bivariate function 
\begin{align*}
  \psi_{(0,j_{1},k_{1})}^{b} \otimes \psi_{(0,j_{2},k_{2})}^{b} ,
\end{align*}
with modifications for $ j_{1} \leq j_{0}  $ and/or $ j_{2} \leq j_{0}  $. Hence, for $j_{1}, j_{2} \geq j_{0} + 1  $, $  k_{1} \in \{ 0, 1, \ldots , 2^{j_{1}} - 1 \} $, $ k_{2} \in \{ 0, 1, \ldots , 2^{j_{2}} - 1 \}    $ and  $ p_{1}, p_{2} \in \mathbb{N}_{0}     $ each bivariate quarklet index $   ( ( p_{1}, j_{1},k_{1}), ( p_{2},  j_{2},k_{2}) )    $ can be associated with a bivariate tensor quarklet  
\begin{align*}
  \psi_{(p_{1},j_{1},k_{1})}^{b} \otimes \psi_{(p_{2},j_{2},k_{2})}^{b} .
\end{align*} 
In the case $ j_{1} < j_{0}  $ and/or $ j_{2} < j_{0}  $ following Section \ref{sec_ref_rec_wavl} we define $  \boldsymbol{\psi}_{( (p_{1},j_{1},k_{1}), (p_{2},j_{2},k_{2}) ) } := 0    $. For $ j_{1} = j_{0}  $ and $ j_{2} > j_{0}  $ as described above the quarklet index $   ( ( p_{1}, j_{0},k_{1}), ( p_{2},  j_{2},k_{2}) ) $ refers to the bivariate tensor quarklet 
\begin{align*}
  \psi_{(p_{1},j_{0},k_{1})}^{b} \otimes \psi_{(p_{2},j_{2},k_{2})}^{b} \qquad \mbox{and to all functions} \qquad \varphi_{(p_{1},j_{0},\tilde{k})} \otimes \psi_{(p_{2},j_{2},k_{2})}^{b} .
\end{align*} 
Again  $ \tilde{k} $ runs through each element of the set  $  \nabla_{j_0-1}   $ such that $ \ell_{\tilde{k}} =  k_{1}  $, see \eqref{eq_def_l_k_opt2}. Here $ \varphi_{(p_{1},j_{0},\tilde{k})}   $ refers to the Schoenberg B-spline quarks given in Definition \ref{def_bound_quark}. The converse case $ j_{1} > j_{0}  $ and $ j_{2} = j_{0}  $ can be treated with similar methods. In order to collect all bivariate quarklet indices we introduce the index set 
\begin{align*}
\boldsymbol{\Lambda} := & \Big  \{ ( (p_{1},j_{1},k_{1}), (p_{2},j_{2},k_{2}) ) : p_{1}, p_{2},  j_{1} , j_{2} \in \mathbb{N}_{0},  k_{1} \in \{ 0, 1, \ldots , 2^{j_{1}} - 1 \} , \\
& k_{2} \in \{ 0, 1, \ldots , 2^{j_{2}} - 1 \}  \Big   \} .
\end{align*}
To gather the enhanced bivariate quarklet indices, which also contain the refinement options of the corresponding wavelet indices, we define
\begin{align*}
\tilde{\boldsymbol{\Lambda}} := &  \Big  \{ ( (p_{1},j_{1},k_{1}), (p_{2},j_{2},k_{2}), \alpha ) : p_{1}, p_{2},  j_{1} , j_{2} \in \mathbb{N}_{0},  k_{1} \in \{ 0, 1, \ldots , 2^{j_{1}} - 1 \}  ,  \\
& k_{2} \in \{ 0, 1, \ldots , 2^{j_{2}} - 1 \}, \alpha \in \{ 0 , 1 , 2   \}  \Big   \} .
\end{align*}
For a given (enhanced) quarklet index $ \tilde{\boldsymbol{\lambda}} = ( (p_{1},j_{1},k_{1}), (p_{2},j_{2},k_{2}), \alpha ) \in \tilde{\boldsymbol{\Lambda}}   $ we use the notation $ |  \tilde{\boldsymbol{\lambda}}   | :=  j_{1} + j_{2}  $. Furthermore we introduce the mapping $\circ : \tilde{\boldsymbol{\Lambda}} \rightarrow \tilde{\boldsymbol{\Lambda}}_0 $ defined by
\begin{equation*}\label{eq:circ}
\tilde{\boldsymbol{\lambda}} \mapsto \tilde{\boldsymbol{\lambda}}^\circ \coloneqq ( (p_{1},j_{1},k_{1}), (p_{2},j_{2},k_{2}), \alpha )^\circ \coloneqq   ( (0,j_{1},k_{1}), (0,j_{2},k_{2}), \alpha )  ,
\end{equation*}
which provides the corresponding bivariate wavelet index for a given bivariate quarklet index. Recall, that we can identify enhanced bivariate wavelet indices  $  ( (j_{1},k_{1}), (j_{2},k_{2}), \alpha ) $ with enhanced bivariate quarklet indices $ ( (0,j_{1},k_{1}), (0,j_{2},k_{2}), \alpha )   $. Consequently we have $ \tilde{\boldsymbol{\Lambda}}_{0} \subset    \tilde{\boldsymbol{\Lambda}} $. Let us consider a bivariate wavelet tree $\mathcal{T} \subset \tilde{\boldsymbol{\Lambda}}_0 \subset \tilde{\boldsymbol{\Lambda}}$ as a set of bivariate quarklet indices. Then there are different options for the refinement of a leaf $\tilde{\boldsymbol{\lambda}} \in \mathcal{V}(\mathcal{T})$.

\vspace{0,2 cm}

\underline{Option 1: Space refinement}

\vspace{0,2 cm}

For a given leaf $ ( (0,j_{1},k_{1}), (0,j_{2},k_{2}), \alpha ) =  \tilde{\boldsymbol{\lambda}} \in \mathcal{V}(\mathcal{T})    $ we can refine in space by using one of the refinement strategies (LSR.a.1), (LSR.a.2), (LSR.b) or (LSR.c) and add the corresponding child wavelet indices to it. Due to Definition \ref{def_tre_terms} this is always possible.

\vspace{0,2 cm}

\underline{Option 2: Increase the polynomial degree}

\vspace{0,2 cm}

For a given leaf $ ( (0,j_{1},k_{1}), (0,j_{2},k_{2}), \alpha ) =  \tilde{\boldsymbol{\lambda}} \in \mathcal{V}(\mathcal{T})    $ it is also possible to increase the polynomial degree by adding certain quarklet indices with $  p_{1} > 0    $ or/and $  p_{2} > 0    $. Here we have different possibilities due to the different Cartesian directions. However, we require a heuristic which tells us which bivariate tensor quarklets should be added. To this end we consider sets $\Upsilon(\tilde{\boldsymbol{\lambda}}) \subset \tilde{\boldsymbol{\Lambda}}_0$ such that for each bivariate wavelet tree $\mathcal{T} \subset \tilde{\boldsymbol{\Lambda}}_0$ the union of the sets $\Upsilon(\tilde{\boldsymbol{\lambda}})$ over all leaves $\tilde{\boldsymbol{\lambda}} \in \mathcal{V}(\mathcal{T})$ provides a disjoint decomposition of $\mathcal{T}$. In order to define such sets $\Upsilon(\tilde{\boldsymbol{\lambda}})$ for each of the refinement strategies (LSR.a.1), (LSR.a.2), (LSR.b) and (LSR.c) we determine a so-called chosen child. There are different possibilities how this can be done. One, which we will always use in our later considerations, is the following.

\begin{itemize}
\item In strategy (LSR.a.1) the chosen child is $   ( (j_{1},k_{1}), (j_{2},k_{2}), 2 ) $.

\item In strategy (LSR.a.2) the chosen child is $ ( (j_{1} , k_{1}), (j_{2} + 1 , 2 k_{2}), 0 )  $.

\item In strategy (LSR.b) the chosen child is $  ( (j_{1} + 1, 2 k_{1}), (j_{2},k_{2}), 0 )   $.

\item In strategy (LSR.c) the chosen child is $ ( (j_{1} , k_{1}), (j_{2} + 1 , 2 k_{2}), 0 )   $.

\end{itemize}

With other words, whenever possible the chosen child refers to a refinement in direction $i = 2$ with even $k_{2} \in \{ 0, 1, \ldots , 2^{j_{2} + 1} - 1 \}  $. In strategy (LSR.a.1) such a child does not exist. Therefore we pick $   ( (j_{1},k_{1}), (j_{2},k_{2}), 2 ) $ as chosen child since it allows us a space refinement in direction $i = 2$ later on. Strategy (LSR.b) is completely devoted to a refinement in direction $i = 1$. Hence in this case also the chosen child has to refer to a refinement in the first Cartesian direction and we select $  ( (j_{1} + 1, 2 k_{1}), (j_{2},k_{2}), 0 )   $. Now we can use the concept of chosen children to define the sets $ \Upsilon(\tilde{\boldsymbol{\lambda}}) $. For that purpose let a wavelet tree $\mathcal{T} \subset \tilde{\boldsymbol{\Lambda}}_0 $ be given. Then for each leave $ \tilde{\boldsymbol{\lambda}} \in \mathcal{V}(\mathcal{T}) $ we can determine the sets $ \Upsilon(\tilde{\boldsymbol{\lambda}}) $ by using the following algorithm.

\vspace{0,4 cm}

\medskip
\noindent\shadowbox{\parbox{0.96\textwidth}{
		\begin{algorithm*} {\bf CREATE\textunderscore $\Upsilon$} $[\mathcal{T}, \tilde{\boldsymbol{\lambda}} \in \mathcal{V}(\mathcal{T})] \mapsto \Upsilon(\tilde{\boldsymbol{\lambda}})$\\[-3ex]
			\noindent{
				\begin{tabbing}
					set $\Upsilon(\tilde{\boldsymbol{\lambda}}) = \{  \tilde{\boldsymbol{\lambda}}  \}$, $\tilde{\boldsymbol{\mu}} := \tilde{\boldsymbol{\lambda}} $ and $ c = 0 $;\\[0.1cm]
					\texttt{\em whi}\=\texttt{\em le} $c = 0$ \\[0.1cm]
					\> take direct ancestor $ \tilde{\boldsymbol{\eta}} \in  \mathcal{T} $ of $ \tilde{\boldsymbol{\mu}}  $; \\[0.1cm]
					\> \texttt{\em if}\={} $ \tilde{\boldsymbol{\mu}} $ is chosen child of $  \tilde{\boldsymbol{\eta}} $\\[0.1cm]
					\> \> add $ \tilde{\boldsymbol{\eta}}  $ to $  \Upsilon(\tilde{\boldsymbol{\lambda}}) $ and put $\tilde{\boldsymbol{\mu}} := \tilde{\boldsymbol{\eta}}$;\\[0.1cm]
					\> \texttt{\em else}\\[0.1cm]
					\> \> put $ c = 1$;\\[0.1cm]
					\> \texttt{\em end if}\\[0.1cm]
					\texttt{\em end while}
				\end{tabbing}
			}
		\end{algorithm*}
}}
\medskip

\vspace{0,4 cm}

When we apply {\bf CREATE\textunderscore $\Upsilon$} for all leaves $ \tilde{\boldsymbol{\lambda}} \in \mathcal{V}(\mathcal{T})$ of a given tree $  \mathcal{T} $ the resulting sets $ \Upsilon(\tilde{\boldsymbol{\lambda}})  $ have the following properties:

\begin{itemize}
	\item[(i)] The sets $\Upsilon(\tilde{\boldsymbol{\lambda}})$ have the form 
	\begin{equation}\label{eq:Upsilon}
\Upsilon(\tilde{\boldsymbol{\lambda}})= \{ \tilde{\boldsymbol{\mu}} \in \mathcal{T} : \tilde{\boldsymbol{\lambda}} \succeq  \tilde{\boldsymbol{\mu}} \succeq \boldsymbol{\tilde{\mu}_{\tilde{\lambda}}} \}
	\end{equation} 
	for some fixed $ \boldsymbol{\tilde{\mu}_{\tilde{\lambda}}} \preceq \tilde{\boldsymbol{\lambda}}$ with  $ \boldsymbol{\tilde{\mu}_{\tilde{\lambda}}} \in \mathcal{T}   $. 
    \item[(ii)] For each tree $\mathcal{T}$ and inner node $\tilde{\boldsymbol{\mu}} \in \mathcal{T}\backslash \mathcal{V}(\mathcal{T})$ with children $\tilde{\boldsymbol{\eta}}_1,\tilde{\boldsymbol{\eta}}_2, \tilde{\boldsymbol{\eta}}_3  \in \mathcal{T}$ it holds $\Upsilon(\tilde{\boldsymbol{\eta}}_1) \cap  \Upsilon(\tilde{\boldsymbol{\eta}}_2) \cap  \Upsilon(\tilde{\boldsymbol{\eta}}_3) = \emptyset$. The same holds if there are only two children.
    \item[(iii)] For each tree $\mathcal{T}$ it holds $\bigcup_{\tilde{\boldsymbol{\lambda}} \in \mathcal{V}(\mathcal{T})} \Upsilon(\tilde{\boldsymbol{\lambda}}) = \mathcal{T}$.
\end{itemize}

\begin{Remark}
It seems to be possible to choose alternative definitions for the sets $ \Upsilon(\tilde{\boldsymbol{\lambda}})  $. For example, when selecting the chosen children we also can favor refinements in direction $i = 1$. And also other selection procedures are conceivable. However, the value of the employed definition for the sets $ \Upsilon(\tilde{\boldsymbol{\lambda}})  $ also depends on the test function which should be approximated. The selection process presented in our algorithm {\bf CREATE\textunderscore $\Upsilon$} is well-suited for the test case given in Section \ref{Sec_aniso_example} below. Therefore we stick with this definition in what follows.  
\end{Remark}

Now we can use the sets $\Upsilon(\tilde{\boldsymbol{\lambda}})$ to introduce the polynomial enrichment of a leaf. We increase the maximal polynomial degree of $\tilde{\boldsymbol{\lambda}} \in \mathcal{V}(\mathcal{T})$ by adding bivariate quarklet indices with either $p_{1} = 1$ or $p_{2} = 1$  to each node $\tilde{\boldsymbol{\mu}} \in \Upsilon(\tilde{\boldsymbol{\lambda}})$. Again there are different possibilities for increasing the polynomial degree due to the different Cartesian directions. So for the first polynomial enrichment of a leaf $ \tilde{\boldsymbol{\lambda}}    $ by $1$ we can put
\begin{equation*}\label{eq:p-refinement_d1}
\mathcal{T} \cup  \bigcup_{\tilde{\boldsymbol{\mu}} = ((i_{1},\ell_{1}),(i_{2},\ell_{2}), \alpha) \in \Upsilon(\tilde{\boldsymbol{\lambda}})} \bigcup_{\substack{p_{1}, p_{2} \in \mathbb{N}_{0} \\ 0 < p_{1} + p_{2} \leq p_{\max} = 1}} ((p_{1}, i_{1}, \ell_{1}),(p_{2} , i_{2}, \ell_{2}), \alpha). 
\end{equation*}
In a next step the process of polynomial enrichment can be repeated with a different leaf or with the same leaf and the next higher polynomial degree $p_{\max} = 2$ (and subsequently also $p_{\max} = 3, 4, 5, \ldots $). Now we are well-prepared to define bivariate quarklet trees. For that purpose at first we require some additional notation. Let $T \subset \tilde{\boldsymbol{\Lambda}} $ be a set of enhanced bivariate quarklet indices. Then by $T^\circ$ we denote the corresponding set of enhanced bivariate wavelet indices, namely
\begin{align*}
T^\circ := \{ \tilde{\boldsymbol{\lambda}}^\circ \in \tilde{\boldsymbol{\Lambda}}_{0} :  \tilde{\boldsymbol{\lambda}} \in T    \} .
\end{align*}
The definition of bivariate quarklet trees reads as follows. 

\begin{Definition}\label{def:quarklet_tree}
Let $T \subset \tilde{\boldsymbol{\Lambda}}$ be a set of enhanced bivariate quarklet indices. For all $\tilde{\boldsymbol{\lambda}}^\circ = ( (0,j_{1},k_{1}), (0,j_{2},k_{2}), \alpha ) \in T^{\circ}$ we put
\begin{align*} 
p_{\max}(\tilde{\boldsymbol{\lambda}}^\circ) \coloneqq p_{\max}(\tilde{\boldsymbol{\lambda}}^\circ,T) \coloneqq \max\{ p_{1} + p_{2} \in \mathbb{N}_0 : ( (p_{1},j_{1},k_{1}), (p_{2},j_{2},k_{2}), \alpha ) \in T\} .
	 \end{align*}
Then $  T $  is called a \emph{bivariate quarklet tree} if the following conditions are fulfilled:
	\begin{itemize}
		\item[(i)]
		The corresponding set $T^\circ \subset \tilde{\boldsymbol{\Lambda}}_{0}$ is a bivariate wavelet tree according to Definition \ref{def_wav_tre}. 
		\item[(ii)] For each $\tilde{\boldsymbol{\lambda}}^\circ \in \mathcal{V}(T^{\circ})$ we have $p_{\max}(\tilde{\boldsymbol{\lambda}}^\circ) = p_{\max}(\tilde{\boldsymbol{\mu}}^\circ)$ for all  $\tilde{\boldsymbol{\mu}}^\circ \in \Upsilon(\tilde{\boldsymbol{\lambda}}^\circ)$. 
		\item[(iii)] For each $\tilde{\boldsymbol{\lambda}}^\circ = ( (0,j_{1},k_{1}), (0,j_{2},k_{2}), \alpha ) \in T^{\circ}$ we have $  ( (p_{1},j_{1},k_{1}), (p_{2},j_{2},k_{2}), \alpha )   \in T$ for all $p_{1}, p_{2} \in \mathbb{N}_{0}$ with $0 < p_{1} + p_{2} \le p_{\max}(\tilde{\boldsymbol{\lambda}}^\circ)$. 
\end{itemize}
\end{Definition}
In other words a bivariate quarklet tree consists of an underlying bivariate wavelet index set possessing a tree structure and moreover the nodes of this tree are enriched with all enhanced bivariate quarklet indices up to a certain polynomial degree $p_{\max}$. When we talk about a leaf, node or root of a quarklet tree, we always mean the corresponding wavelet index, which is guaranteed to be an element of the tree and can be accessed from a suitable enhanced bivariate quarklet index via the mapping $\circ$.  By $|T|$ we denote the number of enhanced bivariate wavelet indices in an (arbitrary) index set $T \subseteq \tilde{\boldsymbol{\Lambda}} $. For a bivariate quarklet tree $T$ we set its cardinality to be the number of quarklet indices in the tree, namely
\begin{equation}\label{eq_card_tree11}
\#T \coloneqq |T| +  \sum_{\tilde{\boldsymbol{\lambda}}^\circ \in T^{\circ}} \Big ( \frac{( p_{\textrm{max}}(\tilde{\boldsymbol{\lambda}}^\circ) + 1)^2 + (  p_{\textrm{max}}(\tilde{\boldsymbol{\lambda}}^\circ) + 1  )}{2} - 1 \Big ) .
\end{equation}
There are two different ways to characterize a bivariate quarklet tree $T$. The first way is to consider a bivariate wavelet tree $\mathcal{T}$ and then fix the maximal polynomial degrees $p_{\max}$ on all leaves. To this end we write
\begin{align*}
P_{\max} \coloneqq \big\{ ( p_{\textrm{max}}(\tilde{\boldsymbol{\lambda}}^{\circ}) )  \big\}_{\tilde{\boldsymbol{\lambda}}^{\circ} \in \mathcal{V}(\mathcal{T})}.
\end{align*}
Then the maximal polynomial degrees on the inner nodes can be determined by using Definition \ref{def:quarklet_tree}. For all $\tilde{\boldsymbol{\lambda}}^{\circ} \in \mathcal{V}(\mathcal{T})$ and $\tilde{\boldsymbol{\mu}}^{\circ} \in \Upsilon(\tilde{\boldsymbol{\lambda}}^{\circ})$ we have to set $p_{\max}(\tilde{\boldsymbol{\mu}}^{\circ}) = p_{\max}(\tilde{\boldsymbol{\lambda}}^{\circ})$ to end up with a bivariate quarklet tree $T$. Therefore the assignment $P_{\max}$ already implies the maximal polynomial degrees on all nodes (and not just on the leaves) and we can write $T = (\mathcal{T}, P_{\max})$ since this notation contains all information to establish a bivariate quarklet tree.

For the second option we consider two bivariate wavelet trees $ \mathcal{T}  $ and $ \mathcal{T}'   $ with $ \mathcal{T} \subset \mathcal{T}'    $. Let $  \tilde{\boldsymbol{\lambda}}^{\circ} \in \mathcal{V}(\mathcal{T})  $ and investigate the set $ ( \mathcal{T}' \setminus \mathcal{T}  ) \cup \{ \tilde{\boldsymbol{\lambda}}^{\circ}  \} $. Let $ R( \mathcal{T} , \mathcal{T}' , \tilde{\boldsymbol{\lambda}}^{\circ})    $ be the largest subset of $ ( \mathcal{T}' \setminus \mathcal{T}  ) \cup \{ \tilde{\boldsymbol{\lambda}}^{\circ}  \}   $ that can be obtained by a sequence of space refinements of the form (LSR.a.1), (LSR.a.2), (LSR.b) and (LSR.c) starting at $  \tilde{\boldsymbol{\lambda}}^{\circ} $ and using an iterative process, whereby in each step (UPC) is fulfilled. Due to Definition \ref{def_wav_tre} such a set always exists. Let $ r(\mathcal{T} , \mathcal{T}' , \tilde{\boldsymbol{\lambda}}^{\circ}) \in \mathbb{N}_{0}  $ be the number of space refinements of the form (LSR.a.1), (LSR.a.2), (LSR.b) and (LSR.c) that are used to obtain the set $ R( \mathcal{T} , \mathcal{T}' , \tilde{\boldsymbol{\lambda}}^{\circ}) $. Now for each $   \tilde{\boldsymbol{\lambda}}^{\circ} \in \mathcal{V}(\mathcal{T})$ we can put $ p_{\max}(\tilde{\boldsymbol{\lambda}}^{\circ}) := r(\mathcal{T} , \mathcal{T}' , \tilde{\boldsymbol{\lambda}}^{\circ})    $. Using Definition \ref{def:quarklet_tree} this implies a quarklet tree $  T = (\mathcal{T}, P_{\max})   $. Consequently we can also write $ T = (\mathcal{T}, P_{\max}) = ( \mathcal{T} , \mathcal{T}'   )   $ since $ P_{\max}  $ is given by $ \mathcal{T}    $ and $ \mathcal{T}'    $. The intuition behind this is that we delete the descendants of $ \tilde{\boldsymbol{\lambda}}^{\circ} \in \mathcal{V}(\mathcal{T})$ and instead employ polynomial enrichment. This process is called trimming and will be used in our adaptive scheme later on.

\subsection{Local Errors, Global Errors and Best Approximation}\label{subsec_localerr}
It is one of the main goals of this paper to construct an adaptive bivariate quarklet algorithm to approximate given bivariate functions $ f \in L_{2}((0,1)^2)   $ in an efficient way. For that purpose in what follows we have to introduce some error functionals. In a first step for each bivariate wavelet index $ \tilde{\boldsymbol{\lambda}} = ( (j_{1},k_{1}), (j_{2},k_{2}), \alpha ) \in \tilde{\boldsymbol{\Lambda}}_0   $ with associated maximal polynomial degree $p_{\max}( \tilde{\boldsymbol{\lambda}} ) \in \mathbb{N}_0$ we investigate local errors $ e_{p_{\max}(\tilde{\boldsymbol{\lambda}})}(\tilde{\boldsymbol{\lambda}}) : \tilde{\boldsymbol{\Lambda}}_0 \rightarrow [0, \infty )    $. They are supposed to satisfy the following two very important properties:  
\begin{itemize}
	\item[(i)] There is a subadditivity for the error of the lowest order. That means for $\tilde{\boldsymbol{\lambda}} \in \tilde{\boldsymbol{\Lambda}}_0$ with children  $\tilde{\boldsymbol{\eta}} \in \mathcal{C}(\tilde{\boldsymbol{\lambda}}, \clubsuit )$ we require
	\begin{equation}\label{eq:1.2}
	e_0(\tilde{\boldsymbol{\lambda}}) \ge \sum_{\tilde{\boldsymbol{\eta}} \in \mathcal{C}(\tilde{\boldsymbol{\lambda}} , \clubsuit)} e_0(\tilde{\boldsymbol{\eta}}) 
	\end{equation}
simultaneously for all $ \clubsuit \in \{ $ (LSR.a.1), (LSR.a.2), (LSR.b), (LSR.c) $  \}  $. 	
	
	\item[(ii)] The error is reduced by increasing the maximal polynomial degree. Namely we have
	\begin{equation} \label{eq:1.3}
	e_{p_{\max}}(\tilde{\boldsymbol{\lambda}}) \ge e_{p_{\max}+1}(\tilde{\boldsymbol{\lambda}}). 
	\end{equation}
\end{itemize}
\begin{Remark}
The local errors $ e_{p_{\max}(\tilde{\boldsymbol{\lambda}})}(\tilde{\boldsymbol{\lambda}})   $ are given in a quite general way. However, in what follows they will be used to provide an adaptive algorithm which allows for bivariate quarklet tree approximation, whereby the output is near-best in the sense of Theorem \ref{theorem:1}. Due to the generality of $ e_{p_{\max}(\tilde{\boldsymbol{\lambda}})}(\tilde{\boldsymbol{\lambda}})   $ the results obtained below can be applied to approximate a broad class of functions. In some cases it can be useful to state the local errors in a more precise fashion. For example, if we want to approximate a function $ f \in  L_2((0,1)^2) $ such that the approximation error is given in terms of the norm  $ \Vert \cdot \vert   L_2((0,1)^2)    \Vert    $, then we can use Definition \ref{def_locerr_L2_new} to define the local errors.  
\end{Remark}
The local errors can be used to define a global error. For a given bivariate quarklet tree $T = (\mathcal{T}, P_{\max})$ we define the global error $\mathcal{E}(T)$ by
\begin{equation}\label{eq_glob_err}
\mathcal{E}(T) \coloneqq \sum_{\tilde{\boldsymbol{\lambda}}^{\circ} \in \mathcal{V}(\mathcal{T})} e_{p_{\max}(\tilde{\boldsymbol{\lambda}}^{\circ})}(\tilde{\boldsymbol{\lambda}}^{\circ}).
\end{equation}
It collects the local errors for all leaves of the tree.  The global error can be used to define the so-called best approximation error. 
\begin{Definition}\label{def_bestappr_err}
The \emph{error of the best bivariate quarklet tree approximation} of cardinality $n \in \mathbb{N}$ is defined by
\begin{align*}
\sigma_n \coloneqq \inf_{T = (\mathcal{T}, P_{\max} )} \ \inf_{\# T \le n} \mathcal{E}(T).
\end{align*}
\end{Definition}
Below we will find an incremental algorithm that for each $N \in \mathbb{N}$ produces a bivariate quarklet tree $  T_N = (\mathcal{T}_N, P_{\max}) $ with $  \# T_N \le \tilde{C} N^{3} $  that provides a near-best quarklet approximation in the sense of
\begin{equation}\label{eq:nearbest}
\mathcal{E}(T_N) \le C \sigma_{c N},
\end{equation}
with independent constants $C \ge 1 $ and $c \in (0,1]$, see Theorem \ref{theorem:1} for the details. 

\section{Adaptive Refinement Strategy} \label{sec:adap_ref}

\subsection{Error Functionals for Adaptive Refinement}\label{Subsec_Mod_Errfunc}
To construct our bivariate near-best quarklet algorithm we need some more error functionals, which will be introduced in the following section. Most of them trace back to the ideas of Binev, see \cite{Bin18}, and also have counterparts for the case of univariate quarklet tree approximation as described in \cite{DaHoRaVo}, see Section 3.1. Below we present three kinds of error functionals. At first we introduce a penalized version for the local errors of the lowest order that can be used to design near-best space adaptive schemes. Second we establish an error functional for the space and polynomial degree adaptive case. Finally we provide two indicators which help to decide where a refinement can be done in the next step of the algorithm.  \\

\vspace{0,2 cm}

\noindent
\underline{Step 1: A penalized version for the local error of the lowest order.} 

\vspace{0,2 cm}

\noindent
Let an enhanced bivariate wavelet index $ \tilde{\boldsymbol{\lambda}} = ( (j_{1},k_{1}), (j_{2},k_{2}), \alpha ) \in \tilde{\boldsymbol{\Lambda}}_0   $ be given. For that we define a modified local error functional denoted by $\tilde{e}(\tilde{\boldsymbol{\lambda}})$ with $ \tilde{e}(\tilde{\boldsymbol{\lambda}}) : \tilde{\boldsymbol{\Lambda}}_0 \rightarrow [0, \infty )    $, which is strongly connected with the local error of the lowest order $ e_{0}(\tilde{\boldsymbol{\lambda}})  $. Below for the sake of convenience sometimes we use $e(\tilde{\boldsymbol{\lambda}}) := e_0(\tilde{\boldsymbol{\lambda}})$. Let a bivariate wavelet tree $ \mathcal{T}  $ be given. Let $ \mathcal{R}  $ be the root of $ \mathcal{T}  $ and $\tilde{\boldsymbol{\mu}} \in  \mathcal{T}   $ be the parent of $\tilde{\boldsymbol{\lambda}} \in  \mathcal{T}   $. Then we define the modified local errors $ \tilde{e}  $ step by step via
\begin{equation}\label{eq:2.3}
\tilde{e}(\mathcal{R})\coloneqq e(\mathcal{R}), \qquad \tilde{e}(\tilde{\boldsymbol{\lambda}}) \coloneqq \frac{e(\tilde{\boldsymbol{\lambda}})\tilde{e}(\tilde{\boldsymbol{\mu}})}{e(\tilde{\boldsymbol{\lambda}})+\tilde{e}(\tilde{\boldsymbol{\mu}})}.
\end{equation}
In the case $e(\tilde{\boldsymbol{\lambda}})=\tilde{e}(\tilde{\boldsymbol{\mu}})=0$ we set $\tilde{e}(\tilde{\boldsymbol{\lambda}}) :=0 $. Equation \eqref{eq:2.3} implies
\begin{equation}\label{eq:2.7}
\frac{1}{\tilde{e}(\tilde{\boldsymbol{\lambda}})} = \frac{1}{e(\tilde{\boldsymbol{\lambda}})}+\frac{1}{\tilde{e}(\tilde{\boldsymbol{\mu}})} \qquad \mbox{and by iteration also} \qquad \frac{1}{\tilde{e}(\tilde{\boldsymbol{\lambda}})} = \sum_{\tilde{\boldsymbol{\mu}} \in \mathcal{A}_{\mathcal{T}}(\tilde{\boldsymbol{\lambda}})} \frac{1}{e(\tilde{\boldsymbol{\mu}})}.
\end{equation}
Here $ \mathcal{A}_{\mathcal{T}}(\tilde{\boldsymbol{\lambda}}) $ refers to the set of all ancestors of $ \tilde{\boldsymbol{\lambda}} $ (including $ \tilde{\boldsymbol{\lambda}} $ itself) according to the given bivariate wavelet tree $ \mathcal{T}  $.

\vspace{0,2 cm}

\noindent
\underline{Step 2: Local errors concerning space refinement and polynomial enrichment.} 

\vspace{0,2 cm}

\noindent
Let a bivariate wavelet tree $ \mathcal{T}  $ and an enhanced bivariate wavelet index $ \tilde{\boldsymbol{\lambda}} = ( (j_{1},k_{1}), (j_{2},k_{2}), \alpha ) \in \mathcal{T} $ be given. Then we introduce an error functional $ E(\tilde{\boldsymbol{\lambda}}) \coloneqq E(\tilde{\boldsymbol{\lambda}}, \mathcal{T})  $ with $ E(\tilde{\boldsymbol{\lambda}}) : \mathcal{T} \rightarrow [0 , \infty )  $. It is defined recursively starting at the leaves of the tree $\mathcal{T}$. Here for $ \tilde{\boldsymbol{\lambda}} \in \mathcal{V}(\mathcal{T})  $ we define $ E(\tilde{\boldsymbol{\lambda}})\coloneqq e(\tilde{\boldsymbol{\lambda}}) = e_0(\tilde{\boldsymbol{\lambda}})$. For the inner nodes of the tree the error functional is defined step by step moving from the leaves towards the root. Let $\tilde{\boldsymbol{\lambda}} \in \mathcal{T} \setminus \mathcal{V}(\mathcal{T})  $ and assume that $E(\tilde{\boldsymbol{\eta}})$ for all $ \tilde{\boldsymbol{\eta}} \in \mathcal{C}(\tilde{\boldsymbol{\lambda}} , \mathcal{T} )  $ are already known. Moreover, let $ r( \mathcal{T}, \tilde{\boldsymbol{\lambda}} ) \in \mathbb{N}_{0}   $ be the number of refinement steps of the form (LSR.a.1), (LSR.a.2), (LSR.b) and (LSR.c) that are required to obtain the greatest possible subtree of $ \mathcal{T}  $ starting at the root $  \tilde{\boldsymbol{\lambda}}  $.  Then for an inner node $ \tilde{\boldsymbol{\lambda}} \in \mathcal{T} \setminus \mathcal{V}(\mathcal{T}) $ we put    
\begin{equation}\label{eq:3.1}
E(\tilde{\boldsymbol{\lambda}}) \coloneqq \min \Big \{ \sum_{\tilde{\boldsymbol{\eta}} \in \mathcal{C}(\tilde{\boldsymbol{\lambda}} , \mathcal{T} ) } E(\tilde{\boldsymbol{\eta}}) , e_{r(\mathcal{T} , \tilde{\boldsymbol{\lambda}})}(\tilde{\boldsymbol{\lambda}}) \Big  \}.
\end{equation}
This refers to the adaptive choice between the two refinement types. Next we want to introduce a modified version of the error functional $E(\tilde{\boldsymbol{\lambda}})$. Therefore we observe that enlarging the tree $\mathcal{T}$ changes the quantity $E(\tilde{\boldsymbol{\lambda}}) = E(\tilde{\boldsymbol{\lambda}}, \mathcal{T})$ only if $ r(\mathcal{T} , \tilde{\boldsymbol{\lambda}})  $ changes. This is a direct consequence of \eqref{eq:3.1}. We use this observation and consider a sequence $\mathcal{T}_1,\mathcal{T}_2,\mathcal{T}_3,\ldots$ of growing trees. With that we mean that each tree $\mathcal{T}_{k+1}$ is derived from $\mathcal{T}_k$  by subdividing a leaf and adding two or three child indices according to the refinement options (LSR.a.1), (LSR.a.2), (LSR.b) or (LSR.c) to it. For a node $\tilde{\boldsymbol{\lambda}}$ and $j \in \mathbb{N}_0$ there might exist multiple trees $\mathcal{T}_\star$ with $\tilde{\boldsymbol{\lambda}} \in \mathcal{T}_\star$ and $r( \mathcal{T}_\star , \tilde{\boldsymbol{\lambda}})=j$ in the sequence $\mathcal{T}_1,\mathcal{T}_2,\mathcal{T}_3,\ldots$ of trees, since there is the possibility to carry out refinement steps in other parts of the tree. This means that the subtree emanating from $\tilde{\boldsymbol{\lambda}}$ stays the same in all the trees $\mathcal{T}_\star$ and consequently the quantity $E(\tilde{\boldsymbol{\lambda}},\mathcal{T}_\star)$ does not change. By using this observation we can let $ j \in \mathbb{N}_0  $ and $\mathcal{T}_\star$ be any of the trees in the sequence $\mathcal{T}_1,\mathcal{T}_2,\mathcal{T}_3,\ldots$ such that $r( \mathcal{T}_\star , \tilde{\boldsymbol{\lambda}})=j$ to define $E_j(\tilde{\boldsymbol{\lambda}}) \coloneqq E(\tilde{\boldsymbol{\lambda}}, \mathcal{T}_\star) $. Using the error functional $  E_j(\tilde{\boldsymbol{\lambda}})  $ as a starting point, we can also define modified errors  $\tilde{E}_j(\tilde{\boldsymbol{\lambda}})$. They have some similarities with the modified local errors given in \eqref{eq:2.3}. For $ j = 0  $ we put $\tilde{E}_0(\tilde{\boldsymbol{\lambda}})\coloneqq \tilde{e}(\tilde{\boldsymbol{\lambda}})$. For $  j \in \mathbb{N}  $ with $ j > 0  $ the error functionals $ \tilde{E}_j(\tilde{\boldsymbol{\lambda}})  $ are defined recursively via
\begin{equation}\label{eq:3.2}
\tilde{E}_j(\tilde{\boldsymbol{\lambda}})\coloneqq \frac{E_j(\tilde{\boldsymbol{\lambda}})\tilde{E}_{j-1}(\tilde{\boldsymbol{\lambda}})}{E_j(\tilde{\boldsymbol{\lambda}})+\tilde{E}_{j-1}(\tilde{\boldsymbol{\lambda}})}.
\end{equation}
In the special case $E_j(\tilde{\boldsymbol{\lambda}})=\tilde{E}_{j-1}(\tilde{\boldsymbol{\lambda}})=0$ we set $\tilde{E}_j(\tilde{\boldsymbol{\lambda}}) := 0$. The error functional $ \tilde{E}_j(\tilde{\boldsymbol{\lambda}})  $  can be reformulated in terms of some of the other error functionals which have been introduced above. For that purpose we use the definition of $\tilde{E}_j(\tilde{\boldsymbol{\lambda}})$ several times and plug in equation \eqref{eq:2.7}. Then we get
\begin{align}
\frac{1}{\tilde{E}_j(\tilde{\boldsymbol{\lambda}})} = \frac{1}{E_j(\tilde{\boldsymbol{\lambda}})}+\frac{1}{\tilde{E}_{j-1}(\tilde{\boldsymbol{\lambda}})} =\sum_{k=1}^j \frac{1}{E_k(\tilde{\boldsymbol{\lambda}})}+\frac{1}{\tilde{E}_0(\tilde{\boldsymbol{\lambda}})}
= \sum_{k=1}^j \frac{1}{E_k(\tilde{\boldsymbol{\lambda}})} + \sum_{ \tilde{\boldsymbol{\mu}} \in \mathcal{A}_{\mathcal{T}}(\tilde{\boldsymbol{\lambda}}) } \frac{1}{e(\tilde{\boldsymbol{\mu}})}. \label{eq:3.3}
\end{align}
Based on the definition of $ \tilde{E}_j(\tilde{\boldsymbol{\lambda}})   $ we can apply \eqref{eq:3.2} with $ j = r( \mathcal{T}, \tilde{\boldsymbol{\lambda}})   $ to define $ \tilde{E}(\tilde{\boldsymbol{\lambda}}) \coloneqq \tilde{E}(\tilde{\boldsymbol{\lambda}}, \mathcal{T}) \coloneqq \tilde{E}_{r( \mathcal{T} ,  \tilde{\boldsymbol{\lambda}})}(\tilde{\boldsymbol{\lambda}})$.

\vspace{0,2 cm}

\noindent
\underline{Step 3: Indicator functions for an adaptive decision.} 

\vspace{0,2 cm}

\noindent
Based on the error functionals we introduced above in what follows we define two indicator functions denoted by $a$ and $b$. They can be used to make an adaptive decision in our algorithm later on. Let $ \mathcal{T}  $ be a bivariate wavelet tree. Then we define a function $a: \mathcal{T} \rightarrow [0,\infty)$. For a leaf $\tilde{\boldsymbol{\lambda}} \in \mathcal{V}(\mathcal{T})$ we put $a(\tilde{\boldsymbol{\lambda}}) \coloneqq \tilde{e}(\tilde{\boldsymbol{\lambda}}) = \tilde{E}_0(\tilde{\boldsymbol{\lambda}})$. Given an inner node $\tilde{\boldsymbol{\lambda}} \in \mathcal{T} \backslash \mathcal{V}(\mathcal{T})$ the function $a$ is defined step by step moving from $ \tilde{\boldsymbol{\lambda}}  $ towards the leaves. We set
\begin{equation}\label{def_func_q1}
a(\tilde{\boldsymbol{\lambda}})  \coloneqq \min \left\{\max_{\tilde{\boldsymbol{\eta}} \in \mathcal{C}(\tilde{\boldsymbol{\lambda}},  \mathcal{T})}a(\tilde{\boldsymbol{\eta}}) , \tilde{E}_{r( \mathcal{T} , \tilde{\boldsymbol{\lambda}})}(\tilde{\boldsymbol{\lambda}})\right\}.
\end{equation}
The function $a$ serves as foundation when it comes to the definition of the decision function $b:\mathcal{T} \rightarrow\mathcal{V}(\mathcal{T})$. It maps each node of a bivariate wavelet tree $ \mathcal{T} $ to a leaf contained in $ \mathcal{V}(\mathcal{T})  $. Given a leaf $\tilde{\boldsymbol{\lambda}} \in \mathcal{V}(\mathcal{T})$ itself we put $b(\tilde{\boldsymbol{\lambda}}) \coloneqq \tilde{\boldsymbol{\lambda}}$. For an inner node $\tilde{\boldsymbol{\lambda}} \in \mathcal{T} \backslash \mathcal{V}(\mathcal{T})$ the function $b$ is defined step by step, whereby also the function $a$ is used. We put
\begin{equation*}\label{eq:3.5}
b(\tilde{\boldsymbol{\lambda}}) \coloneqq b\left( \operatorname{argmax}_{\tilde{\boldsymbol{\eta}} \in \mathcal{C}(\tilde{\boldsymbol{\lambda}}, \mathcal{T} )} a(\tilde{\boldsymbol{\eta}})  \right) .
\end{equation*}
Hence, the decision function $b$ points to that leaf of the investigated subtree with the largest penalized local error.

\subsection{An Algorithm for Adaptive Bivariate Quarklet Tree Approximation}\label{subsec_alg1}

Now we have all tools at hand to state our adaptive quarklet algorithm. As an input for the algorithm we can either use a function $f \in L_2((0,1)^2)$ or a sequence of its quarklet expansion coefficients. Below we use the notation $\mathbf{f}$ which stands for either of these two options. Then our algorithm called {\bf BIVARIATE\textunderscore NEARBEST\textunderscore TREE} adaptively produces a bivariate wavelet tree $  \mathcal{T}'_N  $. Recall that $\mathcal{R}$ stands for the root of the tree. 

\begin{footnotesize}

\vspace{0,4 cm}
\medskip
\noindent\shadowbox{\parbox{0.96\textwidth}{
		\begin{algorithm*} {\bf BIVARIATE\textunderscore NEARBEST\textunderscore TREE} $[\mathbf{f},N_{\max}] \mapsto \mathcal{T}'_N$\\[-3ex]
			\noindent{
				\begin{tabbing}
					set $\mathcal{T}'_0 \coloneqq \{\mathcal{R}\}$, $\tilde{e}(\mathcal{R}) \coloneqq e(\mathcal{R})$, $E_0(\mathcal{R}) \coloneqq e(\mathcal{R})$, $\tilde{E}_0(\mathcal{R}) \coloneqq \tilde{e}(\mathcal{R})$, $a(\mathcal{R}) \coloneqq \tilde{e}(\mathcal{R})$, $b(\mathcal{R}) \coloneqq \mathcal{R}$, $r( \mathcal{T}'_0 ,  \mathcal{R}) \coloneqq 0$;\\[0.1cm]
					\texttt{\em for}\={} $N  = 1$ \texttt{to} $N_{\max}$ \\[0.1cm]
set $\tilde{\boldsymbol{\lambda}}_N \coloneqq b(\mathcal{R})$ and compute $\alpha_{N} \coloneqq \alpha (\tilde{\boldsymbol{\lambda}}_N)$;   \\[0.1cm]
\texttt{\em if} $\alpha_{N} = 0$       \\[0.1cm]
					\texttt{\em for}\={} $i \in \{ 1, 2 \}$  \\[0.1cm]
\>  expand the current tree $\mathcal{T}'_{N-1}$ to $(\mathcal{T}'_{N})_{i}$ by subdividing $\tilde{\boldsymbol{\lambda}}_N = b(\mathcal{R})$ and \\[0.1cm]
					\> adding its available children $\hat{\tilde{\boldsymbol{\eta}}} \in \mathcal{C}(\tilde{\boldsymbol{\lambda}}_N , (LSR.a.i))$ consistent with strategy (LSR.a.i) to it; \\[0.1cm]
\> compute $A_{i}(\tilde{\boldsymbol{\lambda}}_N) \coloneqq \sum_{\hat{\tilde{\boldsymbol{\eta}}} \in \mathcal{C}(\tilde{\boldsymbol{\lambda}}_N, (LSR.a.i))} e_{0}(\hat{\tilde{\boldsymbol{\eta}}})$ or put $A_{i}(\tilde{\boldsymbol{\lambda}}_N)= + \infty $ if (LSR.a.i) is not available;     \\[0.1cm]										\texttt{\em end for} \\[0.1cm]
\texttt{\em if} $ A_{1}(\tilde{\boldsymbol{\lambda}}_N) \leq A_{2}(\tilde{\boldsymbol{\lambda}}_N) $       \\[0.1cm]
\> put $ i^{\star} \coloneqq 1    $; \\[0.1cm]
\texttt{\em else} \\[0.1cm]
\> put $  i^{\star} \coloneqq 2        $; \\[0.1cm]
\texttt{\em end if} \\[0.1cm]
 expand the current tree $\mathcal{T}'_{N-1}$ to $\mathcal{T}'_{N} \coloneqq  (\mathcal{T}'_{N})_{i^{\star}}$ by subdividing $\tilde{\boldsymbol{\lambda}}_N = b(\mathcal{R})$ and \\[0.1cm]
adding its available children $\hat{\tilde{\boldsymbol{\eta}}} \in \mathcal{C}(\tilde{\boldsymbol{\lambda}}_N , (LSR.a.i^{\star}))$ consistent with strategy (LSR.a.$i^{\star}$) to it; \\[0.1cm]
\texttt{\em for}\={} $\tilde{\boldsymbol{\lambda}} = \hat{\tilde{\boldsymbol{\eta}}} \in \mathcal{C}(\tilde{\boldsymbol{\lambda}}_N , (LSR.a.i^{\star})) $ \\
					\>  calculate $\tilde{e}(\tilde{\boldsymbol{\lambda}}) \coloneqq \frac{e(\tilde{\boldsymbol{\lambda}})\tilde{e}(\tilde{\boldsymbol{\lambda}}_N)}{e(\tilde{\boldsymbol{\lambda}}) + \tilde{e} (\tilde{\boldsymbol{\lambda}}_N)}$, $E_0(\tilde{\boldsymbol{\lambda}}) \coloneqq e(\tilde{\boldsymbol{\lambda}})$,  $\tilde{E}_0(\tilde{\boldsymbol{\lambda}}) \coloneqq \tilde{e}(\tilde{\boldsymbol{\lambda}})$, $a(\tilde{\boldsymbol{\lambda}}) \coloneqq \tilde{e}(\tilde{\boldsymbol{\lambda}})$, $b(\tilde{\boldsymbol{\lambda}}) \coloneqq \tilde{\boldsymbol{\lambda}}$, $r( \mathcal{T}'_{N} ,   \tilde{\boldsymbol{\lambda}}) \coloneqq 0$;\\[0.1cm]    
\texttt{\em end for} \\
set $\tilde{\boldsymbol{\lambda}} = \tilde{\boldsymbol{\lambda}}_N$; \\[0.1cm]					
\texttt{\em else if} $  \alpha_{N} \in \{ 1, 2 \}   $     \\[0.1cm]
expand the current tree $\mathcal{T}'_{N-1}$ to $\mathcal{T}'_{N}$ by subdividing $\tilde{\boldsymbol{\lambda}}_N = b(\mathcal{R})$ and \\[0.1cm]
add its children $\hat{\tilde{\boldsymbol{\eta}}} \in \mathcal{C}(\tilde{\boldsymbol{\lambda}}_N , \clubsuit )$ with $ \clubsuit = (LSR.b)   $ if $ \alpha_{N} = 1  $ or $ \clubsuit = (LSR.c)   $ if $ \alpha_{N} = 2  $ to it; \\[0.1cm]
\texttt{\em for}\={} $\tilde{\boldsymbol{\lambda}} = \hat{\tilde{\boldsymbol{\eta}}} \in \mathcal{C}(\tilde{\boldsymbol{\lambda}}_N , \clubsuit ) $ \\
					\>  calculate $\tilde{e}(\tilde{\boldsymbol{\lambda}}) \coloneqq \frac{e(\tilde{\boldsymbol{\lambda}})\tilde{e}(\tilde{\boldsymbol{\lambda}}_N)}{e(\tilde{\boldsymbol{\lambda}}) + \tilde{e} (\tilde{\boldsymbol{\lambda}}_N)}$, $E_0(\tilde{\boldsymbol{\lambda}}) \coloneqq e(\tilde{\boldsymbol{\lambda}})$,  $\tilde{E}_0(\tilde{\boldsymbol{\lambda}}) \coloneqq \tilde{e}(\tilde{\boldsymbol{\lambda}})$, $a(\tilde{\boldsymbol{\lambda}}) \coloneqq \tilde{e}(\tilde{\boldsymbol{\lambda}})$, $b(\tilde{\boldsymbol{\lambda}}) \coloneqq \tilde{\boldsymbol{\lambda}}$, $r( \mathcal{T}'_{N} ,  \tilde{\boldsymbol{\lambda}}) \coloneqq 0$;\\[0.1cm]    
\texttt{\em end for} \\
set $\tilde{\boldsymbol{\lambda}} = \tilde{\boldsymbol{\lambda}}_N$; \\[0.1cm]
					
\texttt{\em end if}\\[0.1cm]
\texttt{\em whi}\=\texttt{\em le} $\tilde{\boldsymbol{\lambda}} \neq \emptyset$ \\[0.1cm]
					 \> set $r( \mathcal{T}'_{N} ,  \tilde{\boldsymbol{\lambda}}) \coloneqq r( \mathcal{T}'_{N - 1} ,  \tilde{\boldsymbol{\lambda}}) + 1 $; calculate $e_{r( \mathcal{T}'_{N} ,  \tilde{\boldsymbol{\lambda}})}(\tilde{\boldsymbol{\lambda}})$; set $\tilde{\boldsymbol{\eta}} \in \mathcal{C}(\tilde{\boldsymbol{\lambda}}, \mathcal{T}'_{N} )$ to be the children of $\tilde{\boldsymbol{\lambda}}$;\\[0.1cm]
					 \> set $E_{r( \mathcal{T}'_{N} ,  \tilde{\boldsymbol{\lambda}})}(\tilde{\boldsymbol{\lambda}}) \coloneqq \min\{ \sum_{\tilde{\boldsymbol{\eta}} \in \mathcal{C}(\tilde{\boldsymbol{\lambda}},  \mathcal{T}'_{N} )} E_{r( \mathcal{T}'_{N} , \tilde{\boldsymbol{\eta}})}(\tilde{\boldsymbol{\eta}}), e_{r( \mathcal{T}'_{N} ,  \tilde{\boldsymbol{\lambda}})}(\tilde{\boldsymbol{\lambda}})\}$; \\[0.1cm]
					 \> set $\tilde{E}_{r( \mathcal{T}'_{N} ,  \tilde{\boldsymbol{\lambda}})}(\tilde{\boldsymbol{\lambda}})\coloneqq \frac{E_{r( \mathcal{T}'_{N} , \tilde{\boldsymbol{\lambda}})}(\tilde{\boldsymbol{\lambda}})\tilde{E}_{{r( \mathcal{T}'_{N} ,  \tilde{\boldsymbol{\lambda}})}-1}(\tilde{\boldsymbol{\lambda}})}{E_{r( \mathcal{T}'_{N} ,  \tilde{\boldsymbol{\lambda}})}(\tilde{\boldsymbol{\lambda}})+\tilde{E}_{{r( \mathcal{T}'_{N} ,  \tilde{\boldsymbol{\lambda}})}-1}(\tilde{\boldsymbol{\lambda}})}$; \\[0.1cm]
					 \> set $\tilde{\boldsymbol{\eta}}^{\star} \coloneqq \operatorname{argmax}_{\tilde{\boldsymbol{\eta}} \in \mathcal{C}(\tilde{\boldsymbol{\lambda}}, \mathcal{T}'_{N}  )} a(\tilde{\boldsymbol{\eta}}) $, $a(\tilde{\boldsymbol{\lambda}}) \coloneqq \min\{a(\tilde{\boldsymbol{\eta}}^{\star}), \tilde{E}_{r( \mathcal{T}'_{N} ,  \tilde{\boldsymbol{\lambda}})}(\tilde{\boldsymbol{\lambda}})\}$ and $b(\tilde{\boldsymbol{\lambda}}) \coloneqq b(\tilde{\boldsymbol{\eta}}^{\star})$; \\[0.1cm]
					 \> replace $\tilde{\boldsymbol{\lambda}}$ with its parent (or $\emptyset$ if $\tilde{\boldsymbol{\lambda}} = \mathcal{R}$); \\[0.1cm]
 \texttt{\em end while}\\[0.1cm]					
					          \texttt{\em end for} 
				\end{tabbing}
			}
		\end{algorithm*}
}}
\medskip
\vspace{0,4 cm}

\end{footnotesize}

The algorithm {\bf BIVARIATE\textunderscore NEARBEST\textunderscore TREE} has many similarities with its univariate forerunner given in \cite{DaHoRaVo}, see Section 3.2. As in the univariate case we start with a tree $\mathcal{T}'_0 \coloneqq \{\mathcal{R}\}$ and expand it step by step. As long as we have $  N \leq   N_{\max}  $ we work with a tree $\mathcal{T}'_{N-1}$ and subdivide its leaf $   b(\mathcal{R}) $ by adding the child nodes according to one of the possible refinement options to it in order to obtain $\mathcal{T}'_{N}$. Then for these children in a for-loop the error functionals are computed. Moreover, we use a while-loop to update all important quantities going from the new leaves back to the root $  \mathcal{R}  $. Here especially the modified error $  \tilde{E}   $ and the functions $a$ and $b$ are essential since they allow for the adaptive decision where to refine in the next step of the algorithm. This part of the algorithm works similar as in the univariate setting and therefore also is explained in detail in \cite{DaHoRaVo}, see Section 3.2. However, there also is an important difference compared to the univariate case. When we subdivide a leaf $   b(\mathcal{R}) $ in the algorithm {\bf BIVARIATE\textunderscore NEARBEST\textunderscore TREE} we have different possibilities according to the refinement strategies (LSR.a.1), (LSR.a.2), (LSR.b) and (LSR.c). For that reason in each step we calculate $  \alpha_{N} = \alpha (   b(\mathcal{R})  )   $. If $  \alpha_{N} \in \{ 1, 2 \}   $ we apply the related refinement strategies (LSR.b) or (LSR.c). This is done in the lower part of the algorithm. Else if $\alpha_{N} = 0$ one of the strategies (LSR.a.1) or (LSR.a.2) has to be used. To select the best possible strategy we compute the resulting local errors of the lowest order for the new children and then choose the option with a smaller cumulated local error. In connection with that the quantities $A_{i}(\tilde{\boldsymbol{\lambda}}_N)$ show up and $ i^{\star}  $ refers to the better option. This decision is described in the upper part of the algorithm. Investigating the complexity of the algorithm {\bf BIVARIATE\textunderscore NEARBEST\textunderscore TREE} we obtain the following lemma. 

\begin{Lemma}\label{lem_complexity}
Let $N \in \mathbb{N}$. Then the algorithm {\bf BIVARIATE\textunderscore NEARBEST\textunderscore TREE} performs $ \sum_{\tilde{\boldsymbol{\lambda}} \in  \mathcal{T}'_{N}} ( r( \mathcal{T}'_{N} ,  \tilde{\boldsymbol{\lambda}}) + 1 )     $ steps to obtain $  \mathcal{T}'_{N}  $.
\end{Lemma}

\begin{proof}
This result can be proved with similar methods as Lemma 3.2 in \cite{Bin18}, see also Lemma 3.3 in \cite{DaHoRaVo}. The number of steps in the algorithm {\bf BIVARIATE\textunderscore NEARBEST\textunderscore TREE} is determined by the outer for-loop where $N$ runs from $1$ to $N_{\max}$ and an inner while-loop in which the calculations at the nodes of the tree, starting at the newly subdivided node and then returning to the root, are performed. For a new node $ \tilde{\boldsymbol{\lambda}}  $ of the tree the quantity $ r( \mathcal{T}'_{N} ,  \tilde{\boldsymbol{\lambda}})   $ is initialized as $0$ and then increased by $1$ whenever the node $ \tilde{\boldsymbol{\lambda}}  $ is revisited in the inner while-loop of the algorithm later on. Consequently the number $( r( \mathcal{T}'_{N} ,  \tilde{\boldsymbol{\lambda}}) + 1 )$ counts how many times the node $  \tilde{\boldsymbol{\lambda}}   $ is visited by the algorithm, whereby each visit is connected with a small number of calculations. Taking the sum over all $ \tilde{\boldsymbol{\lambda}} \in  \mathcal{T}'_{N}   $ therefore delivers the total number of steps performed in the algorithm {\bf BIVARIATE\textunderscore NEARBEST\textunderscore TREE}.
\end{proof}
Lemma \ref{lem_complexity} looks like its univariate counterpart which is given in \cite{DaHoRaVo}, see Lemma 3.3. Nevertheless the algorithm {\bf BIVARIATE\textunderscore NEARBEST\textunderscore TREE} performs a larger number of steps than its univariate forerunner. The main reason for this is the fact that in case of the refinement strategies (LSR.a.1) and (LSR.a.2) three new children are added instead of two in the univariate setting. Consequently the tree $ \mathcal{T}'_{N}  $ consists of more nodes than its univariate counterpart. Moreover, some of the steps in the algorithm {\bf BIVARIATE\textunderscore NEARBEST\textunderscore TREE} are connected with a larger number of computations compared to the univariate algorithm. To see this, recall that if $ \mathcal{T}'_{N}   $ is created by {\bf BIVARIATE\textunderscore NEARBEST\textunderscore TREE}, we have to decide $N$ times which of the strategies (LSR.a.1), (LSR.a.2), (LSR.b) or (LSR.c) is chosen. Each of these choices is connected with a number of calculations.

\subsection{The Process of Trimming}

The tree $ \mathcal{T}'_N  $ produced by the algorithm {\bf BIVARIATE\textunderscore NEARBEST\textunderscore TREE} consists of bivariate wavelet indices only. However, it can be transformed into a bivariate quarklet tree $T_{N} = (\mathcal{T}_N, \mathcal{T}'_N)$ easily. An important step to carry out this transformation is the process of trimming. It is applied in order to obtain the optimal subtree $\mathcal{T}_N$ of $\mathcal{T}'_N$. For that purpose we start with the wavelet tree $\mathcal{T}'_N$. Then we walk from the root $ \mathcal{R}$ towards one of the  leaves $\tilde{\boldsymbol{\eta}} \in \mathcal{V}(\mathcal{T}'_N)$. At the first node where we observe $E(\tilde{\boldsymbol{\lambda}}) = e_{r(  \mathcal{T}'_N , \tilde{\boldsymbol{\lambda}})}(\tilde{\boldsymbol{\lambda}})$ in \eqref{eq:3.1} we trim the tree. Therefore we delete all descendants of $\tilde{\boldsymbol{\lambda}}$. Recall, that by definition we have $E(\tilde{\boldsymbol{\eta}}) = e_0(\tilde{\boldsymbol{\eta}})$ on the leaves $\tilde{\boldsymbol{\eta}} \in \mathcal{V}(\mathcal{T}'_N)$, see Section \ref{Subsec_Mod_Errfunc}. Hence this situation will surely show up after some steps. To continue this procedure is repeated for all remaining paths which have not been treated so far. Consequently, $\mathcal{T}_N$ becomes the minimal tree with $E(\tilde{\boldsymbol{\lambda}}) = e_{r(  \mathcal{T}'_N ,  \tilde{\boldsymbol{\lambda}})}(\tilde{\boldsymbol{\lambda}})$ on all leaves. Now the tree $\mathcal{T}'_N$ and its subtree $ \mathcal{T}_N $ can be used to obtain a bivariate quarklet tree $T_N  = (\mathcal{T}_N, \mathcal{T}'_N) = (\mathcal{T}_N, P_{\max})$ by setting $p_{\max}(\tilde{\boldsymbol{\lambda}}^{\circ}) := r( \mathcal{T}_N , \mathcal{T}'_N, \tilde{\boldsymbol{\lambda}}^{\circ} )$ on each leaf $ \tilde{\boldsymbol{\lambda}}^{\circ} \in \mathcal{V}(\mathcal{T}_N)   $ as explained in Section \ref{subsec_biv_qua_tre}. Recall, that by Definition \ref{def:quarklet_tree} this already implies the polynomial degrees on all nodes of the tree. The following algorithm {\bf BIVARIATE\textunderscore TRIM} provides one possible way to implement the trimming procedure.

\vspace{0,4 cm}

\medskip
\noindent\shadowbox{\parbox{0.96\textwidth}{
		\begin{algorithm*} {\bf BIVARIATE\textunderscore TRIM} $[\mathcal{T}'] \mapsto \mathcal{T}$\\[-3ex]
			\noindent{
				\begin{tabbing}
					set $B = \{\mathcal{R}\}$ and $ \mathcal{T} = \mathcal{T}'$;\\[0.1cm]
					\texttt{\em whi}\=\texttt{\em le} $B \neq \emptyset$ \\[0.1cm]
					\> take $\tilde{\boldsymbol{\lambda}} \in B$; \\[0.1cm]
					\> \texttt{\em if}\={} $E(\tilde{\boldsymbol{\lambda}}) = e_{r( \mathcal{T}' ,  \tilde{\boldsymbol{\lambda}})}(\tilde{\boldsymbol{\lambda}})$\\[0.1cm]
					\> \> remove all descendants from $\tilde{\boldsymbol{\lambda}}$ in $\mathcal{T}$;\\[0.1cm]
					\> \texttt{\em else}\\[0.1cm]
					\> \> add the children $\tilde{\boldsymbol{\eta}} \in \mathcal{C}(\tilde{\boldsymbol{\lambda}} ,  \mathcal{T}' )$  of $\tilde{\boldsymbol{\lambda}}$ according to $  \mathcal{T}' $ to $B$;\\[0.1cm]
					\> \texttt{\em end if}\\[0.1cm]
					\> remove $\tilde{\boldsymbol{\lambda}}$ from $B$;\\[0.1cm]
					\texttt{\em end while}
				\end{tabbing}
			}
		\end{algorithm*}
}}
\medskip

\vspace{0,4 cm}

Recall, that if $\mathcal{T}'_N$ is created by the algorithm {\bf BIVARIATE\textunderscore NEARBEST\textunderscore TREE} the quantities $E(\tilde{\boldsymbol{\lambda}})$ and  $e_{r(  \mathcal{T}'_N ,  \tilde{\boldsymbol{\lambda}})}(\tilde{\boldsymbol{\lambda}})$ already have been computed there. In this case no further calculations are needed to run the algorithm {\bf BIVARIATE\textunderscore TRIM}. In a next step we estimate the cardinality of the bivariate quarklet tree $T_N$ obtained above.

\begin{Lemma} \label{lemma:T_N_cardinality}
	Let $  N \in \mathbb{N}  $ with $N \geq 3$.	Let the bivariate quarklet tree $T_N = (\mathcal{T}_N,\mathcal{T}'_N)$ be created by the algorithm {\bf BIVARIATE\textunderscore NEARBEST\textunderscore TREE} and a subsequent trimming. Then it holds
	\begin{equation}\label{T_N_cardinality}
	2N+1 \le \# T_N \le  \frac{1}{2} N^3   +  \frac{16}{5} N^2  +    \frac{25}{6} N +  3  .
	\end{equation}
\end{Lemma}
\begin{proof}
Let $T$ be a bivariate quarklet tree, $\mathcal{T}$ the underlying wavelet tree and $\mathcal{R}$ its root.	Recall, that each bivariate quarklet tree $T$ can be described by two types of refinement. The first one is a sequence of refinements in space, which can be depicted via a bivariate wavelet tree $\mathcal{T}$. The second one can be expressed by the steps of polynomial enrichment of $\mathcal{T}$, characterized by $\{ p_{\textrm{max}}(\tilde{\boldsymbol{\lambda}}^{\circ})\}_{\tilde{\boldsymbol{\lambda}} ^{\circ}\in \mathcal{V}(\mathcal{T})}$. Now let $N_h \coloneqq r(\mathcal{T}, \mathcal{R})$ be the total number of space refinements of the form (LSR.a.1), (LSR.a.2), (LSR.b) or (LSR.c) that is necessary to create the bivariate wavelet tree $\mathcal{T}$. Moreover, let $N_p \coloneqq \sum_{\tilde{\boldsymbol{\lambda}}^{\circ} \in \mathcal{V}(\mathcal{T})} p_{\max}(\tilde{\boldsymbol{\lambda}}^{\circ})$. Then $N=N_h+N_p$ denotes the total number of refinements in space and polynomial degree, which is necessary to create the bivariate quarklet tree $T$. A refinement in space always increases the cardinality of the tree by two or three depending on the selected strategy (LSR.a.1), (LSR.a.2), (LSR.b) or (LSR.c). On the other hand increasing the polynomial degree on a leaf $\tilde{\boldsymbol{\lambda}}^{\circ}$ enlarges the cardinality depending on the size of the set $\Upsilon( \tilde{\boldsymbol{\lambda}}^{\circ}  )$. Since the set $ \Upsilon(\tilde{\boldsymbol{\lambda}}^{\circ}) $ has the form $ \Upsilon(\tilde{\boldsymbol{\lambda}}^{\circ})= \{ \tilde{\boldsymbol{\mu}} \in \mathcal{T} : \tilde{\boldsymbol{\lambda}}^{\circ} \succeq  \tilde{\boldsymbol{\mu}} \succeq \boldsymbol{\tilde{\mu}_{\tilde{\lambda}}} \} $ with a fixed $ \boldsymbol{\tilde{\mu}_{\tilde{\lambda}}}  \in \mathcal{T}$, see \eqref{eq:Upsilon}, we observe
	\[1 \le \vert \Upsilon( \tilde{\boldsymbol{\lambda}}^{\circ}  ) \vert \le \vert \{ \tilde{\boldsymbol{\mu}}   \in \mathcal{T} : \tilde{\boldsymbol{\lambda}}^{\circ} \succeq \tilde{\boldsymbol{\mu}} \succeq \mathcal{R}\}\vert \le \vert \tilde{\boldsymbol{\lambda}}^{\circ} \vert +1.\]
In what follows we prove the lower estimate in \eqref{T_N_cardinality}. If we refine $N_p = N$ times in polynomial degree on the node $ \tilde{\boldsymbol{\lambda}}^{\circ}  =\mathcal{R}$ with $\vert \Upsilon(\mathcal{R}) \vert = 1$ we obtain
\begin{equation*}
\#T_{N} \coloneqq  \frac{( p_{\textrm{max}}(\mathcal{R}) + 1)^2 + (  p_{\textrm{max}}(\mathcal{R}) + 1  )}{2} =  \frac{( N + 1)^2 + (  N + 1  )}{2}  ,
\end{equation*}
see \eqref{eq_card_tree11}. Else, if we refine $N_{h} = N$ times in space, we get the estimate $ \#T_{N} \geq 1 + 2N $. Here in order to find a lower estimate we assumed that in each refinement step we added exactly two children. Looking at the case $N=N_h+N_p$ we observe that increasing $  p_{\textrm{max}}( \tilde{\boldsymbol{\lambda}}^{\circ} )   $ from $0$ to $1$ on a single node $    \tilde{\boldsymbol{\lambda}}^{\circ}   $ raises the cardinality of the bivariate quarklet tree by two at least. Any further increase of $  p_{\textrm{max}}( \tilde{\boldsymbol{\lambda}}^{\circ} ) \in \mathbb{N}   $ to $  p_{\textrm{max}}( \tilde{\boldsymbol{\lambda}}^{\circ} ) + 1   $  raises the cardinality by $   p_{\textrm{max}}( \tilde{\boldsymbol{\lambda}}^{\circ} ) + 2  $ at least. Consequently,  to obtain the lower estimate in \eqref{T_N_cardinality}, we only refine in space $N_{h} = N$ times. Now we want to prove the upper bound. For that purpose we have to investigate how we can create the bivariate quarklet tree $T$ which maximizes  $\# T$ after $N$ refinement steps. In a single step, the largest increase in cardinality that is possible for a tree $(\mathcal{T},P_{\max})$ of depth $J = \max_{\tilde{\boldsymbol{\lambda}}^{\circ} \in \mathcal{T}} |\tilde{\boldsymbol{\lambda}}^{\circ}|$ by means of polynomial enrichment can show up if there exits a leaf $\tilde{\boldsymbol{\lambda}}^{\circ} \in \mathcal{V}(\mathcal{T})$  with $\vert \tilde{\boldsymbol{\lambda}}^{\circ} \vert = J$ and $\Upsilon( \tilde{\boldsymbol{\lambda}}^{\circ}  ) 
	= \{ \tilde{\boldsymbol{\mu}}  \in \mathcal{T} : \tilde{\boldsymbol{\lambda}}^{\circ} \succeq \tilde{\boldsymbol{\mu}} \succeq \mathcal{R}\}$. In this case one step of polynomial enrichment of $\tilde{\boldsymbol{\lambda}}^{\circ}$ will increase the cardinality of the quarklet tree by $ (J + 1) ( p_{\max}(\tilde{\boldsymbol{\lambda}}^{\circ}) + 1) $. On the other hand we avoid having many leaves on a high level since space refinement increases the cardinality only by two or three. Hence, the largest possible bivariate quarklet tree after $N$ refinement steps consists only of leaves and a single path to a leaf $\tilde{\boldsymbol{\lambda}}^{\circ}$ on a high level with $\Upsilon(\tilde{\boldsymbol{\lambda}}^{\circ}) 
	= \{  \tilde{\boldsymbol{\mu}}  \in \mathcal{T} : \tilde{\boldsymbol{\lambda}}^{\circ} \succeq \tilde{\boldsymbol{\mu}} \succeq \mathcal{R}\}$ and polynomial enrichment is applied only on this leaf. To obtain this situation we first have to employ $N_h$ steps of space refinement along this path such that we have $|\tilde{\boldsymbol{\lambda}}^{\circ}| = N_h$. Then we refine the polynomial degree $N_p$-times on the leaf $ \tilde{\boldsymbol{\lambda}}^{\circ} $. The cardinality of such a tree can be estimated by
	\begin{equation}\label{upper_bound_card}
\#T \leq   1   +  3 N_{h} +    (N_{h} + 1) \Big (  \frac{( N - N_{h} + 1)^2 + (  N - N_{h} + 1  )}{2}  -    1   \Big ) .
\end{equation}
For $N \geq 3$ the right hand side has its maximum in $[0,N]$ at 
\begin{align*}
N_{h}  =   \frac{1}{3} \Big (  -  \sqrt{N^2 + 5 N - 5} +  2 N + 2    \Big ) .
\end{align*}
To obtain an (almost) sharp upper estimate for the cardinality of the bivariate quarklet tree we can plug in $ N_{h} $ into the right hand side of    \eqref{upper_bound_card}. In order to present a result in a clearly arranged way we can further estimate
\begin{align*}
N_{h} \leq   \frac{1}{3} \Big ( N + 2    \Big ) \qquad \qquad \mbox{and} \qquad \qquad N - N_{h} \leq \frac{5}{3} N .
\end{align*}
Using this in combination with \eqref{upper_bound_card} we finally get
\begin{align*}
\#T \leq \frac{1}{2} N^3   +  \frac{16}{5} N^2  +    \frac{25}{6} N +  3 .
\end{align*}
The proof is complete. 
\end{proof}

\subsection{The Bivariate Quarklet Trees \texorpdfstring{$T_N$}{TN} are Near-Best}\label{subsec_alg_nearbest1}

In this section we prove that the bivariate quarklet trees produced by the algorithm {\bf BIVARIATE\textunderscore NEARBEST\textunderscore TREE} with a subsequent trimming are actually near-best in the sense of \eqref{eq:nearbest}. To see this two substeps have to be carried out. At first we establish a lower bound for the best approximation error $\sigma_n$ using the parameter $a_N = a(\mathcal{R})$ for a given wavelet tree with root $ \mathcal{R}  $.
\begin{Lemma}\label{lemma:lower1}
Let $  n,N \in \mathbb{N}  $ with $n \le N$. Let $T^\star = (\mathcal{T}^\star, P_{\max}^\star)$ be the optimal bivariate quarklet tree of cardinality $n $ such that $\sigma_n = \mathcal{E}(T^\star)$. Let the quarklet tree $T_N = (\mathcal{T}_N, \mathcal{T}'_N) $ be created by the algorithm {\bf BIVARIATE\textunderscore NEARBEST\textunderscore TREE} with a subsequent trimming. We define the parameter $a_N\coloneqq a(\mathcal{R})$ for the wavelet tree $\mathcal{T}'_N$ with root $ \mathcal{R} $. Then it holds
	\begin{equation*}\label{eq:lower1}
		\sigma_n \ge   a_N \Big (  N    -  \frac{2}{3}  n + \frac{1}{2} \Big )   .
	\end{equation*}
\end{Lemma}

\begin{proof}
Let $T^\star = (\mathcal{T}^\star, P_{\max}^\star)$ be the optimal bivariate quarklet tree of cardinality $n $ such that $\sigma_n = \mathcal{E}(T^\star)$. In order to obtain a lower estimate for $ \sigma_n  $ we consider the leaves $\tilde{\boldsymbol{\lambda}}^{\circ} \in \mathcal{V}(\mathcal{T}^\star)$ and their orders $ P_{\max}^\star = \{p_{\textrm{max}}^\star( \tilde{\boldsymbol{\lambda}}^{\circ}  )\}_{  \tilde{\boldsymbol{\lambda}}^{\circ}  \in \mathcal{V}(\mathcal{T}^\star)}$. For $r( \mathcal{T}'_N, \tilde{\boldsymbol{\lambda}}^{\circ}  ) \le p_{\max}^\star( \tilde{\boldsymbol{\lambda}}^{\circ} )$ we ignore the contribution of $e_{p_{\max}^\star(\tilde{\boldsymbol{\lambda}}^{\circ})}(\tilde{\boldsymbol{\lambda}}^{\circ})$ to the gobal error $\mathcal{E}(T^\star)$. Then we get
     \begin{equation}\label{eq:split}
     	\sigma_n = \mathcal{E}(T^\star) =  \sum_{\tilde{\boldsymbol{\lambda}}^{\circ} \in \mathcal{V}(\mathcal{T}^\star)} e_{p_{\max}^\star(\tilde{\boldsymbol{\lambda}}^{\circ})}( \tilde{\boldsymbol{\lambda}}^{\circ} ) \ge \sum_{\substack{ \tilde{\boldsymbol{\lambda}}^{\circ}  \in \mathcal{V}(\mathcal{T}^\star), \\ r( \mathcal{T}'_N, \tilde{\boldsymbol{\lambda}}^{\circ} ) > p_{\max}^\star( \tilde{\boldsymbol{\lambda}}^{\circ}  )}} e_{p_{\max}^\star( \tilde{\boldsymbol{\lambda}}^{\circ} )}( \tilde{\boldsymbol{\lambda}}^{\circ}  ).
     \end{equation} 
Here we used \eqref{eq_glob_err}. To continue let $ k \in \mathbb{N}_0  $ with $ k \leq N  $. For the remaining leaves with $r( \mathcal{T}'_N, \tilde{\boldsymbol{\lambda}}^{\circ}  ) > p_{\max}^\star( \tilde{\boldsymbol{\lambda}}^{\circ}  )$ we consider $a_k \coloneqq a(\mathcal{R})$ at the stage $\mathcal{T}'_k$ of growing the wavelet tree $\mathcal{T}'_N$ using the algorithm {\bf BIVARIATE\textunderscore NEARBEST\textunderscore TREE} at the last increase of $r( \mathcal{T}'_N, \tilde{\boldsymbol{\lambda}}^{\circ})$. That means $\mathcal{T}'_k$ is the last bivariate wavelet tree in the sequence of trees $\mathcal{T}'_1,\ldots,\mathcal{T}'_N$ created by the algorithm {\bf BIVARIATE\textunderscore NEARBEST\textunderscore TREE} where a descendant of $\tilde{\boldsymbol{\lambda}}^{\circ}$ is enclosed. Recall, that the numbers $a_k$ are decreasing with $k$, such that we have $a_k \ge a_N$. Using the definitions of the decision functions $a$ and $b$ it follows that at this stage of growing the wavelet tree we have $b(\mathcal{R}) = b(\tilde{\boldsymbol{\lambda}}^{\circ}) = b(\tilde{\boldsymbol{\mu}}^{\circ})$ for all $ \tilde{\boldsymbol{\mu}}^{\circ} \preceq \tilde{\boldsymbol{\lambda}}^{\circ} $ with $ \tilde{\boldsymbol{\mu}}^{\circ}  \in \mathcal{T}'_k  $. Now let $\tilde{\boldsymbol{\mu}}^{\circ}$ be the parent of $\tilde{\boldsymbol{\lambda}}^{\circ}$. Then an application of \eqref{def_func_q1} yields
\begin{align*}
	a(\tilde{\boldsymbol{\mu}}^{\circ}) \coloneqq \min \left\{\max_{\tilde{\boldsymbol{\eta}}^{\circ} \in \mathcal{C}(\tilde{\boldsymbol{\mu}}^{\circ} , \mathcal{T}'_{k}  )}a(\tilde{\boldsymbol{\eta}}^{\circ}) , \tilde{E}_{r( \mathcal{T}'_{k} , \tilde{\boldsymbol{\mu}}^{\circ} )}( \tilde{\boldsymbol{\mu}}^{\circ}  )\right\} = \min \{a(\tilde{\boldsymbol{\lambda}}^{\circ}), \tilde{E}_{r( \mathcal{T}'_k, \tilde{\boldsymbol{\mu}}^{\circ}  )}(\tilde{\boldsymbol{\mu}}^{\circ})\} \le a(\tilde{\boldsymbol{\lambda}}^{\circ}).
\end{align*}
Using this argument several times we obtain $a(\tilde{\boldsymbol{\lambda}}^{\circ}) \ge a_k =a(\mathcal{R})$ and $\tilde{E}_j(\tilde{\boldsymbol{\lambda}}^{\circ}) \ge a(\tilde{\boldsymbol{\lambda}}^{\circ}) \ge a_N$ with $j=r(\mathcal{T}'_N, \tilde{\boldsymbol{\lambda}}^{\circ}  )-1$. Here again \eqref{def_func_q1} has been applied. To continue we can argue as in \eqref{eq:3.3} to get	
\begin{align*}
	\frac{1}{\tilde{E}_j(\tilde{\boldsymbol{\lambda}}^{\circ})} = \sum_{\ell = p_{\max}^\star(\tilde{\boldsymbol{\lambda}}^{\circ})+1}^j \frac{1}{E_{\ell}(\tilde{\boldsymbol{\lambda}}^{\circ})} + \frac{1}{\tilde{E}_{p_{\max}^\star(\tilde{\boldsymbol{\lambda}}^{\circ} )}(\tilde{\boldsymbol{\lambda}}^{\circ})}.
\end{align*}
This also yields
\begin{align*}
E_{p_{\max}^\star(\tilde{\boldsymbol{\lambda}}^{\circ} )}(\tilde{\boldsymbol{\lambda}}^{\circ}) & =  E_{p_{\max}^\star(\tilde{\boldsymbol{\lambda}}^{\circ} )}(\tilde{\boldsymbol{\lambda}}^{\circ})  \tilde{E}_j(\tilde{\boldsymbol{\lambda}}^{\circ})  \left ( \sum_{\ell = p_{\max}^\star(\tilde{\boldsymbol{\lambda}}^{\circ})+1}^j \frac{1}{E_{\ell}(\tilde{\boldsymbol{\lambda}}^{\circ})} + \frac{1}{\tilde{E}_{p_{\max}^\star(\tilde{\boldsymbol{\lambda}}^{\circ} )}(\tilde{\boldsymbol{\lambda}}^{\circ})} \right ) \\
& =   \tilde{E}_j(\tilde{\boldsymbol{\lambda}}^{\circ})  \left ( \sum_{\ell = p_{\max}^\star(\tilde{\boldsymbol{\lambda}}^{\circ})+1}^j \frac{ E_{p_{\max}^\star(\tilde{\boldsymbol{\lambda}}^{\circ} )}(\tilde{\boldsymbol{\lambda}}^{\circ})}{E_{\ell}(\tilde{\boldsymbol{\lambda}}^{\circ})} + \frac{ E_{p_{\max}^\star(\tilde{\boldsymbol{\lambda}}^{\circ} )}(\tilde{\boldsymbol{\lambda}}^{\circ})}{\tilde{E}_{p_{\max}^\star(\tilde{\boldsymbol{\lambda}}^{\circ} )}(\tilde{\boldsymbol{\lambda}}^{\circ})} \right ) \\
& \geq   a_{N}  \left ( \sum_{\ell = p_{\max}^\star(\tilde{\boldsymbol{\lambda}}^{\circ})+1}^j \frac{ E_{p_{\max}^\star(\tilde{\boldsymbol{\lambda}}^{\circ} )}(\tilde{\boldsymbol{\lambda}}^{\circ})}{E_{\ell}(\tilde{\boldsymbol{\lambda}}^{\circ})} + \frac{ E_{p_{\max}^\star(\tilde{\boldsymbol{\lambda}}^{\circ} )}(\tilde{\boldsymbol{\lambda}}^{\circ})}{\tilde{E}_{p_{\max}^\star(\tilde{\boldsymbol{\lambda}}^{\circ} )}(\tilde{\boldsymbol{\lambda}}^{\circ})} \right ) .
\end{align*}
Recall, that the $  E_{\ell}(\tilde{\boldsymbol{\lambda}}^{\circ})  $ are nonincreasing in $\ell$ and $  E_{p_{\max}^\star(\tilde{\boldsymbol{\lambda}}^{\circ} )}(\tilde{\boldsymbol{\lambda}}^{\circ}) \geq \tilde{E}_{p_{\max}^\star(\tilde{\boldsymbol{\lambda}}^{\circ} )}(\tilde{\boldsymbol{\lambda}}^{\circ})    $, see \eqref{eq:3.2}. Consequently, we find
\begin{align*}
E_{p_{\max}^\star(\tilde{\boldsymbol{\lambda}}^{\circ} )}(\tilde{\boldsymbol{\lambda}}^{\circ})
& \geq   a_{N}  \left ( \sum_{\ell = p_{\max}^\star(\tilde{\boldsymbol{\lambda}}^{\circ})+1}^j 1 + 1 \right ) =  a_{N}   (  j - p_{\max}^\star(\tilde{\boldsymbol{\lambda}}^{\circ}) + 1  )  .
\end{align*}
Using this in combination with $j=r(\mathcal{T}'_N, \tilde{\boldsymbol{\lambda}}^{\circ}  )-1 $ we obtain
\begin{equation} \label{eq:3.6}
	e_{p_{\max}^\star(\tilde{\boldsymbol{\lambda}}^{\circ} )}(\tilde{\boldsymbol{\lambda}}^{\circ} ) \ge E_{p_{\max}^\star( \tilde{\boldsymbol{\lambda}}^{\circ})}(\tilde{\boldsymbol{\lambda}}^{\circ} ) \ge a_N \max\{r( \mathcal{T}'_N,  \tilde{\boldsymbol{\lambda}}^{\circ}   ) - p_{\max}^\star( \tilde{\boldsymbol{\lambda}}^{\circ}   ),0\}.
	\end{equation}
In a next step we estimate  $r( \mathcal{T}'_N , \tilde{\boldsymbol{\lambda}}^{\circ}   ) - p_{\max}^\star( \tilde{\boldsymbol{\lambda}}^{\circ} )$ for all $  \tilde{\boldsymbol{\lambda}}^{\circ} \in \mathcal{V}(\mathcal{T}^\star)$. Recall, that $p_{\max}^\star(  \tilde{\boldsymbol{\lambda}}^{\circ} ) $ for $ \tilde{\boldsymbol{\lambda}}^{\circ} \in \mathcal{T}^\star \backslash \mathcal{V}(\mathcal{T}^\star)$ is determined by $P_{\max}^{\star}$. Therefore we get
\begin{equation}\label{eq:A4}
\sum_{  \tilde{\boldsymbol{\lambda}}^{\circ}  \in \mathcal{V}(\mathcal{T}^\star)} \max\{r( \mathcal{T}'_N , \tilde{\boldsymbol{\lambda}}^{\circ}  )- p_{\max}^\star(  \tilde{\boldsymbol{\lambda}}^{\circ}  ),0\}\ge\sum_{ \tilde{\boldsymbol{\lambda}}^{\circ}  \in \mathcal{V}(\mathcal{T}^\star \cap \mathcal{T}'_N)} r( \mathcal{T}'_N,  \tilde{\boldsymbol{\lambda}}^{\circ}  )- p_{\max}^\star(   \tilde{\boldsymbol{\lambda}}^{\circ} ). 
\end{equation}
To further estimate \eqref{eq:A4} on the one hand we observe
		\begin{align*}
	\sum_{ \tilde{\boldsymbol{\lambda}}^{\circ}  \in \mathcal{V}(\mathcal{T}^\star \cap \mathcal{T}'_N)} p_{\max}^\star( \tilde{\boldsymbol{\lambda}}^{\circ}  ) & \le \sum_{ \tilde{\boldsymbol{\lambda}}^{\circ} \in \mathcal{T}^\star} p_{\max}^\star( \tilde{\boldsymbol{\lambda}}^{\circ}  ) \\
	&  = \frac{2}{3} \sum_{ \tilde{\boldsymbol{\lambda}}^{\circ} \in \mathcal{T}^\star} \frac{3}{2} p_{\max}^\star( \tilde{\boldsymbol{\lambda}}^{\circ}  )   \\
	& \le \frac{2}{3}  \sum_{\tilde{\boldsymbol{\lambda}}^{\circ} \in \mathcal{T}^\star} \Big ( \frac{( p_{\textrm{max}}^{\star}(\tilde{\boldsymbol{\lambda}}^{\circ}) + 1)^2 + (  p_{\textrm{max}}^{\star}(\tilde{\boldsymbol{\lambda}}^{\circ}) + 1  )}{2} - 1 \Big )   \\
	& = \frac{2}{3} ( \# T^{\star}   - |T^\star| ).
	\end{align*}
Here in the last step we applied \eqref{eq_card_tree11}. Recall, that the bivariate quarklet tree $  T^{\star} $ has cardinality $n$, which means $ \# T^{\star} = n  $. Moreover, of course we have $ |\mathcal{T}^\star| = |  T^{\star}    |  $. Consequently, we find
\begin{equation}\label{eq_proof_plug_togh1}
\sum_{ \tilde{\boldsymbol{\lambda}}^{\circ}  \in \mathcal{V}(\mathcal{T}^\star \cap \mathcal{T}'_N)} p_{\max}^\star( \tilde{\boldsymbol{\lambda}}^{\circ}  )  \le \frac{2}{3} ( n   - |\mathcal{T}^\star| ).
\end{equation}
On the other hand, to deal with \eqref{eq:A4}, recall that the wavelet tree $ \mathcal{T}'_N  $ is resulting out of the algorithm {\bf BIVARIATE\textunderscore NEARBEST\textunderscore TREE} after $ N $ refinement steps. Moreover, we use that in each refinement step at least two children are added to the current wavelet tree. Hence we get
\begin{align*}
\sum_{\tilde{\boldsymbol{\lambda}}^{\circ} \in \mathcal{V}(\mathcal{T}^\star \cap \mathcal{T}'_N)} r( \mathcal{T}'_N, \tilde{\boldsymbol{\lambda}}^{\circ}  ) & = r( \mathcal{T}'_N, \mathcal{R}  ) - r (  \mathcal{T}^\star \cap \mathcal{T}'_N   ,  \mathcal{R}  ) \\
& = N - r (  \mathcal{T}^\star \cap \mathcal{T}'_N   ,  \mathcal{R}  ) \\
& \ge N - \frac{|\mathcal{T}^\star \cap \mathcal{T}'_N|-1}{2} \\
& \ge N - \frac{|\mathcal{T}^\star |-1}{2} .
\end{align*}
Thus it follows
\begin{equation}\label{eq_proof_plug_togh2}
\sum_{\tilde{\boldsymbol{\lambda}}^{\circ} \in \mathcal{V}(\mathcal{T}^\star \cap \mathcal{T}'_N)} r( \mathcal{T}'_N, \tilde{\boldsymbol{\lambda}}^{\circ}  ) \ge N - \frac{1}{2} |\mathcal{T}^\star | + \frac{1}{2} .
\end{equation}
Now a combination of \eqref{eq:A4} with \eqref{eq_proof_plug_togh1} and \eqref{eq_proof_plug_togh2} yields
\begin{align*}
\sum_{  \tilde{\boldsymbol{\lambda}}^{\circ}  \in \mathcal{V}(\mathcal{T}^\star)} \max\{r( \mathcal{T}'_N , \tilde{\boldsymbol{\lambda}}^{\circ}  )- p_{\max}^\star(  \tilde{\boldsymbol{\lambda}}^{\circ}  ),0\} & \ge \sum_{ \tilde{\boldsymbol{\lambda}}^{\circ}  \in \mathcal{V}(\mathcal{T}^\star \cap \mathcal{T}'_N)} r( \mathcal{T}'_N,  \tilde{\boldsymbol{\lambda}}^{\circ}  )- \sum_{ \tilde{\boldsymbol{\lambda}}^{\circ}  \in \mathcal{V}(\mathcal{T}^\star \cap \mathcal{T}'_N)} p_{\max}^\star(   \tilde{\boldsymbol{\lambda}}^{\circ} ) \\ 
& \ge  N - \frac{1}{2} |\mathcal{T}^\star | + \frac{1}{2}   -  \frac{2}{3} ( n   - |\mathcal{T}^\star| )  \\ 
& =  N  + \frac{1}{2}   -  \frac{2}{3}  n   + \frac{1}{6} |\mathcal{T}^\star|   \\
& \ge  N  + \frac{1}{2}   -  \frac{2}{3}  n   .
\end{align*}  
Finally, to complete the proof, this estimate in conjunction with \eqref{eq:split} and \eqref{eq:3.6} implies
	\begin{align*}
	\sigma_n &    \ge \sum_{\substack{ \tilde{\boldsymbol{\lambda}}^{\circ}  \in \mathcal{V}(\mathcal{T}^\star), \\ r( \mathcal{T}'_N, \tilde{\boldsymbol{\lambda}}^{\circ} ) > p_{\max}^\star( \tilde{\boldsymbol{\lambda}}^{\circ}  )}} e_{p_{\max}^\star( \tilde{\boldsymbol{\lambda}}^{\circ} )}( \tilde{\boldsymbol{\lambda}}^{\circ}  )   \\
	&   \ge \sum_{\substack{ \tilde{\boldsymbol{\lambda}}^{\circ}  \in \mathcal{V}(\mathcal{T}^\star), \\ r( \mathcal{T}'_N, \tilde{\boldsymbol{\lambda}}^{\circ} ) > p_{\max}^\star( \tilde{\boldsymbol{\lambda}}^{\circ}  )}}   a_N \max\{r( \mathcal{T}'_N,  \tilde{\boldsymbol{\lambda}}^{\circ}   ) - p_{\max}^\star( \tilde{\boldsymbol{\lambda}}^{\circ}   ),0\}  \\
	&   = a_{N} \sum_{ \tilde{\boldsymbol{\lambda}}^{\circ}  \in \mathcal{V}(\mathcal{T}^\star)}    \max\{r( \mathcal{T}'_N,  \tilde{\boldsymbol{\lambda}}^{\circ}   ) - p_{\max}^\star( \tilde{\boldsymbol{\lambda}}^{\circ}   ),0\}  \\
	&\ge a_N \Big (  N    -  \frac{2}{3}  n + \frac{1}{2} \Big ). 
	\end{align*}
This is the desired result.
\end{proof}

To continue we deduce an upper estimate for the global error using the parameter $a_N = a(\mathcal{R})$ for a given wavelet tree with root $ \mathcal{R}  $.

\begin{Lemma}\label{lemma:upper}
Let $  N \in \mathbb{N}  $ and let $T_N = (\mathcal{T}_N, \mathcal{T}'_N) $ be the bivariate quarklet tree created by the algorithm {\bf BIVARIATE\textunderscore  NEARBEST\textunderscore TREE} with a subsequent trimming. We define the parameter $a_N\coloneqq a(\mathcal{R})$ for the wavelet tree $\mathcal{T}'_N$ with root $ \mathcal{R} $. Then for the global error it holds
	\begin{equation}\label{lem_upper:estimate}
		\mathcal{E}(T_N) \le  a_N \left(3N+1\right).
	\end{equation}
\end{Lemma}

\begin{proof}
For the proof let $ \mathcal{T}'_N  $ be the bivariate wavelet tree produced by the algorithm {\bf BIVARIATE\textunderscore  NEARBEST\textunderscore TREE}. Let $L$ be the set of nodes $\tilde{\boldsymbol{\lambda}}^{\circ} \in \mathcal{T}'_N  $ for that $a(\tilde{\boldsymbol{\lambda}}^{\circ}) = \tilde{E}_{r( \mathcal{T}'_N,  \tilde{\boldsymbol{\lambda}}^{\circ}  )}( \tilde{\boldsymbol{\lambda}}^{\circ} )$ holds in \eqref{def_func_q1}. Furthermore, $Q$ is the maximal subtree of $\mathcal{T}'_N$ with root $\mathcal{R}$ for that we observe $L \cap Q = \mathcal{V}(Q)$. With other words the set $L$ does not contain any inner node of $Q$. We observe that $ \tilde{\boldsymbol{\lambda}}^{\circ}  \in \mathcal{V}(Q)  $ yields $ \tilde{\boldsymbol{\lambda}}^{\circ} \in L   $, which implies
\begin{equation}\label{eq:B3}
		\tilde{E}_{r( \mathcal{T}'_N, \tilde{\boldsymbol{\lambda}}^{\circ}   )}( \tilde{\boldsymbol{\lambda}}^{\circ}   ) = a(  \tilde{\boldsymbol{\lambda}}^{\circ}   ) \le  a(\mathcal{R}) = a_N.
\end{equation}
To continue let $ \tilde{\boldsymbol{\lambda}}^{\circ}  \in \mathcal{V}(Q)$ and $ r( \mathcal{T}'_N, \tilde{\boldsymbol{\lambda}}^{\circ}  ) = j$. We can apply \eqref{eq:3.3} to get
\begin{align*}
\frac{1}{\tilde{E}_j(\tilde{\boldsymbol{\lambda}}^{\circ})}  = \sum_{k=1}^j \frac{1}{E_k(\tilde{\boldsymbol{\lambda}}^{\circ})} + \sum_{ \tilde{\boldsymbol{\mu}}^{\circ} \in \mathcal{A}_{\mathcal{T}'_{N}}(\tilde{\boldsymbol{\lambda}}^{\circ}) } \frac{1}{e(\tilde{\boldsymbol{\mu}}^{\circ})} .
\end{align*}
Consequently, we also find
\begin{align*}
E_j(\tilde{\boldsymbol{\lambda}}^{\circ})   =  \tilde{E}_j(\tilde{\boldsymbol{\lambda}}^{\circ}) \left ( \sum_{k=1}^j \frac{E_j(\tilde{\boldsymbol{\lambda}}^{\circ})}{E_k(\tilde{\boldsymbol{\lambda}}^{\circ})} + \sum_{ \tilde{\boldsymbol{\mu}}^{\circ} \in \mathcal{A}_{\mathcal{T}'_{N}}(\tilde{\boldsymbol{\lambda}}^{\circ}) } \frac{E_j(\tilde{\boldsymbol{\lambda}}^{\circ})}{e(\tilde{\boldsymbol{\mu}}^{\circ})} \right )  .
\end{align*}
Recall, that the $  E_{k}(\tilde{\boldsymbol{\lambda}}^{\circ})  $ are nonincreasing in $k$. Using this in combination with \eqref{eq:B3}, we obtain 
\begin{align*}
E_j(\tilde{\boldsymbol{\lambda}}^{\circ})   \leq  a_{N} \left ( \sum_{k=1}^j 1 + \sum_{ \tilde{\boldsymbol{\mu}}^{\circ} \in \mathcal{A}_{\mathcal{T}'_{N}}(\tilde{\boldsymbol{\lambda}}^{\circ}) } \frac{E_j(\tilde{\boldsymbol{\lambda}}^{\circ})}{e(\tilde{\boldsymbol{\mu}}^{\circ})} \right )   = a_{N} \left ( j + \sum_{ \tilde{\boldsymbol{\mu}}^{\circ} \in \mathcal{A}_{\mathcal{T}'_{N}}(\tilde{\boldsymbol{\lambda}}^{\circ}) } \frac{E_j(\tilde{\boldsymbol{\lambda}}^{\circ})}{e(\tilde{\boldsymbol{\mu}}^{\circ})} \right ). 
\end{align*}
Moreover, an application of $ E_j( \tilde{\boldsymbol{\lambda}}^{\circ}  ) \leq E_0( \tilde{\boldsymbol{\lambda}}^{\circ}  ) = e( \tilde{\boldsymbol{\lambda}}^{\circ} )   $ implies
\begin{equation}\label{eq:A}
E_j(\tilde{\boldsymbol{\lambda}}^{\circ})   \leq a_{N} \left ( j + \sum_{ \tilde{\boldsymbol{\mu}}^{\circ} \in \mathcal{A}_{\mathcal{T}'_{N}}(\tilde{\boldsymbol{\lambda}}^{\circ}) } \frac{ e( \tilde{\boldsymbol{\lambda}}^{\circ} )  }{e(\tilde{\boldsymbol{\mu}}^{\circ})} \right ). 
\end{equation}
To continue we distinguish two different cases for a leaf $ \tilde{\boldsymbol{\lambda}}^{\circ}  \in \mathcal{V}(Q)$. First we consider the case that there exists a $ \tilde{\boldsymbol{\nu}}^{\circ} \in \mathcal{V}(\mathcal{T}_N)$, such that $ \tilde{\boldsymbol{\lambda}}^{\circ}  \preceq  \tilde{\boldsymbol{\nu}}^{\circ} $. Here $  \mathcal{T}_N $ is the subtree of $ \mathcal{T}'_{N}  $ resulting out of the trimming process. Let $  \mathcal{T}'_{\tilde{\boldsymbol{\lambda}}^{\circ}}  $ be the maximal subtree of $ \mathcal{T}'_N  $ with root $ \tilde{\boldsymbol{\lambda}}^{\circ}  $. Using \eqref{eq:3.1} and the characteristic property of the trimmed tree $ \mathcal{T}_N  $ we find
	\begin{equation}\label{eq:desc}
		E_{r( \mathcal{T}'_N, \tilde{\boldsymbol{\lambda}}^{\circ}   )}( \tilde{\boldsymbol{\lambda}}^{\circ}  ) = \sum_{\tilde{\boldsymbol{\eta}}^{\circ} \in \mathcal{V}(\mathcal{T}'_{\tilde{\boldsymbol{\lambda}}^{\circ}} \cap \mathcal{T}_N)} E_{r( \mathcal{T}'_N,  \tilde{\boldsymbol{\eta}}^{\circ}  )}( \tilde{\boldsymbol{\eta}}^{\circ} ) = \sum_{\tilde{\boldsymbol{\eta}}^{\circ} \in \mathcal{V}(\mathcal{T}'_{\tilde{\boldsymbol{\lambda}}^{\circ}} \cap \mathcal{T}_N)} e_{r( \mathcal{T}'_N,  \tilde{\boldsymbol{\eta}}^{\circ}  )}(\tilde{\boldsymbol{\eta}}^{\circ}).
\end{equation}
The second case is that there is a $\tilde{\boldsymbol{\nu}}^{\circ} \in \mathcal{V}(\mathcal{T}_N)$ with $\tilde{\boldsymbol{\lambda}}^{\circ} \succ \tilde{\boldsymbol{\nu}}^{\circ}$. Then the characteristic property of the trimmed tree $ \mathcal{T}_N   $ and an application of \eqref{eq:3.1} yields
\begin{equation}\label{eq:prec}
e_{r( \mathcal{T}'_N,  \tilde{\boldsymbol{\nu}}^{\circ}  )}( \tilde{\boldsymbol{\nu}}^{\circ}  ) = E_{r( \mathcal{T}'_N , \tilde{\boldsymbol{\nu}}^{\circ}   )}(  \tilde{\boldsymbol{\nu}}^{\circ} ) \le  \sum_{  \tilde{\boldsymbol{\eta}}^{\circ}  \in \mathcal{V}(\mathcal{T}'_{\tilde{\boldsymbol{\nu}}^{\circ}} \cap Q)} E_{r( \mathcal{T}'_N , \tilde{\boldsymbol{\eta}}^{\circ}  )}(  \tilde{\boldsymbol{\eta}}^{\circ} ). 
\end{equation}
Here $  \mathcal{T}'_{\tilde{\boldsymbol{\nu}}^{\circ}}  $ is the maximal subtree of $ \mathcal{T}'_N  $ with root $ \tilde{\boldsymbol{\nu}}^{\circ}  $. To continue we can divide the set of leaves $\mathcal{V}(\mathcal{T}_N)$ into two groups according to the two cases explained above. Then we can use the definition of the global error, see \eqref{eq_glob_err}, and afterwards our observations \eqref{eq:desc} as well as  \eqref{eq:prec} to find 
\begin{equation}\label{eq:B_start_prof_B}
		\mathcal{E}(T_N) = \sum_{ \tilde{\boldsymbol{\lambda}}^{\circ}  \in \mathcal{V}(\mathcal{T}_N)} e_{r( \mathcal{T}'_N , \tilde{\boldsymbol{\lambda}}^{\circ}   )}( \tilde{\boldsymbol{\lambda}}^{\circ}   ) \le \sum_{  \tilde{\boldsymbol{\lambda}}^{\circ} \in \mathcal{V}(Q)} E_{r( \mathcal{T}'_N, \tilde{\boldsymbol{\lambda}}^{\circ}   )}( \tilde{\boldsymbol{\lambda}}^{\circ}  ).
\end{equation}
For the previous step it was essential that the bivariate quarklet tree $  T_{N}  $ is defined by $   T_N = (\mathcal{T}_N, \mathcal{T}'_N)$. To further estimate \eqref{eq:B_start_prof_B} we plug in \eqref{eq:A} with  $ j = r( \mathcal{T}'_N, \tilde{\boldsymbol{\lambda}}^{\circ}  ) $  and obtain
\begin{align}
\mathcal{E}(T_N) &\le  \sum_{  \tilde{\boldsymbol{\lambda}}^{\circ} \in \mathcal{V}(Q)} E_{r( \mathcal{T}'_N, \tilde{\boldsymbol{\lambda}}^{\circ}   )}( \tilde{\boldsymbol{\lambda}}^{\circ}  )  \notag \\
& \le \sum_{  \tilde{\boldsymbol{\lambda}}^{\circ} \in \mathcal{V}(Q)} a_N \left(  r( \mathcal{T}'_N, \tilde{\boldsymbol{\lambda}}^{\circ}  )   +  \sum_{ \tilde{\boldsymbol{\mu}}^{\circ} \in \mathcal{A}_{\mathcal{T}_{N}'}(\tilde{\boldsymbol{\lambda}}^{\circ}) }   \frac{e(  \tilde{\boldsymbol{\lambda}}^{\circ} )}{e( \tilde{\boldsymbol{\mu}}^{\circ}  )} \right) \notag \\
& =   a_N \left(  \sum_{  \tilde{\boldsymbol{\lambda}}^{\circ} \in \mathcal{V}(Q)}   r( \mathcal{T}'_N, \tilde{\boldsymbol{\lambda}}^{\circ}  )   + \sum_{  \tilde{\boldsymbol{\lambda}}^{\circ} \in \mathcal{V}(Q)}  \sum_{ \tilde{\boldsymbol{\mu}}^{\circ} \in \mathcal{A}_{\mathcal{T}_{N}'}(\tilde{\boldsymbol{\lambda}}^{\circ}) }   \frac{e(  \tilde{\boldsymbol{\lambda}}^{\circ} )}{e( \tilde{\boldsymbol{\mu}}^{\circ}  )} \right) \notag \\
& =   a_N \left(  \sum_{  \tilde{\boldsymbol{\lambda}}^{\circ} \in \mathcal{V}(Q)}   r( \mathcal{T}'_N, \tilde{\boldsymbol{\lambda}}^{\circ}  )   +  \sum_{\tilde{\boldsymbol{\mu}}^{\circ} \in  \mathcal{V}(Q)} \frac{  e(  \tilde{\boldsymbol{\mu}}^{\circ} )}{e( \tilde{\boldsymbol{\mu}}^{\circ}  )}   +  \sum_{\tilde{\boldsymbol{\mu}}^{\circ} \in (Q \backslash \mathcal{V}(Q))} \frac{ \sum\limits_{ \tilde{\boldsymbol{\lambda}}^{\circ} \in \mathcal{V}(\mathcal{D}_{\mathcal{T}_{N}'}(\tilde{\boldsymbol{\mu}}^{\circ}) \cap Q)}  e(  \tilde{\boldsymbol{\lambda}}^{\circ} )}{e( \tilde{\boldsymbol{\mu}}^{\circ}  )} \right) \notag \\
& =   a_N \left(  \sum_{  \tilde{\boldsymbol{\lambda}}^{\circ} \in \mathcal{V}(Q)}   r( \mathcal{T}'_N, \tilde{\boldsymbol{\lambda}}^{\circ}  )   +  \sum_{\tilde{\boldsymbol{\mu}}^{\circ} \in  \mathcal{V}(Q)} 1   +  \sum_{\tilde{\boldsymbol{\mu}}^{\circ} \in (Q \backslash \mathcal{V}(Q))} \frac{ \sum\limits_{ \tilde{\boldsymbol{\lambda}}^{\circ} \in \mathcal{V}(\mathcal{D}_{\mathcal{T}_{N}'}(\tilde{\boldsymbol{\mu}}^{\circ}) \cap Q)}  e(  \tilde{\boldsymbol{\lambda}}^{\circ} )}{e( \tilde{\boldsymbol{\mu}}^{\circ}  )} \right) ,
\label{eq:global_err_est}  
\end{align}
whereby in the penultimate step we changed the order of summation. Here $  \mathcal{D}_{\mathcal{T}_{N}'}(\tilde{\boldsymbol{\mu}}^{\circ})  $ refers to the set of all descendants of $ \tilde{\boldsymbol{\mu}}^{\circ}  $ including $ \tilde{\boldsymbol{\mu}}^{\circ}  $ itself according to the bivariate wavelet tree $  \mathcal{T}_{N}'  $. Recall, that by the subadditivity for the error of the lowest order, see $\eqref{eq:1.2}$, we find
\begin{equation*} \label{eq:weakSub}
e(\tilde{\boldsymbol{\lambda}}^{\circ}) \ge  \sum_{\tilde{\boldsymbol{\eta}}^{\circ} \in \mathcal{V}( \mathcal{D}_{\mathcal{T}}(\tilde{\boldsymbol{\lambda}}^{\circ}) )} e(\tilde{\boldsymbol{\eta}}^{\circ})
\end{equation*} 
for all bivariate wavelet trees $\mathcal{T}$ as given in Definition \ref{def_wav_tre}. Using this to further estimate the last sum in \eqref{eq:global_err_est}, we obtain
\begin{align}
\mathcal{E}(T_N) &\le   a_N \left(  \sum_{  \tilde{\boldsymbol{\lambda}}^{\circ} \in \mathcal{V}(Q)}  ( r( \mathcal{T}'_N, \tilde{\boldsymbol{\lambda}}^{\circ}  )  +1 ) +  \sum_{\tilde{\boldsymbol{\mu}}^{\circ} \in (Q \backslash \mathcal{V}(Q))} \frac{ \sum\limits_{ \tilde{\boldsymbol{\lambda}}^{\circ} \in \mathcal{V}(\mathcal{D}_{\mathcal{T}_{N}'}(\tilde{\boldsymbol{\mu}}^{\circ}) \cap Q)}  e(  \tilde{\boldsymbol{\lambda}}^{\circ} )}{e( \tilde{\boldsymbol{\mu}}^{\circ}  )} \right) \notag \\
&  \le   a_N \left(  \sum_{  \tilde{\boldsymbol{\lambda}}^{\circ} \in \mathcal{V}(Q)}   r( \mathcal{T}'_N, \tilde{\boldsymbol{\lambda}}^{\circ}  )  + \sum_{  \tilde{\boldsymbol{\lambda}}^{\circ} \in \mathcal{V}(Q)}  1  +  \sum_{\tilde{\boldsymbol{\mu}}^{\circ} \in (Q \backslash \mathcal{V}(Q))} 1 \right)  \notag \\
&  =   a_N \left(  \sum_{  \tilde{\boldsymbol{\lambda}}^{\circ} \in \mathcal{V}(Q)}   r( \mathcal{T}'_N, \tilde{\boldsymbol{\lambda}}^{\circ}  )  +  | \mathcal{V}(Q) |  + | Q \backslash \mathcal{V}(Q) | \right) .
\label{eq:A2}  
\end{align}
When we apply that $Q$ is a subtree of $\mathcal{T}'_N$, and that in each refinement step either two or three children are added, we see, that the number of nodes in $\mathcal{T}'_N$ can be estimated by
	\begin{equation*}\label{eq:B2}
		|\mathcal{T}'_N| \geq |Q| + 2\sum_{\tilde{\boldsymbol{\lambda}}^{\circ} \in \mathcal{V}(Q)} r( \mathcal{T}'_N  , \tilde{\boldsymbol{\lambda}}^{\circ}  ).
	\end{equation*}
Hence we also find	
\begin{align*}
\sum_{\tilde{\boldsymbol{\lambda}}^{\circ} \in \mathcal{V}(Q)} r( \mathcal{T}'_N  , \tilde{\boldsymbol{\lambda}}^{\circ}  )  \leq \frac{ |\mathcal{T}'_N| - |Q|}{2} .
\end{align*}
In combination with \eqref{eq:A2} this yields
\begin{align*}
\mathcal{E}(T_N) &  \leq   a_N \left(  \sum_{  \tilde{\boldsymbol{\lambda}}^{\circ} \in \mathcal{V}(Q)}   r( \mathcal{T}'_N, \tilde{\boldsymbol{\lambda}}^{\circ}  )  +  | \mathcal{V}(Q) |  + | Q \backslash \mathcal{V}(Q) | \right) \notag \\
& =  a_N \left(  \sum_{  \tilde{\boldsymbol{\lambda}}^{\circ} \in \mathcal{V}(Q)}   r( \mathcal{T}'_N, \tilde{\boldsymbol{\lambda}}^{\circ}  )  +  | Q | \right) \notag \\
& \leq  a_N \left(  \frac{ |\mathcal{T}'_N| - |Q|}{2}  +  | Q | \right) \notag \\
& =  a_N \left(  \frac{ |\mathcal{T}'_N| + |Q|}{2} \right) . 
\end{align*}
To continue again we use that $ Q $ is a subtree of $ \mathcal{T}'_N   $. Moreover, recall that $ \mathcal{T}'_N  $ is a wavelet tree resulting out of the algorithm {\bf BIVARIATE\textunderscore  NEARBEST\textunderscore TREE} after $N$ refinement steps. We know, that in each refinement step either two or three children are added to the current wavelet tree. Consequently we observe
\begin{align*}
 |Q| \le |\mathcal{T}'_N| \le  3N+1 .
\end{align*} 
Hence we get 
\begin{align*}
\mathcal{E}(T_N)   \leq a_N \left(  \frac{ |\mathcal{T}'_N| + |Q|}{2} \right) \leq a_N  |\mathcal{T}'_N| \leq a_{N} (3 N + 1) .
\end{align*}
So the proof of \eqref{lem_upper:estimate} is complete.
\end{proof}

Now we have all tools at hand to show that the algorithm  {\bf BIVARIATE\textunderscore  NEARBEST\textunderscore TREE} assembles a bivariate quarklet tree which has the property \eqref{eq:nearbest} and therefore can be seen as near-best.

\begin{Theorem}\label{theorem:1}
Let $ n,N \in \mathbb{N} $ with $n\le N$ and let $ e_{p}(\tilde{\boldsymbol{\lambda}})  $ be local errors that fulfill \eqref{eq:1.2} and \eqref{eq:1.3}. Then the algorithm {\bf BIVARIATE\textunderscore  NEARBEST\textunderscore TREE} with a subsequent trimming provides a bivariate quarklet tree $T_N = (\mathcal{T}_N, \mathcal{T}'_N)$ such that the corresponding approximation in terms of bivariate tensor quarklets is near-best in the sense
	\begin{equation}\label{eq:result1}
		\mathcal{E}(T_N) \le \frac{3N+1}{ N    -  \frac{2}{3}  n + \frac{1}{2}}\sigma_n.
	\end{equation}
\end{Theorem}

\begin{proof}
For the proof at first we recall that Lemma \ref{lemma:lower1} yields
\begin{align*}
a_{N} \leq \frac{\sigma_{n}}{ N    -  \frac{2}{3}  n + \frac{1}{2}}.
\end{align*}
In combination with Lemma \ref{lemma:upper} we find
\begin{align*}
\mathcal{E}(T_N) \le  a_N ( 3N+1  ) \le \frac{3N+1}{ N    -  \frac{2}{3}  n + \frac{1}{2}}\sigma_n .
\end{align*}
The proof is complete.
\end{proof}

\begin{Remark}
We can use Theorem \ref{theorem:1} in order to see \eqref{eq:nearbest}. For that purpose let $ M \in \mathbb{N}  $. Then Lemma \ref{lemma:T_N_cardinality} yields that we can run $ N = (\frac{15}{163})^{1/3} M^{1/3}     $ steps of the algorithm {\bf BIVARIATE\textunderscore NEARBEST\textunderscore TREE} while guaranteeing that $   \# T_N \le  M   $ for the resulting bivariate quarklet tree. Now let $  n = \frac{N}{2}   $. Then the constant showing up in \eqref{eq:result1} can be estimated by 
\begin{align*}
\frac{3N+1}{ N    -  \frac{2}{3}  n + \frac{1}{2}} = \frac{3  (\frac{15}{163})^{1/3} M^{1/3}     +1}{  (\frac{15}{163})^{1/3} M^{1/3}     -  \frac{1}{3}  (\frac{15}{163})^{1/3} M^{1/3}    + \frac{1}{2}}  \leq 10. 
\end{align*}
Consequently, \eqref{eq:result1} becomes 
\begin{align*}
		\mathcal{E}(T_N) \le \frac{3N+1}{ N    -  \frac{2}{3}  n + \frac{1}{2}}\sigma_n   \leq  10 \sigma_{\frac{1}{2}  (\frac{15}{163})^{1/3} M^{1/3} }  .
	\end{align*}
This also implies 
\begin{align*}
		\mathcal{E}(T_N)  \leq  C \sigma_{c M^{1/3} }  
	\end{align*}
with independent constants $  C > 1   $ and $  c \in (0,1]    $. Hence, \eqref{eq:nearbest} follows. 
\end{Remark}

\section{Approximation in \texorpdfstring{$\ell_2$}{l2}} \label{sec_practice}

Below we describe how we can use the algorithm {\bf BIVARIATE\textunderscore NEARBEST\textunderscore TREE} to approximate functions $ f \in L_2((0,1)^2)  $ via bivariate tensor quarklets. For that purpose in a first step we see that each bivariate quarklet tree $ T $ refers to a function $ f_{T} \in L_{2}((0,1)^2)   $ which consists of a sum of bivariate quarklets only. 

\subsection{Bivariate Quarklet Trees and Functions in $L_{2}((0,1)^2)$}\label{sec:coeff_approach}

Let $ f \in L_{2}((0,1)^2)    $ be given. In Section \ref{subsec_biv_tensorwavelets} we have seen that there exists at least one sequence $ \{  c_{\boldsymbol{\lambda}}  \}_{\boldsymbol{\lambda} \in  \boldsymbol{\nabla} } \in \ell_{2}( \boldsymbol{\nabla} ) $ such that 
\begin{equation}\label{sec_how_biv_equ1}
f = \sum_{ \boldsymbol{\lambda} \in  \boldsymbol{\nabla}  }   c_{\boldsymbol{\lambda}} w_{\boldsymbol{\lambda}}^{-1}  \boldsymbol{\psi_{\lambda}} , 
\end{equation} 
see \eqref{eq_canon_dual_frame_2_vorl1}. Recall, that each bivariate quarklet tree consists of enhanced bivariate quarklet indices $  \tilde{\boldsymbol{\lambda}} = ( (p_{1},j_{1},k_{1}), (p_{2},j_{2},k_{2}), \alpha ) \in \tilde{\boldsymbol{\Lambda}}     $. Here $  \alpha \in \{ 0, 1, 2 \}   $ refers to the refinement options as described in Section \ref{subsec_wav_tree}. In order to incorporate the different refinement options to the quarklet indices collected in the set $  \boldsymbol{\nabla}    $ we define
\begin{equation}\label{eq_extend_frame_tilde1}
\boldsymbol{\tilde{\nabla}} := \{ \tilde{\boldsymbol{\lambda}} = ( (p_{1},j_{1},k_{1}), (p_{2},j_{2},k_{2}), \alpha )  :      \boldsymbol{\lambda} = ( (p_{1},j_{1},k_{1}), (p_{2},j_{2},k_{2}) ) \in  \boldsymbol{\nabla}, \alpha \in \{  0, 1, 2 \} \}.
\end{equation}
Moreover, we put
\begin{equation}\label{eq_extend_frame_tilde2}
\boldsymbol{\psi_{\tilde{\lambda}}} := \boldsymbol{\psi_{\lambda}} ,
\end{equation}
see \eqref{def_tensor_quarklet}. That means, enhanced bivariate quarklet indices that only differ in the refinement option parameter $\alpha$ refer to the same bivariate tensor quarklet. Consequently, \eqref{sec_how_biv_equ1} implies that there always exists a sequence $ \{  c_{\tilde{\boldsymbol{\lambda}}}  \}_{\tilde{\boldsymbol{\lambda}} \in  \boldsymbol{\tilde{\nabla}} } \in \ell_{2}( \boldsymbol{\tilde{\nabla}} ) $ such that 
\begin{equation}\label{eq_canon_dual_frame_2}
f = \sum_{ \tilde{\boldsymbol{\lambda}} \in  \boldsymbol{\tilde{\nabla}}  }   c_{\tilde{\boldsymbol{\lambda}}} w_{\tilde{\boldsymbol{\lambda}}}^{-1}  \boldsymbol{\psi_{\tilde{\lambda}}} . 
\end{equation}
Notice, that the representation given in \eqref{eq_canon_dual_frame_2} is not unique. Now let a bivariate quarklet tree $ T = ( \mathcal{T} , \mathcal{T}'   )   $ produced by the algorithm {\bf BIVARIATE\textunderscore NEARBEST\textunderscore TREE} and a subsequent trimming be given. Let $\mathcal{R} :=  ((j_{1},k_{1}), (j_{2},k_{2}), \alpha )  $ be the root of this tree. Then for the case $j_{1} < j_{0}$ or $ j_{2} < j_{0}    $ the corresponding quarklet tree also contains some artificial nodes, which do not refer to bivariate tensor quarklets that are frame elements, but stand for zero functions, see \eqref{eq_zerowav_put}. Consequently, in what follows we want to transform the quarklet tree $T$ into a slightly modified index set $\tilde{T}$ which fulfills $ \tilde{T} \subset \boldsymbol{\tilde{\nabla}}   $. For that purpose we use an idea from \cite{VoDiss}, see equation (4.39) in Chapter 4.5.1. We put
\begin{small}
\begin{align*}
\tilde{T} & :=  \{ ((p_{1},j_{0}-1,k_{1}), (p_{2},j_{2},k_{2}), \alpha ) \in        \boldsymbol{\tilde{\nabla}} : \tilde{\boldsymbol{\lambda}} = ((0,j_{0},\ell_{k_{1}}), (0,j_{2},k_{2}), \alpha ) \in T, p_{1} + p_{2} \leq p_{\max}(\tilde{\boldsymbol{\lambda}}^{\circ})   \} \\
  & \cup \{ ((p_{1},j_{1},k_{1}), (p_{2},j_{0}-1,k_{2}), \alpha ) \in        \boldsymbol{\tilde{\nabla}} : \tilde{\boldsymbol{\lambda}} = ((0,j_{1},k_{1}), (0,j_{0},\ell_{k_{2}}), \alpha ) \in T, p_{1} + p_{2} \leq p_{\max}(\tilde{\boldsymbol{\lambda}}^{\circ})   \} \\
  & \cup \{ ((p_{1},j_{1},k_{1}), (p_{2},j_{2},k_{2}), \alpha ) \in T : j_{1} \geq j_{0}, j_{2} \geq j_{0}   \} .
\end{align*}
\end{small}
In other words, we keep the bivariate quarklet indices from the intersection of $T$ with $ \boldsymbol{\tilde{\nabla}}   $  and add the nodes referring to functions consisting of generators up to the polynomial degree of the assigned wavelet node. We observe $ \# \tilde{T} \leq C \# T    $. Here the constant depends on the maximal number of functions consisting of generators assigned to a single wavelet node on level $j_{1}=j_{0}$ or $j_{2}=j_{0}$. Now let $ f \in L_{2}((0,1)^2)    $ in the form \eqref{eq_canon_dual_frame_2} be given. Moreover, let $T$ be the bivariate quarklet tree resulting out of the algorithm {\bf BIVARIATE\textunderscore NEARBEST\textunderscore TREE} and $   \tilde{T}$ its modified version as described above. Then the quarklet tree approximation $ f_{T}   $ of $ f $ is defined as 
\begin{equation}\label{eq_func_quarktree_approx_mod}
f_{T} = \sum_{ \tilde{\boldsymbol{\lambda}} \in \tilde{T}  }     c_{\tilde{\boldsymbol{\lambda}}} w_{\tilde{\boldsymbol{\lambda}}}^{-1}  \boldsymbol{\psi_{\tilde{\lambda}}} . 
\end{equation}

\subsection{An Approach to Local Errors for Functions in $L_{2}((0,1)^{2})$}

In what follows we want to apply the algorithm {\bf BIVARIATE\textunderscore NEARBEST\textunderscore TREE} to functions $  f \in L_{2}((0,1)^{2})    $. For that purpose we have to find a precise definition for the local errors $  e_{p_{\max}}(\tilde{\boldsymbol{\lambda}})  $. Recall, that in Section \ref{subsec_localerr} we already identified some properties that the local errors necessarily have to fulfill, see \eqref{eq:1.2} and \eqref{eq:1.3}. Hence there are several restrictions we have to consider when looking for a possible definition for the local errors. At first let us recall that each function $ f \in L_{2}((0,1)^2)    $ can be written in the form \eqref{eq_canon_dual_frame_2} with a coefficient sequence $   \{  c_{\tilde{\boldsymbol{\lambda}}}  \}_{\tilde{\boldsymbol{\lambda}} \in  \boldsymbol{\tilde{\nabla}} } \in \ell_{2}( \boldsymbol{\tilde{\nabla}} )  $. For each (enhanced) quarklet index $ \tilde{\boldsymbol{\lambda}} = ( (p_{1},j_{1},k_{1}), (p_{2},j_{2},k_{2}), \alpha ) \in \tilde{\boldsymbol{\Lambda}}   $ we define a number $  d_{\tilde{\boldsymbol{\lambda}}}   $ by
\begin{equation}\label{eq_def_d_lambda}
d_{\tilde{\boldsymbol{\lambda}}}^{2}  := \left\{ \begin{array}{lll}
 0   & \quad & \mbox{for} \qquad j_{1} < j_{0} \; \mbox{or/and} \; j_{2} < j_{0} ;
\\  
c_{\tilde{\boldsymbol{\lambda}}}^{2} + \sum_{\ell_{1} \in \square_{j_{0}, k_{1} }} c^{2}_{( (p_{1},j_{0}-1,\ell_{1}), (p_{2},j_{2},k_{2}), \alpha )} & \quad & \mbox{for}\qquad j_{1} = j_{0} \; \mbox{and} \; j_{2} > j_{0}  ;
\\  
c_{\tilde{\boldsymbol{\lambda}}}^{2} + \sum_{\ell_{2} \in \square_{j_{0}, k_{2} }} c^{2}_{( (p_{1},j_{1},k_{1}), (p_{2},j_{0}-1,\ell_{2}), \alpha )} & \quad & \mbox{for}\qquad j_{2} = j_{0} \; \mbox{and} \; j_{1} > j_{0}  ;
\\
c_{\tilde{\boldsymbol{\lambda}}}^{2} & \quad & \mbox{for}\qquad j_{1} > j_{0} \; \mbox{and} \; j_{2} > j_{0} .
\\
\end{array}
\right.
\end{equation}
We can collect these numbers in a sequence $   \{  d_{\tilde{\boldsymbol{\lambda}}}  \}_{\tilde{\boldsymbol{\lambda}} \in  \tilde{\boldsymbol{\Lambda}}  } \in \ell_{2}( \tilde{\boldsymbol{\Lambda}}  )     $. Now we are well-prepared to give a precise definition for the local errors $   e_{p_{\max}}(\tilde{\boldsymbol{\lambda}}^{\circ})  $.
\begin{Definition}\label{def_locerr_L2_new}
	Let $ f \in  L_2((0,1)^2) $ be given in the form \eqref{eq_canon_dual_frame_2}. Then for each node $ \tilde{\boldsymbol{\lambda}}  \in \tilde{\boldsymbol{\Lambda}}_{0}  $ and $p_{\max} \in \mathbb{N}_0$ we define the local errors $  e_{p_{\max}}(\tilde{\boldsymbol{\lambda}}) $ via 
	\begin{align*}
	e_{p_{\max}}(  \tilde{\boldsymbol{\lambda}}   ) & \coloneqq \sum_{((j_{1}, k_{1}),(j_{2}, k_{2}), \alpha) \in \Upsilon(  \tilde{\boldsymbol{\lambda}}  )} \sum_{p_{1} + p_{2} > p_{\max}} \vert d_{((p_{1},j_{1},k_{1}),(p_{2},j_{2},k_{2}), \alpha)} \vert^2 \\
	& \qquad \qquad  + \sum_{\substack{((j_{1},k_{1}),(j_{2},k_{2}), \alpha )  \in  \mathcal{J}_{\tilde{\boldsymbol{\lambda}}} \\ ((j_{1},k_{1}),(j_{2},k_{2}),\alpha) \not = \tilde{\boldsymbol{\lambda}}  }} \sum_{p_{1} + p_{2} \ge 0} \vert d_{((p_{1},j_{1},k_{1}),(p_{2},j_{2},k_{2}), \alpha)} \vert^2.
	\end{align*}
\end{Definition}
Recall, that the set $ \Upsilon(  \tilde{\boldsymbol{\lambda}}  )   $ was defined in Section \ref{subsec_biv_qua_tre}. A definition for $ \mathcal{J}_{\tilde{\boldsymbol{\lambda}}}   $    can be found in Section \ref{subsec_wav_tree}. Having a closer look at Definition \ref{def_locerr_L2_new} it turns out that the first sum refers to $ \tilde{\boldsymbol{\lambda}}  $ and a subset of its ancestors, whereby the cumulated polynomial degree is greater than the maximal degree $p_{\max}$. The second sum gathers the descendants of $ \tilde{\boldsymbol{\lambda}}  $ according to $  \mathcal{J}_{\tilde{\boldsymbol{\lambda}}}   $ for all possible polynomial degrees.

\begin{Remark}
The local errors given in Definition \ref{def_locerr_L2_new} are inspired by Definition 4.2 in \cite{DaHoRaVo}. There local errors in the context of univariate quarklet tree approximation for functions $ f \in L_{2}((0,1))   $ have been established. Similar approaches already have been applied successfully for tree approximation using adaptive wavelet schemes, see \cite{bib:CDD03} and \cite{kappei2011adaptive}. 
\end{Remark}

To continue we verify that for the local errors given in Definition \ref{def_locerr_L2_new} the conditions \eqref{eq:1.2} and \eqref{eq:1.3} are fulfilled.

\begin{Lemma}\label{lem_loc_err_prop1_new}
Let $ f \in L_2((0,1)^2)  $ be given in the form \eqref{eq_canon_dual_frame_2}. Let $ \tilde{\boldsymbol{\lambda}}  \in \tilde{\boldsymbol{\Lambda}}_{0}    $ and $p_{\max} \in \mathbb{N}_0$.  Then the local errors $ e_{p_{\max}}( \tilde{\boldsymbol{\lambda}}  )  $ formulated in Definition \ref{def_locerr_L2_new} satisfy the properties \eqref{eq:1.2} and \eqref{eq:1.3}.
\end{Lemma}

\begin{proof}
For the proof we assume that $ \tilde{\boldsymbol{\lambda}}  \in \tilde{\boldsymbol{\Lambda}}_{0}    $ has exactly three children. The case that $ \tilde{\boldsymbol{\lambda}} $ has two children can be treated with similar methods. 
 
\textit{Step 1.}
	At first we prove \eqref{eq:1.2}. Therefore let $p_{\max}=0$ and $ \tilde{\boldsymbol{\eta}}_1, \tilde{\boldsymbol{\eta}}_2, \tilde{\boldsymbol{\eta}}_3 $ be the children of $ \tilde{\boldsymbol{\lambda}} \in \tilde{\boldsymbol{\Lambda}}_{0}  $ according to the refinement strategies (LSR.a.1) or (LSR.a.2), whereby the condition (UPC) is fulfilled. Recall, that by Section \ref{subsec_wav_tree} we find
\begin{align*}
( \mathcal{J}_{\tilde{\boldsymbol{\lambda}}} \setminus  \{   \tilde{\boldsymbol{\lambda}} \} ) \supset \{ \tilde{\boldsymbol{\eta}}_{1} , \tilde{\boldsymbol{\eta}}_{2}, \tilde{\boldsymbol{\eta}}_{3}  \} \cup (\mathcal{J}_{\tilde{\boldsymbol{\eta}}_{1}} \setminus  \{ \tilde{\boldsymbol{\eta}}_{1}   \}   ) \cup (\mathcal{J}_{\tilde{\boldsymbol{\eta}}_{2}} \setminus  \{ \tilde{\boldsymbol{\eta}}_{2}   \}   ) \cup (\mathcal{J}_{\tilde{\boldsymbol{\eta}}_{3}} \setminus  \{ \tilde{\boldsymbol{\eta}}_{3}   \}   ) .
\end{align*}
Using this in combination with Definition \ref{def_locerr_L2_new} for the local error of the lowest order concerning $ \tilde{\boldsymbol{\lambda}}  $ we observe
\begin{align*}
	e_{0}(  \tilde{\boldsymbol{\lambda}}   ) &  \geq \sum_{((j_{1}, k_{1}),(j_{2}, k_{2}),\alpha) \in \Upsilon(  \tilde{\boldsymbol{\lambda}}  )} \sum_{p_{1} + p_{2} > 0} \vert d_{((p_{1},j_{1},k_{1}),(p_{2},j_{2},k_{2}),\alpha)} \vert^2   \\
	& \qquad \qquad +  \sum_{((j_{1}, k_{1}),(j_{2}, k_{2}),\alpha) \in \{ \tilde{\boldsymbol{\eta}}_{1} , \tilde{\boldsymbol{\eta}}_{2}, \tilde{\boldsymbol{\eta}}_{3}  \}} \sum_{p_{1} + p_{2} > 0} \vert d_{((p_{1},j_{1},k_{1}),(p_{2},j_{2},k_{2}),\alpha)} \vert^2      \\
	& \qquad \qquad   + \sum_{\substack{((j_{1},k_{1}),(j_{2},k_{2}),\alpha)  \in  \mathcal{J}_{\tilde{\boldsymbol{\eta}}_{1}} \\ ((j_{1},k_{1}),(j_{2},k_{2}),\alpha) \not = \tilde{\boldsymbol{\eta}}_{1}  }} \sum_{p_{1} + p_{2} \ge 0} \vert d_{((p_{1},j_{1},k_{1}),(p_{2},j_{2},k_{2}),\alpha)} \vert^2 \\
	& \qquad \qquad   + \sum_{\substack{((j_{1},k_{1}),(j_{2},k_{2}),\alpha)  \in  \mathcal{J}_{\tilde{\boldsymbol{\eta}}_{2}} \\ ((j_{1},k_{1}),(j_{2},k_{2}),\alpha) \not = \tilde{\boldsymbol{\eta}}_{2}  }} \sum_{p_{1} + p_{2} \ge 0} \vert d_{((p_{1},j_{1},k_{1}),(p_{2},j_{2},k_{2}),\alpha)} \vert^2 \\
	& \qquad \qquad   + \sum_{\substack{((j_{1},k_{1}),(j_{2},k_{2}),\alpha)  \in  \mathcal{J}_{\tilde{\boldsymbol{\eta}}_{3}} \\ ((j_{1},k_{1}),(j_{2},k_{2}),\alpha) \not = \tilde{\boldsymbol{\eta}}_{3}  }} \sum_{p_{1} + p_{2} \ge 0} \vert d_{((p_{1},j_{1},k_{1}),(p_{2},j_{2},k_{2}),\alpha)} \vert^2  .
\end{align*}
Notice that by definition of the sets $\Upsilon$ we have 
\begin{align*}
 \Upsilon(\tilde{\boldsymbol{\lambda}}) \cup \{\tilde{\boldsymbol{\eta}}_1, \tilde{\boldsymbol{\eta}}_2, \tilde{\boldsymbol{\eta}}_3   \} = \Upsilon(\tilde{\boldsymbol{\eta}}_1) \cup \Upsilon(\tilde{\boldsymbol{\eta}}_2) \cup \Upsilon(\tilde{\boldsymbol{\eta}}_3) ,
\end{align*}
see the explanations below \eqref{eq:Upsilon}. Hence, we get
\begin{align*}
e_{0}(  \tilde{\boldsymbol{\lambda}}   )	&  \geq  \sum_{((j_{1}, k_{1}),(j_{2}, k_{2}),\alpha) \in \Upsilon(  \tilde{\boldsymbol{\eta}}_{1}  )} \sum_{p_{1} + p_{2} > 0} \vert d_{((p_{1},j_{1},k_{1}),(p_{2},j_{2},k_{2}),\alpha)} \vert^2 \\
	& \qquad \qquad   + \sum_{\substack{((j_{1},k_{1}),(j_{2},k_{2}),\alpha)  \in  \mathcal{J}_{\tilde{\boldsymbol{\eta}}_{1}} \\ ((j_{1},k_{1}),(j_{2},k_{2}),\alpha) \not = \tilde{\boldsymbol{\eta}}_{1}  }} \sum_{p_{1} + p_{2} \ge 0} \vert d_{((p_{1},j_{1},k_{1}),(p_{2},j_{2},k_{2}),\alpha)} \vert^2 \\
	& \qquad \qquad + \sum_{((j_{1}, k_{1}),(j_{2}, k_{2}),\alpha) \in \Upsilon(  \tilde{\boldsymbol{\eta}}_{2}  )} \sum_{p_{1} + p_{2} > 0} \vert d_{((p_{1},j_{1},k_{1}),(p_{2},j_{2},k_{2}),\alpha)} \vert^2 \\
	& \qquad \qquad   + \sum_{\substack{((j_{1},k_{1}),(j_{2},k_{2}),\alpha)  \in  \mathcal{J}_{\tilde{\boldsymbol{\eta}}_{2}} \\ ((j_{1},k_{1}),(j_{2},k_{2}),\alpha) \not = \tilde{\boldsymbol{\eta}}_{2}  }} \sum_{p_{1} + p_{2} \ge 0} \vert d_{((p_{1},j_{1},k_{1}),(p_{2},j_{2},k_{2}),\alpha)} \vert^2 \\ 
	& \qquad \qquad + \sum_{((j_{1}, k_{1}),(j_{2}, k_{2}),\alpha) \in \Upsilon(  \tilde{\boldsymbol{\eta}}_{3}  )} \sum_{p_{1} + p_{2} > 0} \vert d_{((p_{1},j_{1},k_{1}),(p_{2},j_{2},k_{2}),\alpha)} \vert^2 \\
	& \qquad \qquad   + \sum_{\substack{((j_{1},k_{1}),(j_{2},k_{2}),\alpha)  \in  \mathcal{J}_{\tilde{\boldsymbol{\eta}}_{3}} \\ ((j_{1},k_{1}),(j_{2},k_{2}),\alpha) \not = \tilde{\boldsymbol{\eta}}_{3}  }} \sum_{p_{1} + p_{2} \ge 0} \vert d_{((p_{1},j_{1},k_{1}),(p_{2},j_{2},k_{2}),\alpha)} \vert^2 \\
	& = e_{0}(  \tilde{\boldsymbol{\eta}}_{1}   ) +  e_{0}(  \tilde{\boldsymbol{\eta}}_{2}   ) + e_{0}(  \tilde{\boldsymbol{\eta}}_{3}   )  .
\end{align*}	
Here in the last step again we used Definition \ref{def_locerr_L2_new}. Consequently, \eqref{eq:1.2} is fulfilled.	

\textit{Step 2.}
Now we verify property \eqref{eq:1.3}. For that purpose let $ \tilde{\boldsymbol{\lambda}}  \in \tilde{\boldsymbol{\Lambda}}_{0} $ and $p_{\max} \in \mathbb{N}_{0}$. Then Definition \ref{def_locerr_L2_new} yields 
\begin{align*}
	e_{p_{\max}}(  \tilde{\boldsymbol{\lambda}}   ) & = \sum_{((j_{1}, k_{1}),(j_{2}, k_{2}),\alpha) \in \Upsilon(  \tilde{\boldsymbol{\lambda}}  )} \sum_{p_{1} + p_{2} > p_{\max}} \vert d_{((p_{1},j_{1},k_{1}),(p_{2},j_{2},k_{2}),\alpha)} \vert^2 \\
	& \qquad \qquad  + \sum_{\substack{((j_{1},k_{1}),(j_{2},k_{2}),\alpha)  \in  \mathcal{J}_{\tilde{\boldsymbol{\lambda}}} \\ ((j_{1},k_{1}),(j_{2},k_{2}),\alpha) \not = \tilde{\boldsymbol{\lambda}}  }} \sum_{p_{1} + p_{2} \ge 0} \vert d_{((p_{1},j_{1},k_{1}),(p_{2},j_{2},k_{2}),\alpha)} \vert^2 \\
	& \geq \sum_{((j_{1}, k_{1}),(j_{2}, k_{2}),\alpha) \in \Upsilon(  \tilde{\boldsymbol{\lambda}}  )} \sum_{p_{1} + p_{2} > p_{\max}+1} \vert d_{((p_{1},j_{1},k_{1}),(p_{2},j_{2},k_{2}),\alpha)} \vert^2 \\
	& \qquad \qquad  + \sum_{\substack{((j_{1},k_{1}),(j_{2},k_{2}),\alpha)  \in  \mathcal{J}_{\tilde{\boldsymbol{\lambda}}} \\ ((j_{1},k_{1}),(j_{2},k_{2}),\alpha) \not = \tilde{\boldsymbol{\lambda}}  }} \sum_{p_{1} + p_{2} \ge 0} \vert d_{((p_{1},j_{1},k_{1}),(p_{2},j_{2},k_{2}),\alpha)} \vert^2 \\
	& = e_{p_{\max}+1}(  \tilde{\boldsymbol{\lambda}}   ) .
	\end{align*}
Hence, \eqref{eq:1.3} is fulfilled. Moreover, as already mentioned none of the arguments we used above depends on the question whether $  \tilde{\boldsymbol{\lambda}}  $ has three or two children. Consequently, the case that there are only two children can be handled with similar methods and the proof is complete.
\end{proof}

To continue let $ f \in L_{2}((0,1)^2)   $ in the form \eqref{eq_canon_dual_frame_2} be given. We run the algorithm {\bf BIVARIATE\textunderscore NEARBEST\textunderscore TREE} with a subsequent trimming using the routine {\bf BIVARIATE\textunderscore TRIM} in order to find an approximation in terms of bivariate tensor quarklets which has the form \eqref{eq_func_quarktree_approx_mod}. Below we show that then the global error describes the quality of the resulting approximation.

\begin{Lemma}\label{lem_L2_globerr}
Let $ m \geq 2   $ and $ \tilde{m} \in \mathbb{N} $  with $  \tilde{m} \geq m  $ and $ m + \tilde{m} \in 2 \mathbb{N}   $. Let $ f \in L_2((0,1)^2)  $ be given in the form \eqref{eq_canon_dual_frame_2}. Let the local errors $ e_{p_{\max}}(\tilde{\boldsymbol{\lambda}}^{\circ})   $ be defined as in Definition \ref{def_locerr_L2_new}.  For $ N \in \mathbb{N} $ by $T_{N} = (\mathcal{T}_N, P_{\max})$ we denote the bivariate quarklet tree resulting out of the algorithm {\bf BIVARIATE\textunderscore NEARBEST\textunderscore TREE} with a subsequent trimming. Let $\mathcal{R} =  ((0,0), (0,0), 0 )  $ be the root of $  \mathcal{T}_N $. The corresponding quarklet tree approximation $ f_{T_N} $ of $f$ is given by
	\begin{equation}\label{approx_fN_L2}
	f_{T_N} = \sum_{\tilde{\boldsymbol{\lambda}} \in \tilde{T}_{N} \subset \boldsymbol{\tilde{\nabla}}  }  c_{\tilde{\boldsymbol{\lambda}}} w_{\tilde{\boldsymbol{\lambda}}}^{-1}  \boldsymbol{\psi_{\tilde{\lambda}}}  ,
	\end{equation} 
whereby $ \tilde{T}_{N}  $ is defined by \eqref{eq_func_quarktree_approx_mod}.	Then there exists a constant $ C > 0 $ independent of $ f $ and $ N $, such that for the global error we observe
	\begin{align}\label{eq:reliable}
	\Vert f - f_{T_N} \vert L_{2}((0,1)^2) \Vert^2 \leq C \mathcal{E}(T_{N}) .
	\end{align}
\end{Lemma} 

\begin{proof}
For the proof let $ f \in L_2((0,1)^2)  $ in the form \eqref{eq_canon_dual_frame_2} be given. For $ N \in \mathbb{N} $ by $T_{N} = (\mathcal{T}_N, P_{\max})$ we denote the bivariate quarklet tree produced by {\bf BIVARIATE\textunderscore NEARBEST\textunderscore TREE}. In a first step we apply the definition of the global error as given in \eqref{eq_glob_err}. When we combine it with Definition \ref{def_locerr_L2_new}, we find 
	\begin{align*}
	\mathcal{E}(T_{N})  &= \sum_{\tilde{\boldsymbol{\lambda}}^{\circ} \in \mathcal{V}(\mathcal{T}_N)} e_{p_{\max}(\tilde{\boldsymbol{\lambda}}^{\circ})}(\tilde{\boldsymbol{\lambda}}^{\circ})  \\
	&= \sum_{\tilde{\boldsymbol{\lambda}}^{\circ} \in \mathcal{V}(\mathcal{T}_{N})} \Big( \sum_{((j_{1}, k_{1}),(j_{2}, k_{2}),\alpha) \in \Upsilon(  \tilde{\boldsymbol{\lambda}}^{\circ}  )} \sum_{p_{1} + p_{2} > p_{\max}(\tilde{\boldsymbol{\lambda}}^{\circ})} \vert d_{((p_{1},j_{1},k_{1}),(p_{2},j_{2},k_{2}), \alpha)} \vert^2 \\
	& \qquad \qquad \qquad \qquad  + \sum_{\substack{((j_{1},k_{1}),(j_{2},k_{2}),\alpha)  \in  \mathcal{J}_{\tilde{\boldsymbol{\lambda}}^{\circ}} \\ ((j_{1},k_{1}),(j_{2},k_{2}),\alpha) \not = \tilde{\boldsymbol{\lambda}}^{\circ}  }} \sum_{p_{1} + p_{2} \ge 0} \vert d_{((p_{1},j_{1},k_{1}),(p_{2},j_{2},k_{2}), \alpha)} \vert^2 \Big )   .
	\end{align*} 
Recall, that by definition of the sets $\Upsilon(\tilde{\boldsymbol{\lambda}}^{\circ})$ we get $ \bigcup_{ \tilde{\boldsymbol{\lambda}}^{\circ} \in \mathcal{V}(\mathcal{T}_N)} \Upsilon(\tilde{\boldsymbol{\lambda}}^{\circ}) = \mathcal{T}_N$, see the explanations below \eqref{eq:Upsilon}. Consequently, we obtain 
\begin{align*}
	\mathcal{E}(T_{N}) &=   \sum_{((j_{1}, k_{1}),(j_{2}, k_{2}),\alpha) \in  \mathcal{T}_N   } \sum_{p_{1} + p_{2} > p_{\max}(\tilde{\boldsymbol{\lambda}}^{\circ})} \vert d_{((p_{1},j_{1},k_{1}),(p_{2},j_{2},k_{2}), \alpha)} \vert^2 \\
	& \qquad \qquad  + \sum_{\tilde{\boldsymbol{\lambda}}^{\circ} \in \mathcal{V}(\mathcal{T}_{N})}  \sum_{\substack{((j_{1},k_{1}),(j_{2},k_{2}),\alpha)  \in  \mathcal{J}_{\tilde{\boldsymbol{\lambda}}^{\circ}} \\ ((j_{1},k_{1}),(j_{2},k_{2}),\alpha) \not = \tilde{\boldsymbol{\lambda}}^{\circ}  }} \sum_{p_{1} + p_{2} \ge 0} \vert d_{((p_{1},j_{1},k_{1}),(p_{2},j_{2},k_{2}), \alpha)} \vert^2    .
	\end{align*} 	
In the expression above the first sum runs through all modified quarklet coefficients coming from the sequence $   \{  d_{\tilde{\boldsymbol{\lambda}}}  \}_{\tilde{\boldsymbol{\lambda}} \in  \tilde{\boldsymbol{\Lambda}}  } \in \ell_{2}( \tilde{\boldsymbol{\Lambda}}  )     $ whose corresponding wavelet indices are nodes of the tree $ \mathcal{T}_{N} $. But nevertheless these coefficients do not belong to the tree $ T_{N}  $ since the cumulated polynomial degrees of the corresponding bivariate tensor quarklets are too large. The second sum collects all modified quarklet coefficients with quarklet indices that do not belong to  $ T_{N}  $ since their corresponding wavelet indices are descendants of the leaves of $ \mathcal{T}_{N} $. Hence, using this in combination with $\mathcal{R} =  ((0,0), (0,0), 0 )  $, we can also write
\begin{align*}
	\mathcal{E}(T_{N}) =   \sum_{ \tilde{\boldsymbol{\lambda}} \in \tilde{\boldsymbol{\Lambda}}  } \vert  d_{\tilde{\boldsymbol{\lambda}}} \vert^{2}    -   \sum_{ \tilde{\boldsymbol{\lambda}} \in T_{N} \subset \tilde{\boldsymbol{\Lambda}} }   \vert d_{\tilde{\boldsymbol{\lambda}}} \vert^{2}   .
\end{align*}
Recall, that for any $  \tilde{\boldsymbol{\lambda}} \in \tilde{\boldsymbol{\Lambda}}  $ with $  j_{1} < j_{0}  $ or/and $ j_{2} < j_{0}    $ we have $  \vert d_{\tilde{\boldsymbol{\lambda}}} \vert^{2} = 0   $, see \eqref{eq_def_d_lambda}. Consequently, this also can be expressed as 
\begin{equation}\label{eq_proof_d_c_rewrite11}
	\mathcal{E}(T_{N}) =   \sum_{ \tilde{\boldsymbol{\lambda}} \in  \boldsymbol{\tilde{\nabla}} } \vert  d_{\tilde{\boldsymbol{\lambda}}} \vert^{2}    -   \sum_{ \tilde{\boldsymbol{\lambda}} \in T_{N} \cap \boldsymbol{\tilde{\nabla}} }   \vert d_{\tilde{\boldsymbol{\lambda}}} \vert^{2}   .
\end{equation}
To continue we observe that for $ \tilde{\boldsymbol{\lambda}} \in  \boldsymbol{\tilde{\nabla}}    $ with $ j_{1} > j_{0}   $ and $  j_{2} > j_{0} $ we have $  \vert d_{\tilde{\boldsymbol{\lambda}}} \vert^{2} =  \vert c_{\tilde{\boldsymbol{\lambda}}} \vert^{2}    $. To rewrite the first sum in \eqref{eq_proof_d_c_rewrite11} investigate $  \tilde{\boldsymbol{\lambda}} = ( (p_{1},j_{0},k_{1}), (p_{2},j_{2},k_{2}), \alpha ) \in  \boldsymbol{\tilde{\nabla}}    $ with fixed $p_{1}, p_{2} \in \mathbb{N}_{0}$, $j_{2} > j_{0}$, $  k_{2} \in  \nabla_{j_{2}} $ and $  \alpha \in \{ 0 ,1 ,2 \}   $.  $k_{1}$ runs through $ \nabla_{j_{0}} =
  \{ 0, 1, \ldots , 2^{j_{0}} - 1 \}   $. For these $ \tilde{\boldsymbol{\lambda}}   $ we find 
\begin{align*}
& \sum_{ k_{1} \in \nabla_{j_{0}} } \vert  d_{( (p_{1},j_{0},k_{1}), (p_{2},j_{2},k_{2}), \alpha )} \vert^{2} \\
& \qquad \qquad = \sum_{ k_{1} \in \nabla_{j_{0}} } \vert  c_{( (p_{1},j_{0},k_{1}), (p_{2},j_{2},k_{2}), \alpha )} \vert^{2} +  \sum_{ k_{1} \in \nabla_{j_{0}} } \sum_{\ell_{1} \in \square_{j_{0}, k_{1} }} |c_{( (p_{1},j_{0}-1,\ell_{1}), (p_{2},j_{2},k_{2}), \alpha )} |^{2} \\
& \qquad \qquad = \sum_{ k_{1} \in \nabla_{j_{0}} } \vert  c_{( (p_{1},j_{0},k_{1}), (p_{2},j_{2},k_{2}), \alpha )} \vert^{2} +  \sum_{ \tilde{k}_{1} \in \nabla_{j_{0}-1} }  |c_{( (p_{1},j_{0}-1,\tilde{k}_{1}), (p_{2},j_{2},k_{2}), \alpha )} |^{2} .
\end{align*}
Here in the last step we used \eqref{eq_def_l_k_opt2} and \eqref{index_square_dissdor}. This argument can be repeated with any possible  combination of $p_{1}, p_{2} \in \mathbb{N}_{0}$, $j_{2} > j_{0}$, $  k_{2} \in  \nabla_{j_{2}} $ and $  \alpha \in \{ 0 ,1 ,2 \}   $. Moreover, for given  $  \tilde{\boldsymbol{\lambda}} = ( (p_{1},j_{1},k_{1}), (p_{2},j_{0},k_{2}), \alpha ) \in  \boldsymbol{\tilde{\nabla}}    $ with fixed $p_{1}, p_{2} \in \mathbb{N}_{0}$, $j_{1} > j_{0}$, $  k_{1} \in  \nabla_{j_{1}} $, $  \alpha \in \{ 0 ,1 ,2 \}   $ and  $k_{2}$ running through $ \nabla_{j_{0}} =
  \{ 0, 1, \ldots , 2^{j_{0}} - 1 \}   $, a similar computation yields
\begin{align*}
& \sum_{ k_{2} \in \nabla_{j_{0}} } \vert  d_{( (p_{1},j_{1},k_{1}), (p_{2},j_{0},k_{2}), \alpha )} \vert^{2} \\
& \qquad \qquad = \sum_{ k_{2} \in \nabla_{j_{0}} } \vert  c_{( (p_{1},j_{1},k_{1}), (p_{2},j_{0},k_{2}), \alpha )} \vert^{2} +  \sum_{ \tilde{k}_{2} \in \nabla_{j_{0}-1} }  |c_{( (p_{1},j_{1},k_{1}), (p_{2},j_{0}-1,\tilde{k}_{2}), \alpha )} |^{2} .
\end{align*} 
Using this in combination with $ \vert  d_{( (p_{1},j_{1},k_{1}), (p_{2},j_{2},k_{2}), \alpha )} \vert^{2} = 0   $  for $j_{1} = j_{0} - 1$ or/and $j_{2} = j_{0}-1$ we get
\begin{equation}\label{eq_proof_d_c_rewrite22}
\sum_{ \tilde{\boldsymbol{\lambda}} \in  \boldsymbol{\tilde{\nabla}} } \vert  d_{\tilde{\boldsymbol{\lambda}}} \vert^{2}   =   \sum_{ \tilde{\boldsymbol{\lambda}} \in  \boldsymbol{\tilde{\nabla}} } \vert  c_{\tilde{\boldsymbol{\lambda}}} \vert^{2}  . 
\end{equation} 
To deal with the second sum in \eqref{eq_proof_d_c_rewrite11} let $ \tilde{\boldsymbol{\lambda}} = ( (p_{1},j_{1},k_{1}), (p_{2},j_{2},k_{2}), \alpha ) \in T_{N} \cap \boldsymbol{\tilde{\nabla}}    $. Then due to \eqref{eq_def_d_lambda} for $ j_{1} > j_{0}   $ and $  j_{2} > j_{0} $ we have $  \vert d_{\tilde{\boldsymbol{\lambda}}} \vert^{2} =  \vert c_{\tilde{\boldsymbol{\lambda}}} \vert^{2}    $. As before for $j_{1} = j_{0} -1 $ and/or $j_{2} = j_{0} - 1$ we have $ \vert d_{\tilde{\boldsymbol{\lambda}}} \vert^{2} = 0   $. It remains to deal with the case $j_{1} = j_{0}$ and $j_{2} > j_{0}$. Here we observe
\begin{align*}
\vert d_{\tilde{\boldsymbol{\lambda}}} \vert^{2} =  \vert c_{\tilde{\boldsymbol{\lambda}}} \vert^{2} + \sum_{\ell_{1} \in \square_{j_{0}, k_{1} }} \vert c_{( (p_{1},j_{0}-1,\ell_{1}), (p_{2},j_{2},k_{2}), \alpha )} \vert^{2} .
\end{align*}
Let us recall $  \square_{j_{0}, k_{1}} := \{ k \in \nabla_{j_0-1} : \ell_k = k_{1}     \}    $. On the other hand the transformed tree $ \tilde{T}_{N} \subset \boldsymbol{\tilde{\nabla}}      $ contains the nodes $  ( (p_{1},j_{0},k_{1}), (p_{2},j_{2},k_{2}), \alpha ) \in T_{N}  $ and in addition also the nodes
\begin{small}
\begin{align*}
\{ ((p_{1},j_{0}-1,k), (p_{2},j_{2},k_{2}), \alpha ) \in        \boldsymbol{\tilde{\nabla}} : \tilde{\boldsymbol{\lambda}} = ((0,j_{0},\ell_{k}), (0,j_{2},k_{2}), \alpha ) \in T_{N}, p_{1} + p_{2} \leq p_{\max}(\tilde{\boldsymbol{\lambda}}^{\circ})  \} 
\end{align*}
\end{small}
with $  \ell_{k} = k_{1}   $. A similar observation can be made for the case $ j_{1} > j_{0}   $ and $  j_{2} = j_{0}   $. Consequently, we can write
\begin{equation}\label{eq_proof_d_c_rewrite33}
\sum_{ \tilde{\boldsymbol{\lambda}} \in T_{N} \cap \boldsymbol{\tilde{\nabla}} }   \vert d_{\tilde{\boldsymbol{\lambda}}} \vert^{2}  = \sum_{ \tilde{\boldsymbol{\lambda}} \in \tilde{T}_{N} \subset \boldsymbol{\tilde{\nabla}} }   \vert c_{\tilde{\boldsymbol{\lambda}}} \vert^{2}  .
\end{equation}
Now a combination of \eqref{eq_proof_d_c_rewrite11} with \eqref{eq_proof_d_c_rewrite22} and \eqref{eq_proof_d_c_rewrite33} yields 
\begin{equation}\label{eq_proof_d_c_rewrite44}
	\mathcal{E}(T_{N}) =    \sum_{ \tilde{\boldsymbol{\lambda}} \in  \boldsymbol{\tilde{\nabla}} } \vert  c_{\tilde{\boldsymbol{\lambda}}} \vert^{2}    -  \sum_{ \tilde{\boldsymbol{\lambda}} \in \tilde{T}_{N} \subset \boldsymbol{\tilde{\nabla}} }   \vert c_{\tilde{\boldsymbol{\lambda}}} \vert^{2}    .
\end{equation}
To continue we apply Theorem \ref{thm_biv_tensor_frame_L2}. It shows that the family 
\begin{align*}
\boldsymbol{\Psi}_{L_{2}((0,1)^2)} =  \Big \{  w_{\boldsymbol{\lambda}}^{-1} \boldsymbol{\psi_{\lambda}} : \boldsymbol{\lambda}  \in  \boldsymbol{\nabla} := \nabla \times \nabla \Big \} 
\end{align*}
is a quarkonial tensor frame for $ L_{2}((0,1)^2)  $. The weights $ w_{\boldsymbol{\lambda}}   $ are given by \eqref{tensor_frame_weights1} with $   \delta > 1$. Due to \eqref{eq_extend_frame_tilde1} and \eqref{eq_extend_frame_tilde2} also the bivariate tensor quarklets corresponding to the index set $ \boldsymbol{\tilde{\nabla}}   $ are a frame for $ L_{2}((0,1)^2)  $. Next we use the lower estimate given in Proposition 2.4 in \cite{DaHoRaVo} to find
\begin{equation}\label{eq_proof_d_c_rewrite55}
\mathcal{E}(T_{N}) \gtrsim \Big \Vert   \sum_{ \tilde{\boldsymbol{\lambda}} \in  \boldsymbol{\tilde{\nabla}}    } c_{\tilde{\boldsymbol{\lambda}}} w_{\tilde{\boldsymbol{\lambda}}}^{-1} \boldsymbol{\psi_{\tilde{\lambda}}}  -    \sum_{ \tilde{\boldsymbol{\lambda}} \in \tilde{T}_{N} \subset \boldsymbol{\tilde{\nabla}} }  c_{\tilde{\boldsymbol{\lambda}}} w_{\tilde{\boldsymbol{\lambda}}}^{-1} \boldsymbol{\psi_{\tilde{\lambda}}} \Big \vert  L_{2}((0,1)^2) \Big \Vert^2 .
\end{equation}
Recall, that $ f \in L_2((0,1)^2)  $ has the form \eqref{eq_canon_dual_frame_2}. Moreover, the quarklet tree approximation $ f_{T_N} $ of $f$ is given by \eqref{approx_fN_L2}. Hence, \eqref{eq_proof_d_c_rewrite55} becomes 
\begin{align*}
\mathcal{E}(T_{N}) \gtrsim \Big \Vert   f  -  f_{T_N}  \Big \vert  L_{2}((0,1)^2) \Big \Vert^2 .
\end{align*}
The proof is complete.
\end{proof} 

Now we show that for $  f \in L_2((0,1)^2)  $ given by \eqref{eq_canon_dual_frame_2} the quarklet tree approximation $ f_{T_N} $ produced by the algorithm {\bf BIVARIATE\textunderscore NEARBEST\textunderscore TREE} is near-best. More precisely, there is the following result. 

\begin{Theorem}\label{lem_L2_bestappr}
Let $ m \geq 2   $ and $ \tilde{m} \in \mathbb{N} $  with $  \tilde{m} \geq m  $ and $ m + \tilde{m} \in 2 \mathbb{N}   $. Let $ f \in L_2((0,1)^2)  $ be given in the form \eqref{eq_canon_dual_frame_2}.  Let the local errors $ e_{p_{\max}}(\tilde{\boldsymbol{\lambda}}^{\circ})   $ be defined as in Definition \ref{def_locerr_L2_new}. For $ N \in \mathbb{N} $ by $T_N = (\mathcal{T}_{N},P_{\max})$ we denote the bivariate quarklet tree resulting out of the algorithm {\bf BIVARIATE\textunderscore NEARBEST\textunderscore TREE} with a subsequent trimming. Let $\mathcal{R} =  ((0,0), (0,0), 0 )  $ be the root of $  \mathcal{T}_N $. As already seen in Lemma \ref{lemma:T_N_cardinality} the cardinality of $  T_{N} $ fulfills   $\#T_{N} \lesssim N^{3} $. The corresponding quarklet tree approximation $ f_{T_N} $ of $f$ is given by
\begin{align*}
f_{T_N} = \sum_{\tilde{\boldsymbol{\lambda}} \in \tilde{T}_{N} \subset \boldsymbol{\tilde{\nabla}}  }  c_{\tilde{\boldsymbol{\lambda}}} w_{\tilde{\boldsymbol{\lambda}}}^{-1}  \boldsymbol{\psi_{\tilde{\lambda}}} ,
\end{align*} 
whereby $ \tilde{T}_{N}  $ is defined by \eqref{eq_func_quarktree_approx_mod}. Let $ n \le N  $. Then there exists a constant $ C > 0 $ independent of $ f $, $ N $ and $ n $, such that 
	\begin{align*}
	\Vert f - f_{T_N} \vert L_{2}((0,1)^2) \Vert^2 \leq C  \frac{3N+1}{ N    -  \frac{2}{3}  n + \frac{1}{2}}  \sigma_n  .
	\end{align*}
\end{Theorem}

\begin{proof}
For the proof at first we observe that all conditions stated in Lemma \ref{lem_L2_globerr} are fulfilled. Consequently, it can be applied, and we get
	\begin{align*}
	\Vert f - f_{T_N} \vert L_{2}((0,1)^2)  \Vert^2 \leq C \mathcal{E}(T_{N}) 
	\end{align*} 
with a constant $C > 0$ independent of $f$ and $N$. To continue we want to use Theorem \ref{theorem:1}. To this end recall that the local errors $ e_{p_{\max}}(\tilde{\boldsymbol{\lambda}})  $ given in Definition \ref{def_locerr_L2_new} fulfill the properties \eqref{eq:1.2} and \eqref{eq:1.3}. This observation already has been made in Lemma \ref{lem_loc_err_prop1_new}. Hence, an application of Theorem \ref{theorem:1} yields 
	\begin{align*}
	\Vert f - f_{T_N}  \vert L_{2}((0,1)^2) \Vert^2 \leq C  \frac{3N+1}{ N    -  \frac{2}{3}  n + \frac{1}{2}}  \sigma_n  .
	\end{align*}
Here $  \sigma_n  $ is the error of the best bivariate quarklet tree approximation of cardinality $n \in \mathbb{N}$ as defined in Definition \ref{def_bestappr_err}. The proof is complete. 
\end{proof}

\begin{Remark}
It seems to be possible to prove counterparts of Theorem \ref{lem_L2_bestappr} for Sobolev functions $  f \in H^{s}_{2}((0,1)^2)   $ with $ s>0 $, where the approximation error is measured in terms of a Sobolev norm. For the univariate setting such results can be found in \cite{VoDiss}, see Chapter 4.5.3. 
\end{Remark}

\begin{Remark}
When we look at Definition \ref{def_locerr_L2_new} and Theorem \ref{lem_L2_bestappr} it becomes clear that the machinery presented in this section only can be used if $ f \in L_2((0,1)^2)  $ is given in the form \eqref{eq_canon_dual_frame_2}. In some applications this assumption might be fulfilled automatically, for example if $f$ is the result of prior computations and represents the current approximation of an unknown solution in terms of elements of a bivariate quarklet frame. However, sometimes a representation of the form \eqref{eq_canon_dual_frame_2} will not be given. Then we have to find a sequence $ \{  c_{\tilde{\boldsymbol{\lambda}}}  \}_{\tilde{\boldsymbol{\lambda}} \in  \boldsymbol{\tilde{\nabla}} } \in \ell_{2}( \boldsymbol{\tilde{\nabla}} ) $. How this can be done highly depends on the problem we have in mind. If $f$ is explicitly known, we can find a representation \eqref{eq_canon_dual_frame_2} by solving a matrix-vector equation similar to (35) in \cite{DaHoRaVo}. Remark 4.8 in \cite{DaHoRaVo} provides much more information concerning this topic. The strategies discussed there refer to the univariate case, but can be modified in order to work in the bivariate setting. On the other hand, if $f$ is the unknown solution of a linear elliptic variational problem, for instance, then approximations for $f$  in the form \eqref{eq_canon_dual_frame_2} can be obtained by using the damped Richardson iteration or variations thereof. Much more details concerning this issue can be found in \cite{DaFKRaa}, see Section 4.1. 
\end{Remark}

\subsection{A Test Case achieving inverse-exponential Rates of Convergence}\label{Sec_aniso_example}

To reinforce the results of this paper below we investigate a test function which can be approximated via adaptive bivariate quarklet tree approximation in a very efficient way. For $  x = ( x_{1}, x_{2} ) \in (0,1)^2   $ we deal with the function 
\begin{align*}
f_{\alpha}(x_{1},x_{2}) = x^{\alpha}_{1} , \qquad \qquad \alpha > \frac{1}{2} .
\end{align*}
In \cite{DaRaaS} it has been observed that anisotropic singularities of the form $  f_{\alpha}  $ can be approximated by anisotropic tensor product quarklets very well achieving inverse-exponential rates of convergence. In what follows we will see that the approximating function constructed in \cite{DaRaaS} fits into the setting of adaptive bivariate quarklet tree approximation presented in the current paper. For that purpose we briefly describe the construction provided in \cite{DaRaaS} and explain the connections to bivariate quarklet tree approximation. To start the approximation procedure we use the root $\mathcal{R} =  ((0,0), (0,0), 0 )  $. Then we carry out $ L \in \mathbb{N}   $ refinement steps in direction $i = 1$ using the refinement strategy (LSR.a.1) for the leftmost reference rectangle. For each $\ell \in \{   1, 2, \ldots , L \}$ we obtain the children
\begin{equation}\label{example_eq_111}
\{ ( ( \ell , 0), (0,0), 0 ) ,  ( ( \ell ,  1), (0,0), 0 )  ,  ( (\ell - 1,0), (0,0), 2 ) \} .
\end{equation} 
This refers to the reference rectangles 
\begin{equation}\label{example_eq_222}
[ 0 , 2^{-\ell}  ) \times  [ 0 , 1 )  \qquad \mbox{and} \qquad [ 2^{-\ell}  , 2^{-\ell +1}  ) \times  [ 0 , 1 ) .
\end{equation}
Using a slightly different notation they also can be found in Theorem 4.9 in \cite{DaRaaS}. Notice that this refinement strategy goes along with condition (UPC). By Section \ref{sec_ref_rec_wavl} it is clear which bivariate tensor wavelets can be associated with the reference rectangles given in \eqref{example_eq_222}. To continue, for each refinement level $\ell \in \{   1, 2, \ldots , L \}$ we have to add some more nodes of the form
\begin{align*}
( ( \ell , k), (0,0), 0 ) , \qquad \qquad 1 < k < C(m,\ell), 
\end{align*}
whereby $ C(m, \ell ) \in \mathbb{N}   $ depends on $m$ and $\ell$ and is explicitly given in \cite{DaRaaS}, see Theorem 4.5 and Theorem 3.9. To obtain these wavelet indices only refinements of type (LSR.a.1) must be used. Now to each node of the bivariate wavelet tree we have to assign a maximal polynomial degree $ p_{\max}   $. Following Section 4 in \cite{DaRaaS} the polynomial degree depends on the refinement level $\ell \in \{   1, 2, \ldots , L \}$ and is given by
\begin{align*}
p_{\max}(\ell) := L - \ell + m - 3.
\end{align*}
Since for the creation of the underlying wavelet tree we only used refinements of the form (LSR.a.1) this choice goes along with Definition \ref{def:quarklet_tree}. Consequently, we obtain a collection of enhanced bivariate quarklet indices $ \tilde{\boldsymbol{\nabla}}_{L} \subset \tilde{\boldsymbol{\Lambda}}    $ depending on $L$ which is a bivariate quarklet tree. In Theorem 4.10 in \cite{DaRaaS} it is proved that the bivariate quarks and quarklets addressed by $ \tilde{\boldsymbol{\nabla}}_{L} $ can be applied to approximate $ f_{\alpha}  $ in a very efficient way. More precisely, it is shown that for $  N \in \mathbb{N}  $ with $ N \sim L^5   $  there exists a sequence $ \{  c_{\tilde{\boldsymbol{\lambda}}}  \}_{\tilde{\boldsymbol{\lambda}} \in  \tilde{\boldsymbol{\nabla}}_{L} } \in \ell_{2}( \tilde{\boldsymbol{\nabla}}_{L} ) $ such that 
\begin{equation}\label{example_eq_approximant11}
g = \sum_{ \tilde{\boldsymbol{\lambda}} \in \tilde{\boldsymbol{\nabla}}_{L} : \# \tilde{\boldsymbol{\nabla}}_{L} \leq N  }   c_{\tilde{\boldsymbol{\lambda}}} w_{\tilde{\boldsymbol{\lambda}}}^{-1}  \boldsymbol{\psi_{\tilde{\lambda}}}  
\end{equation}
fulfills
\begin{equation}\label{example_eq_approximant22}
\Big \Vert \frac{\partial}{\partial x_{1}} \Big [ f_{\alpha}(x_{1},x_{2}) - g(x_{1},x_{2}) \Big ]  \Big \vert L_{2}((0,1)^2) \Big \Vert^{2} \lesssim e^{- \min(2, 2 \alpha - 1) \ln(2) N^{\frac{1}{5}} } .
\end{equation}
Here \eqref{example_eq_approximant22} is formulated in terms of a Sobolev seminorm in order to directly apply Theorem 4.10 in \cite{DaRaaS}. In \cite{DaRaaS} it was assumed $ j_{0} = 0   $ which is possible as long as only inner quarks and quarklets are used. Of course there is no guarantee that the algorithm {\bf BIVARIATE\textunderscore NEARBEST\textunderscore TREE} exactly produces the index set $  \tilde{\boldsymbol{\nabla}}_{L}  $. However, since it is near-best in the sense of Theorem \ref{theorem:1}, it surely provides a bivariate quarklet tree whose approximation properties are comparable or even better.

\vspace{0,3 cm}

\textbf{Acknowledgment.} This work partly has been supported by Deutsche Forschungsgemeinschaft (DFG), grant HO 7444/1-1 with project number 528343051. The author would like to thank Stephan Dahlke and Dorian Vogel for several valuable discussions.

\end{document}